\documentclass[reqno,centertags]{amsart}
\usepackage{hyperref}
\usepackage{amsmath,amsthm,amscd,amssymb}
\usepackage{latexsym}

\newcommand{\bbC}{{\mathbb{C}}}
\newcommand{\bbD}{{\mathbb{D}}}

\newcommand{\bbR}{{\mathbb{R}}}

\newcommand{\bbZ}{{\mathbb{Z}}}

\newcommand{\calC}{{\mathcal{C}}}
\newcommand{\calD}{{\mathcal{D}}}
\newcommand{\calE}{{\mathcal{E}}}
\newcommand{\calF}{{\mathcal F}}
\newcommand{\calG}{{\mathcal G}}
\newcommand{\calH}{{\mathcal H}}

\newcommand{\calJ}{{\mathcal J}}

\newcommand{\calL}{{\mathcal L}}
\newcommand{\calM}{{\mathcal M}}
\newcommand{\calR}{{\mathcal R}}
\newcommand{\calS}{{\mathcal S}}
\newcommand{\calT}{{\mathcal T}}
\newcommand{\calU}{{\mathcal U}}
\newcommand{\calZ}{{\mathcal Z}}
\newcommand{\bdone}{{\boldsymbol{1}}}
\newcommand{\bdzero}{{\boldsymbol{0}}}
\newcommand{\bddot}{{\boldsymbol{\cdot}}}
\newcommand{\llangle}{\langle\!\langle}
\newcommand{\rrangle}{\rangle\!\rangle}

\newcommand{\lb}{\label}
\newcommand{\f}{\frac}

\newcommand{\ol}{\overline}
\newcommand{\ti}{\tilde  }
\newcommand{\wti}{\widetilde  }

\newcommand{\Lt}{\text{\rm{L}}}

\newcommand{\unif}{\text{\rm{unif}}}
\newcommand{\tr}{\text{\rm{Tr}}}

\newcommand{\dist}{\text{\rm{dist}}}
\newcommand{\loc}{\text{\rm{loc}}}

\newcommand{\rank}{\text{\rm{rank}}}

\newcommand{\ess}{\text{\rm{ess}}}
\newcommand{\ac}{\text{\rm{ac}}}

\newcommand{\s}{\text{\rm{s}}}

\newcommand{\supp}{\text{\rm{supp}}}

\newcommand{\bi}{\bibitem}

\newcommand{\beq}{\begin{equation}}
\newcommand{\eeq}{\end{equation}}
\newcommand{\ba}{\begin{align}}
\newcommand{\ea}{\end{align}}
\newcommand{\veps}{\varepsilon}
\newcommand{\mol}{\eta}

\newcounter{smalllist}
\newenvironment{SL}{\begin{list}{{\rm(\roman{smalllist})}}{%
\setlength{\topsep}{0mm}\setlength{\parsep}{0mm}\setlength{\itemsep}{0mm}%
\setlength{\labelwidth}{2em}\setlength{\leftmargin}{2em}\usecounter{smalllist}%
}}{\end{list}}
\newenvironment{SSL}{\begin{list}{{\rm(\roman{smalllist})}}{%
\setlength{\topsep}{1ex}\setlength{\parsep}{0mm}\setlength{\itemsep}{1ex}%
\setlength{\labelwidth}{2em}\setlength{\leftmargin}{2em}\usecounter{smalllist}%
}}{\end{list}}

\newcommand{\bigtimes}{\mathop{\mathchoice%
{\smash{\vcenter{\hbox{\LARGE$\times$}}}\vphantom{\prod}}%
{\smash{\vcenter{\hbox{\Large$\times$}}}\vphantom{\prod}}%
{\times}%
{\times}%
}\displaylimits}

\DeclareMathOperator{\Real}{Re}
\DeclareMathOperator{\Ima}{Im}

\DeclareMathOperator*{\wlim}{w-lim}

\DeclareMathOperator*{\Capac}{Cap}
\allowdisplaybreaks
\numberwithin{equation}{section}

\newtheorem{theorem}{Theorem}[section]

\newtheorem*{p2.1}{Proposition 2.1}
\newtheorem{proposition}[theorem]{Proposition}
\newtheorem{lemma}[theorem]{Lemma}
\newtheorem{corollary}[theorem]{Corollary}
\theoremstyle{definition}
\newtheorem*{definition}{Definition}
\newtheorem{example}[theorem]{Example}

\theoremstyle{remark}
\newtheorem*{remark}{Remark}
\newtheorem*{remarks}{Remarks}

\newcommand{\abs}[1]{\lvert#1\rvert}

\begin{document}
\title[Perturbations of OP With Periodic Coefficients]
{Perturbations of Orthogonal Polynomials With Periodic Recursion Coefficients}
\author[D.~Damanik, R.~Killip, and B.~Simon]
{David Damanik$^{1}$, Rowan Killip$^{2}$, and
Barry Simon$^3$}

\thanks{$^1$ Department of Mathematics, Rice University, Houston, TX 77005.
E-mail: damanik@rice.edu. Supported in part by NSF grants DMS-0500910 and DMS--0653720.}

\thanks{$^2$  Department of Mathematics, UCLA, Los Angeles, CA 90095.
E-mail: killip@math.ucla.edu. Supported in part by NSF grant DMS-0401277
and a Sloan Foundation Fellowship.}

\thanks{$^3$ Mathematics 253-37, California Institute of Technology, Pasadena, CA 91125.
E-mail: bsimon@caltech.edu. Supported in part by NSF grant DMS-0140592 and
U.S.--Israel Binational Science Foundation (BSF) Grant No.\ 2002068}

\date{December 5, 2008}

\begin{abstract} We extend the results of Denisov--Rakhmanov, Szeg\H{o}--Shohat--Nevai,
and Killip--Simon from asymptotically constant orthogonal polynomials on the real line
(OPRL) and unit circle (OPUC) to asymptotically periodic OPRL and OPUC. The key tool is
a characterization of the isospectral torus that is well adapted to the study of perturbations.
\end{abstract}

\maketitle

\section{Introduction} \lb{s1}

This is a paper about the spectral theory of orthogonal polynomials
on the real line (OPRL) and orthogonal polynomials on the unit circle
(OPUC), that is, the connection of the underlying (spectral) measure
and the recursion coefficients.

Specifically, given a probability measure, $d\mol$, on $\bbR$ with bounded
but infinite support, the orthonormal polynomials, $p_n(x)$, obey a recursion
relation
\begin{equation} \lb{1.1}
xp_n(x) =a_{n+1} p_{n+1}(x) + b_{n+1} p_n(x) + a_n p_{n-1}(x)
\end{equation}
where the Jacobi parameters $\{a_n,b_n\}_{n=1}^\infty$ obey $b_j\in\bbR$,
$a_j\in (0,\infty)$. As is well known (see, e.g., \cite[Section~1.3]{OPUC1}),
\eqref{1.1} sets up a one-one correspondence between uniformly bounded
$\{a_n,b_n\}_{n=1}^\infty$ and such measures, $d\mol$ (this is sometimes
called Favard's theorem).

Similarly, probability measures, $d\mu$, on $\partial\bbD$ which are nontrivial
(i.e., their support is not a finite set of points) are in one-one correspondence
with sequences $\{\alpha_n\}_{n=0}^\infty$ of Verblunsky coefficients in $\bbD
\equiv \{z : \abs{z}<1\}$ via the recursion relation of the orthonormal
polynomials $\varphi_n(z)$, namely,
\begin{equation} \lb{1.2}
z\varphi_n(z) = \rho_n \varphi_{n+1}(z) + \bar\alpha_n \varphi_n^*(z)
\end{equation}
where
\begin{equation} \lb{1.3}
\varphi_n^*(z) = z^n \, \ol{\varphi_n (1/\bar z)} \qquad
\rho_n = (1-\abs{\alpha_n}^2)^{1/2}
\end{equation}

Underlying the association of measures and recursion coefficients are matrix
representations. For OPRL, we take the matrix for multiplication by $x$
in the $\{p_n\}_{n=0}^\infty$ basis of $L^2 (\bbR,d\mol)$, which is the
tridiagonal Jacobi matrix
\begin{equation} \lb{1.3a}
J = \begin{pmatrix}
b_1 & a_1 & 0 & \vphantom{\ddots} \\
a_1 & b_2 & a_2 & \ddots  \\
0 & a_2 & \ddots & \ddots   \\
{} & \ddots & \ddots & \ddots
\end{pmatrix}
\end{equation}
For OPUC, one takes the matrix, $\calC$, for multiplication by $z$ in
the basis obtained by orthonormalizing $\{1,z,z^{-1},z^2, z^{-2}, \dots\}$
in $L^2 (\partial\bbD, d\mu)$. This CMV matrix  (see
\cite[Section~4.2]{OPUC1}) has the form
\begin{align}
\calC &=\calL\calM \lb{1.3b} \\
\calL &= \Theta(\alpha_0) \oplus\Theta(\alpha_2)\oplus\cdots \lb{1.3c} \\
\calM &= \bdone_{1\times 1} \oplus\Theta(\alpha_1)\oplus\Theta(\alpha_3)\oplus \cdots \lb{1.3d} \\
\Theta(\alpha) &= \begin{pmatrix}
\bar\alpha & \rho \\
\rho & -\alpha \end{pmatrix} \lb{1.3e}
\end{align}
where $\rho\equiv(1-|\alpha|^2)^{1/2}$. Note that $\calC$ is unitary, while $J$ is self-adjoint.

As a model for what we wish to prove, let us briefly survey some of the main results
relating to (slowly decaying) perturbations of the free case, that is,
$a_n\equiv 1$, $b_n\equiv 0$ for OPRL and $\alpha_n\equiv 0$ for OPUC.

\smallskip
\noindent (1) {\it Weyl's Theorem} \cite{Weyl,Blu,AK,GBk1,GSprep}.
If $a_n\to 1$, $b_n\to 0$, then $\sigma_\ess (d\mol)=[-2,2]$ and
if $\alpha_n\to 0$, then $\sigma_\ess (d\mu) = \partial\bbD$. Here
$\sigma_\ess (d\mol)$ is the (topological) support of the measure,
$d\mol$, with all isolated points removed.

\smallskip
\noindent (2) {\it Denisov--Rakhmanov Theorem} \cite{Rakh77,Rakh83,MNT85a,DenPAMS,NTppt}.
If $\sigma_\ess (d\mol)=\Sigma_\ac (d\mol) =[-2,2]$, then $a_n\to 1$ and $b_n\to 0$.
If $\sigma_\ess (d\mu)= \Sigma_\ac (d\mu) =\partial\bbD$, then $\alpha_n\to 0$.
Here $\Sigma_\ac (d\mol)$ is defined as follows: let $d\mol=d\mol_\ac + d\mol_\s$
with $d\mol_\s$ singular and $d\mol_\ac=f(x)dx$, then $\Sigma_\ac (d\mol)=\{x : f(x)\neq 0\}$ as a measure class,
that is, modulo sets of Lebesgue measure zero.

\smallskip
\noindent (3) {\it Szeg\H{o}'s Theorem} \cite{Sz15,Sz20,Sh,Nev79,KS}.
In the OPUC case, define
\begin{equation} \lb{1.4}
Z(d\mu) \equiv -\int \log \biggl( \f{d\mu_\ac}{d\theta}\biggr)\,
\f{d\theta}{2\pi}
\end{equation}
Then $Z(d\mu)<\infty$ if and only if
\begin{equation} \lb{1.5}
\sum_{j=0}^\infty \, \abs{\alpha_j(d\mu)}^2 <\infty
\end{equation}
In the OPRL case, define
\begin{equation} \lb{1.6}
Z(d\mol) \equiv -\int_{-2}^2 \log \biggl( 2\pi(4-E^2)^{1/2}\,
\f{d\mol_\ac}{dE}\biggr) \, \f{dE}{2\pi (4-E^2)^{1/2}}
\end{equation}
Then, if we assume $\supp (d\mol)\subset [-2,2]$, we have
\begin{equation} \lb{1.7}
Z(d\rho) <\infty \iff \limsup \sum_{j=1}^N \log (a_j) >-\infty
\end{equation}
and if that holds, then
\begin{equation} \lb{1.8}
\sum_{n=1}^\infty (a_n-1)^2 + b_n^2 <\infty
\end{equation}
and
\begin{equation} \lb{1.9}
\sum (a_n-1) \quad\text{and}\quad \sum b_n
\end{equation}
are conditionally convergent to finite numbers.

\smallskip
\noindent (4) {\it Killip--Simon Theorem} \cite{KS}. For OPRL, define
\begin{equation} \lb{1.10}
Q(d\mol) =-\int_{-2}^2 \log\biggl( \pi (4-E^2)^{-1/2}\,
\f{d\mol_\ac}{dE}\biggr) \, \f{(4-E^2)^{1/2}\, dE}{\pi}
\end{equation}
and let $\{E_j\}$ be the point masses of $d\mol$ (eigenvalues of $J$) outside $[-2,2]$. Then
\eqref{1.8} holds if and only if $\sigma_\ess (d\mol) =[-2,2]$, $Q(d\mol)<\infty$, and $\sum_j (\abs{E_j}-2)^{3/2} <\infty$.

\smallskip
\noindent (5) {\it Nevai's Conjecture} \cite{Nev92,KS}. For OPRL, if $\sum_{n=1}^\infty
\abs{a_n-1} +\abs{b_n}<\infty$, then $Z(d\rho)<\infty$ ($Z$ given by \eqref{1.6}).

\smallskip
The five results listed above capture different aspects of the philosophy that the measure
is close to the free case if and only if the coefficients are asymptotic to the free ones.
In this paper, we study extensions of all these results to perturbations of a
periodic sequence of Jacobi or Verblunsky coefficients, that is,
\begin{equation} \lb{1.11}
a_{n+p}^{(0)} =a_n^{(0)} \qquad b_{n+p}^{(0)} = b_n^{(0)} \qquad n\geq 1
\end{equation}
or
\begin{equation} \lb{1.12}
\alpha_{n+p}^{(0)} =\alpha_n^{(0)} \qquad n\geq 0
\end{equation}
and some fixed $p\geq 1$. Note that $p=1$ is the perturbation of the free case considered above.
For simplicity in the OPUC case, we will normally suppose $p$ is even---indeed, the shape of a CMV
matrix repeats itself only after shifting by two rows/columns.  As explained in Section~\ref{s13},
the situation when $p$ is odd can be reduced to this using sieving. For OPRL, $p$ is arbitrary.

The philosophy described above becomes more subtle when we move to the periodic setting;
rather than having a single `free operator' we have a manifold of them (the isospectral torus).
Nevertheless---and this is the main thrust of the paper---spectral measures that are
close to those of the isospectral torus correspond to coefficients that approach the
isospectral torus.   One of the key obstructions here is that a sequence of coefficients
may approach the isospectral torus without converging to any particular point therein.

In order to make these heuristics precise, we need to make a few definitions.  To keep the presentation
as coherent as possible, we will focus our attention on the OPRL/Jacobi case for the remainder of the introduction.

To any pair of $p$-periodic sequences, $\{a_n^{(0)},b_n^{(0)}\}_{n\in\bbZ}$, we can associate a
two-sided Jacobi matrix $J_0$.  Two such pairs of sequences are termed \emph{isospectral} if the corresponding
Jacobi matrices have the same spectrum.  We write $\calT_{J_0}$ for the set of $p$-periodic sequences that
are isospectral to $J_0$.   Topologically, this is a torus as explained in Subsection~\ref{ssIT} below.

Given two bounded sequences $\{a_n,b_n\}_{n=1}^\infty$ and $\{a'_n,b'_n\}_{n=1}^\infty$, we define
\begin{equation} \lb{1.43}
d_m ((a,b), (a', b')) = \sum_{k=0}^\infty e^{-k} [\abs{a_{m+k} -a'_{m+k}} + \abs{b_{m+k} - b'_{m+k}}]
\end{equation}
which is a metric for the product topology on $\bigtimes_m^\infty \bigl((0,R]\times [-R,R]\bigr)$.  The
OPUC analog is
\begin{equation} \lb{1.44}
d_m ((\alpha), (\alpha')) = \sum_{k=0}^\infty e^{-k} \abs{\alpha_{m+k} -\alpha'_{m+k}}
\end{equation}
a metric for $\bigtimes_m^\infty \bbD$.  The distance from a point to a set is defined in the usual way:
\begin{equation} \lb{1.45}
d_m ((a,b),\calT) = \inf \{d_m ((a,b), (a',b')) : (a',b')\in\calT\}
\end{equation}
and similarly in the OPUC case.

We begin with the periodic analog of Weyl's Theorem.

\begin{theorem}[Last--Simon \cite{LS_jdam}] \label{T1.1} Let $J_0$ be a two-sided periodic Jacobi
matrix and $J$ a one-sided Jacobi matrix with Jacobi parameters $\{a_m, b_m\}_{m=1}^\infty$. If
\begin{equation} \label{1.46}
d_m ((a,b),\calT_{J_0})\to 0
\end{equation}
then\footnote{Recall that $\sigma_\ess (J)$ is obtained from the
spectrum of the Jacobi matrix $J$ by removing all isolated
points.}
\begin{equation} \lb{1.47}
\sigma_\ess (J) = \sigma(J_0)
\end{equation}
\end{theorem}

As indicated, this result first appeared in \cite{LS_jdam}. It is derived from a theorem that
had earlier been proven with different methods by others \cite{GI1,M1,Rab5}.   The
inclusion $\sigma_\ess (J) \supset \sigma(J_0)$ follows easily using trial vectors;
the reverse seems to be more sophisticated.  In Section~\ref{s7}, we prove this using
the methods of this paper.  The OPUC version appears here as Theorem~\ref{T13.1};
it was also proved in \cite{LS_jdam}.

Note that \eqref{1.46} does not imply that there is a sequence $\{(a'_n,b'_n)\}\in \calT_{J_0}$ such that
$$
d_m ((a,b),(a',b'))\to 0
$$
It is much weaker.  Equality of essential spectra under this stronger hypothesis follows
immediately from Weyl's original theorem on compact perturbations.

Our first major new result is an analog of the Denisov--Rakhmanov Theorem.

\begin{theorem}\lb{T1.2} Let $J_0$ be a two-sided periodic Jacobi matrix and $J$ a
one-sided Jacobi matrix with Jacobi parameters $\{a_m,b_m\}_{m=1}^\infty$.
If $\sigma_\ess (J) = \sigma(J_0)$ and
\begin{equation} \lb{1.48}
\Sigma_\ac (J) = \sigma (J_0)
\end{equation}
then $d_m ((a,b),\calT_{J_0})\to 0$.
\end{theorem}

\begin{remark}  Using Theorem~\ref{T1.4} below, we will show that the hypotheses of this theorem can hold
while $(a,b)$ only approaches $\calT_{J_0}$ without actually having a limit.
\end{remark}

A two-sided $p$-periodic Jacobi matrix is said to have all gaps
open if the spectrum has exactly $p$ connected components---the
largest number  possible.  As explained in Section~\ref{s1a}, this
holds generically (indeed, on a dense open set).

Our next new result is

\begin{theorem}\lb{T1.3} Let $J_0$ be a two-sided periodic Jacobi matrix with all gaps open
and parameters $\{a_n^{(0)}, b_n^{(0)}\}_{n=-\infty}^\infty$.
Also, let $J$ be a one-sided Jacobi matrix with parameters $\{a_n,
b_n\}_{n=1}^\infty$ and spectral measure $d\mol$.  We assume that
$\sigma_\ess (J) = \sigma(J_0)$ and
\begin{equation} \lb{1.49}
\sum_j \dist (E_j,\sigma_\ess (J))^{1/2} <\infty
\end{equation}
where $\{E_j\}$ enumerates the eigenvalues of $J$ outside $\sigma(J_0)$.
Then
\begin{equation} \lb{1.50}
-\int_{\sigma(J_0)} \log \biggl( \f{d\mol_\ac}{dx}\biggr) \dist (x, \bbR\setminus
\sigma(J_0))^{-1/2}\, dx <\infty
\end{equation}
implies
\begin{equation}\lb{1.58a}
\lim_{N\to\infty}\, \sum_{j=1}^{pN} \log \biggl( \f{a_j}{a_j^{(0)}}\biggr)
\end{equation}
exists and lies in $(-\infty,\infty)$. Conversely, \eqref{1.50} holds so long as
\begin{equation} \lb{1.51}
\limsup_{N\to\infty} \, \sum_{j=1}^N \log \biggl( \f{a_j}{a_j^{(0)}}\biggr) >-\infty
\end{equation}
and in this case, the limit in \eqref{1.58a} exists and lies in $(-\infty,\infty)$.

Lastly, if \eqref{1.50} or \eqref{1.51} holds, then
\begin{equation} \lb{1.52}
\sum_{m=0}^\infty d_m ((a,b),\calT_{J_0})^2 <\infty
\end{equation}
and there exists $J_1\in\calT_{J_0}$, so that
\begin{equation}\lb{1.60a}
d_m (J,J_1)\to 0
\end{equation}
\end{theorem}

\begin{remarks} 1. Thus, when \eqref{1.49}--\eqref{1.51} hold, $J$ has a limit $J_1$
in $\calT_{J_0}$. In the normal direction to $\calT_{J_0}$, the convergence is $\ell^2$
(in the sense of \eqref{1.52}). But in the tangential direction, we only prove it has
a limit. It would be interesting to know what can be said about how slowly \eqref{1.60a}
can occur and to know if there are examples where \eqref{1.49}--\eqref{1.52} hold but
\begin{equation} \lb{1.60b}
\sum_{m=0}^\infty d_m (J,J_1)^2 =\infty
\end{equation}

2. Notice that \eqref{1.51} will only fail if the partial sums converge to $-\infty$.

3. The final statement that there exists $J_1\in\calT_{J_0}$ with
\eqref{1.60a} is not our result but one of Peherstorfer--Yuditskii
\cite{PY}. With our methods, we can prove that the $a$'s and $b$'s
approach a periodic limit only if we replace \eqref{1.49} with the
stronger assumption that the discrete spectrum is finite.

4. By \eqref{1.24newA}, all $a_j^{(0)}$ in \eqref{1.58a} and \eqref{1.51} can be replaced by
$\Capac(\sigma(J_0))$, the logarithmic capacity of the spectrum of $J_0$.
\end{remarks}

Our third new result is

\begin{theorem}\lb{T1.4} Let $J_0$ be a two-sided periodic Jacobi matrix with all
gaps open and parameters $\{a_n^{(0)}, b_n^{(0)}\}_{n\in\bbZ}$.
Let $J$ be a Jacobi matrix with parameters $\{a_n, b_n\}_{n=1}^\infty$ and spectral measure $d\mol$.
Then
\begin{equation}\label{1.52b}
\sum_{m=0}^\infty d_m ((a,b),\calT_{J_0})^2 <\infty
\end{equation}
if and only if
\begin{SSL}
\item[{\rm{(i)}}] $\sigma_\ess(J)=\sigma(J_0)$,
\item[{\rm{(ii)}}] $\displaystyle \sum_j \dist (E_j,\sigma_\ess (J))^{3/2} <\infty$, and
\item[{\rm{(iii)}}] $\displaystyle -\int_{\sigma(J_0)} \log \biggl( \f{d\mol_\ac}{dx}\biggr) \dist (x, \bbR\setminus
\sigma(J_0))^{1/2}\, dx <\infty$.
\end{SSL}
Here $\{E_j\}$ enumerates the {\rm{(}}discrete{\rm{)}} spectrum of $J$ outside $\sigma(J_0)$.
\end{theorem}

\begin{remarks} 1. Since (i)--(iii) are equivalent to \eqref{1.52b}, one may easily
construct examples where (i)--(iii) hold, but there is no $J_1$
with \eqref{1.60a}. This provides the examples promised in the
remark after Theorem~\ref{T1.2}. It also shows a stark difference
between \eqref{1.49}--\eqref{1.50} and (ii)--(iii).  In terms of
the spectral measure, this difference is reflected only in the
behavior near the band edges.

2. As we will see (Section~\ref{s10}), there are results even if all gaps are not open,
but for Theorem~\ref{T1.4} they are not so easy to express directly in terms of the
$a$'s and $b$'s.

3. A special case of part of Theorem~\ref{T1.4} is known, namely, Killip \cite{Kil02}
proved that $\Sigma_\ac (J) = \sigma (J_0)$ for $\{a_n, b_n\}_{n=0}^\infty$ obeying
$\sum_{n=1}^\infty \abs{a_n - a_n^{(0)}}^2 + \abs{b_n -b_n^{(0)}}^2 <\infty$ (which
is a strictly stronger hypothesis than \eqref{1.52b}).
\end{remarks}

\begin{theorem}\lb{T1.5} Let $J_0$ be a two-sided periodic Jacobi matrix and $J$ a Jacobi matrix
with Jacobi parameters $\{a_n, b_n\}_{n=1}^\infty$ and spectral measure $d\mol$. Suppose
\begin{equation} \lb{1.53}
\sum_{n=1}^\infty\, \abs{a_n-a_n^{(0)}} + \abs{b_n - b_n^{(0)}} <\infty
\end{equation}
Then \eqref{1.50} holds.
\end{theorem}

\begin{remarks}
1. Condition \eqref{1.53} can be replaced by
\begin{equation} \lb{1.55}
\sum_{n=1}^\infty d_n ((a,b), \calT_{J_0}) <\infty
\end{equation}
Indeed, if \eqref{1.55} holds, then \eqref{1.53} holds with
$\{a_n^{(0)}, b_n^{(0)}\}$ replaced by some fixed sequence in $\calT_{J_0}$.

2. As we will show (see Proposition~\ref{P2.5}), the theorems above continue to hold if $d_m$ is replaced by
\begin{equation} \lb{1.54}
\ti d_m ((a,b), (a',b'))=\sum_{k=0}^{p-1} (\abs{a_{m+k}-a'_{m+k}} +
\abs{b_{m+k} - b'_{m+k}})
\end{equation}
For OPUC, we need to sum $k$ from $0$ to $p$ in order to get an equivalence; see
the discussion at the end of Section~\ref{s3}.
\end{remarks}

In Section~\ref{s13}, we prove an OPUC analog of each of these theorems.
We need to replace the all-gaps-open hypothesis with a stronger one (that holds generically).
The deficiency is not so much with our method, but rather that an independent question, which is
known in the Jacobi case (independence of the Toda Hamiltonians), is currently unresolved in the
CMV case.  Our results confirm Conjectures~12.2.3 and 12.2.4 of \cite{OPUC2} as well as
Conjectures~12.2.5 and~12.2.6 in the (generic) special case that all gaps are open.

For the case of OPUC with a single gap, the analog of Theorem~\ref{T1.2} is known
and motivated Simon's conjectures in \cite{OPUC2}. In that case, the isospectral tori
are labelled by $a\in (0,1)$ and consist of $\{\alpha^{(\lambda)} :
\lambda\in\partial\bbD\}$ where $\alpha_n^{(\lambda)}=\lambda a$. Then $d_m (\alpha,
\calT) \to 0$ is equivalent to $\abs{\alpha_n}\to a$ and $\alpha_{n+1}/\alpha_n\to 1$.
This is often called the L\'opez condition. Bello--L\'opez \cite{BHLL} proved the
OPUC analog of Theorem~\ref{T1.2} for this case if $\sigma_\ess (J)=\sigma(J_0)$ is
strengthened to $\sigma (J)=\sigma (J_0)$ (the analog of Rakhmanov's result). The full
analog for this special case appears in Simon \cite{OPUC2}, Alfaro et al.\ \cite{ABMV},
and Barrios et al.\ \cite{BCL}.

Associated to each two-sided $p$-periodic Jacobi matrix, $J_0$, is a polynomial, $\Delta_{J_0}$,
of degree $p$, known as the discriminant.  This is a classical object described in detail in the next
section.  It is usually defined as the trace of the one-period transfer matrix.  It
is also the unique polynomial (with positive leading coefficient) such that
$$
\sigma(J_0) = \{ x : \Delta_{J_0}(x) \in [-2,2] \}
$$
In particular, two sequences of coefficients are isospectral if and only if they give rise
to the same discriminant.

The key to the proofs of our results is what we call the magic formula. Let $J$ be
a {\it two-sided\/} Jacobi matrix, then
\begin{equation} \lb{1.56}
\Delta_{J_0} (J) = S^p + S^{-p}
\end{equation}
if and only if $J\in\calT_{J_0}$.  Here $S$ is the right shift (cf. \eqref{1.35}).
In particular, \eqref{1.56} already implies that $J$ is periodic!
In the OPUC case, $\Delta$ is a polynomial in $z$ and $z^{-1}$.  It turns out that $\Delta(\calC)$ is always
self-adjoint; moreover,
\begin{equation} \lb{1.57}
\Delta_{\calC_0}(\calC) =S^p + S^{-p}
\end{equation}
if and only if $\calC\in\calT_{\calC_0}$.

It has been previously noted that for periodic $J_0$, one has
\[
\Delta_{J_0} (J_0) = S^p + S^{-p}
\]
That this holds for some polynomial in $J_0$ is in Na\u\i man \cite{Nai,Nai2}. That the
polynomial is the discriminant was found by Sebbar--Falliero \cite{SF}. After learning
of our results, L.~Golinskii has kindly pointed out to us that Na\u\i man \cite{Nai2}
also has a theorem which implies any $J$ obeying \eqref{1.56} is periodic, the core
of proving the converse. We will discuss this further in Section~\ref{s2}.

Nonetheless, the two facts that make \eqref{1.56} magical to us---that it characterizes the isospectral
torus and that it is ideal for the study of perturbations---seem to have escaped prior notice.

While $J$ may be tridiagonal and $\calC$ five-diagonal, both $\Delta(J)$ and $\Delta(\calC)$
are $2p+1$-diagonal, that is, vanishing except for the main diagonal and $p$ diagonals
above and below.  Thus, both will be tridiagonal if written as $p\times p$ blocks.
The key to our proofs will be to extend results from the $a_n\equiv 1$,
$b_n\equiv 0$ case to block tridiagonal matrices, and then use \eqref{1.56} or \eqref{1.57}
to study perturbations of the periodic case.

The magic formula is very powerful and opens up many other avenues for study:
\begin{SL}
\item[(a)] Szeg\H{o} and Jost asymptotics for periodic perturbations and, in particular,
the analogs of Damanik--Simon \cite{Jost1}.
\item[(b)] Periodic analogs of the results of Nevai--Totik \cite{NT89} and its various
recent extensions \cite{Jost2,Sim303,Jost3}.
\item[(c)] Analogs of the Strong Szeg\H{o} Theorem for the periodic case following
Ryckman's paper \cite{Ryck} for the Jacobi case.
\end{SL}

\smallskip
We should point out a major limitation of our results. If $B$ is a
disjoint finite union of closed intervals in $\bbR$ (or $\partial\bbD$), one can
construct an isospectral torus of Jacobi (or CMV) matrices whose recursion coefficients
are almost periodic. As discussed in Section~\ref{s1a}, these are strictly periodic if and only if
the harmonic measure of each interval is rational. There are obvious potential extensions
of Theorems~\ref{T1.1}--\ref{T1.5} to this setting, but except for Theorem~\ref{T1.1}
(where the method of \cite{LS_jdam} applies) and Theorem~\ref{T1.2} (where Section~\ref{new-s9}
has some extensions), we do not know how to prove them (or even if they are true).
There is no analog of $\Delta$ in the almost periodic case, so our method does not work directly.

Here is the plan of this paper. Section~\ref{s1a} reviews the theory of the (unperturbed) periodic
problem.  In Section~\ref{s2}, we prove the magic formula
for OPRL, and in Section~\ref{s3}, the magic formula for OPUC. While we will not
discuss Schr\"odinger operators in detail here, we discuss the magic formula
for such operators in Section~\ref{s4}. As we have mentioned, the magic formula brings
block Jacobi matrices into play, so Section~\ref{s5} discusses matrix-valued OPRL and
OPUC---mainly setting up notation. Section~\ref{s6} uses known results on Rakhmanov's
theorem for matrix-valued orthogonal polynomials to prove a Denisov-type
extension which we use in Section~\ref{s7} to prove Theorem~\ref{T1.2};
the section also proves half of Theorem~\ref{T1.1}. Section~\ref{new-s9} provides two
results that go beyond the periodic case to prove Denisov--Rakhmanov-type theorems
for special almost periodic situations. Section~\ref{s8}, following
\cite{Sim288}, proves the $P_2$ sum rule of Killip--Simon \cite{KS} and the $C_1$
sum rule for matrix-valued measures, and Section~\ref{s9} uses these results to prove
Theorems~\ref{T1.3} and \ref{T1.4}. Section~\ref{s10} explores what we can say if
gaps are closed. Section~\ref{s11} proves analogs to the Lieb--Thirring bounds of
Hundertmark--Simon \cite{HunS} as preparation for proving Theorem~\ref{T1.5}
in Section~\ref{s12}. Finally, Section~\ref{s13} discusses the OPUC results.

\smallskip
\subsection*{Acknowledgements}
It is a pleasure to thank Leonid Golinskii, Irina Nenciu,
Leonid Pastur, and Peter Yuditskii for useful discussions.

\smallskip
\subsection*{Note Added August, 2008}
During the refereeing of this paper, Remling (in \cite{R07}),
motivated in part by this paper, found a positive resolution of
the conjecture that, in the language of our Theorem~\ref{Tn9.5},
every set in $\calG$ is a Denisov-Rahkmanov set. His analysis
depends on a very interesting theorem on right limits of Jacobi
matrices with absolutely continuous spectrum -- it provides a new
approach to Denisov-Rahkmanov theorems.

\section{Review of the Periodic Problem} \label{s1a}

In this section, we collect some of the major elements in the strictly periodic case.  As this is textbook
material, we forgo proofs and historical discussion.
Full details can be found, for example, in \cite{CodLev,East,MagWin,OPUC2,SimonNew,Teschl,Toda}
and the references therein.

To discuss the strictly periodic case, we need to extend our operators to be two-sided, that is,
to act on $\ell^2(\bbZ)$.  In the Jacobi/OPRL case, we simply continue the tridiagonal pattern with parameters
$\{a_n,b_n\}_{n\in\bbZ}$. Two-sided (or extended) CMV matrices are formed as $\calC=\calL\calM$, where $\calL$ and $\calM$
are doubly infinite direct sums
\begin{align}\lb{3.12new}
\calL &= \cdots \oplus \Theta_{-2} (\alpha_{-2}) \oplus \Theta_0 (\alpha_0) \oplus \Theta_2 (\alpha_2) \oplus \cdots \\
\lb{3.13new}
\calM &= \hphantom{\oplus\Theta_0(}\cdots\oplus \Theta_{-1} (\alpha_{-1}) \oplus\Theta_1 (\alpha_1) \oplus\cdots
\end{align}
that are misaligned by one row/column, just as in \eqref{1.3b}--\eqref{1.3e}.

We adopt the convention of indexing the elements of matrices so that
\begin{equation} \label{zeroentry}
J_{11} = b_1 \qquad \calL_{00}=\bar\alpha_0 \qquad \calM_{00}=-\alpha_{-1}
\end{equation}
except $\calM_{00} = 1$ in the one-sided case.

\subsection{Transfer Matrices.}  Let $J$ be a two-sided Jacobi matrix.  A sequence $\{u_n\}$
obeys $(J-x)u\equiv 0$ if and only if
\begin{equation} \lb{1.27}
a_{n} u_{n+1} + (b_{n}-x) u_{n} + a_{n-1} u_{n-1} =0
\end{equation}
or, what is equivalent,
\begin{equation} \lb{1.28}
\begin{pmatrix} u_{n+1} \\ a_{n} {u_n} \end{pmatrix} = \Lambda_n \begin{pmatrix} u_{n} \\ a_{n-1} u_{n-1} \end{pmatrix}
\end{equation}
with
\begin{equation} \lb{1.14}
\Lambda_n(x) = \frac{1}{a_{n}} \begin{pmatrix} x-b_{n} &  -1 \\ a_{n}^2 & 0 \end{pmatrix}
\end{equation}
Note that the desire to have $\Lambda_n$ depend only on one pair $(a_{n},b_{n})$ and have determinant equal to one
resulted in the factors $a_{n}$ and $a_{n-1}$ appearing in \eqref{1.28}.  (The same price is usually paid when writing
Sturm--Liouville equations as first-order systems.) The choice is not the most common one
(although it is used in Pastur--Figotin \cite{PF}), but we feel it is the `right' one since,
in particular, $\det (\Lambda_n(x))=1$.

In the OPUC case we define
\begin{equation}\lb{1.15}
M_n(z) = \rho_n^{-1}
\begin{pmatrix}
z & -\bar\alpha_n  \\
-\alpha_n z & 1
\end{pmatrix}
\end{equation}
which encodes the recurrence relation \eqref{1.2}:
\begin{equation} \lb{1.28C}
\begin{pmatrix} \varphi_{n+1}(z) \\ \varphi^*_{n+1}(z) \end{pmatrix}
    = M_n(z) \begin{pmatrix} \varphi_n(z) \\ \varphi^*_n(z) \end{pmatrix}
\end{equation}
We will now explain the link to formal (i.e., not necessarily
$\ell^2$) eigenvectors of a two-sided CMV matrix, $\calC$. It is
not as simple as \eqref{1.28}.

\begin{lemma}
Suppose $(\calC-z)u=0$ with $z\neq 0$ and let $v=\calZ^{-1}\calM u$ where $\calZ$ denotes the diagonal matrix with entries
$$
\calZ_{jj} = \begin{cases} z &: \text{ $j$ odd} \\ 1 & : \text{ $j$ even }\end{cases}
$$
and $\calM$ is as in \eqref{3.13new}. Then
\begin{equation} \label{CMVtransfer}
z \begin{pmatrix} u_{2n+2} \\ v_{2n+2} \end{pmatrix} = z
M_{2n+1}(z) \begin{pmatrix} v_{2n+1} \\ u_{2n+1} \end{pmatrix} =
M_{2n+1}(z) M_{2n}(z) \begin{pmatrix} u_{2n} \\ v_{2n}
\end{pmatrix}
\end{equation}
\end{lemma}

\begin{proof}
The key observation used to verify \eqref{CMVtransfer} is
\begin{equation} \label{ThetaA}
\begin{pmatrix} zy \\ y' \end{pmatrix} = \Theta(\alpha_n) \begin{pmatrix} x \\ x' \end{pmatrix}
\iff
\begin{pmatrix} x' \\ y' \end{pmatrix} = M_n(z) \begin{pmatrix} y \\ x \end{pmatrix}
\end{equation}
This follows by simple algebraic manipulations:
\begin{align*}
&\!\!\!\begin{pmatrix} zy \\ y' \end{pmatrix} = \Theta(\alpha_n) \begin{pmatrix} x \\ x' \end{pmatrix} \\
\iff\ & \bar\alpha_n x + \rho_n x' = zy \quad\text{and}\quad \rho_n x - \alpha_n x' = y' \\
\iff\ & x' = \rho_n^{-1}(zy - \bar\alpha_n x)  \quad\text{and}\quad y' = \rho_n x - \alpha_n x' \\
\iff\ & x' = \rho_n^{-1}(zy - \bar\alpha_n x)  \quad\text{and}\quad y' = \rho_n^{-1}  (-\alpha_n zy + x)  \\
\iff\ &\!\!\!\begin{pmatrix} x' \\ y' \end{pmatrix} = M_{n}(z) \begin{pmatrix} y \\ x \end{pmatrix}
\end{align*}

With \eqref{ThetaA} now in hand, we may argue as follows:
\begin{align*}
& (\calC-z)u=0 \\
\iff\ & zu = \calL\calM u \\
\iff\ & v:=\calZ^{-1} \calM u \quad\text{obeys}\quad z u = \calL\calZ v \\
\iff\ & \calZ v = \calM u \quad\text{and}\quad z u = \calL\calZ v \\
\iff\ &\!\!\! \begin{pmatrix} z v_{2n-1} \\ v_{2n} \end{pmatrix} = \Theta(\alpha_{2n-1})
    \begin{pmatrix} u_{2n-1} \\ u_{2n} \end{pmatrix} \quad\text{and}\quad \begin{pmatrix} z u_{2n} \\ z u_{2n+1} \end{pmatrix}
    = \Theta(\alpha_{2n}) \begin{pmatrix} v_{2n} \\ z v_{2n+1} \end{pmatrix} \\
\iff\ &\!\!\! \begin{pmatrix} u_{2n} \\ v_{2n} \end{pmatrix} = M_{2n-1}(z) \begin{pmatrix} v_{2n-1} \\ u_{2n-1} \end{pmatrix}
    \quad\text{ and }\quad \begin{pmatrix} z v_{2n+1} \\ z u_{2n+1} \end{pmatrix} = M_{2n}(z) \begin{pmatrix} u_{2n} \\ v_{2n} \end{pmatrix} \\
\iff\ &\!\!\! \begin{pmatrix} u_{2n+2} \\ v_{2n+2} \end{pmatrix} = M_{2n+1}(z) \begin{pmatrix} v_{2n+1} \\ u_{2n+1} \end{pmatrix}
    \quad\text{and}\quad z \begin{pmatrix} v_{2n+1} \\ u_{2n+1} \end{pmatrix} = M_{2n}(z) \begin{pmatrix} u_{2n} \\ v_{2n}
\end{pmatrix}
\end{align*}
which are the two parts of \eqref{CMVtransfer}.
\end{proof}

\subsection{The Discriminant.}  As in Sturm--Liouville theory, the discriminant is defined as
the trace of the one-period transfer matrix:
\begin{equation} \lb{1.21new}
\Delta(z) = \tr( T(z) )
\end{equation}
where
\begin{equation} \lb{1.21Anew}
T(z) = \begin{cases}
\ \Lambda_{p}(z) \cdots \Lambda_2(z) \Lambda_1(z) &\text{ (OPRL)} \\[0.25ex]
\ z^{-p/2} M_{p-1}(z) \cdots M_1(z) M_0(z) & \text{ (OPUC)}
\end{cases}
\end{equation}
In the OPUC case, $p$ is even.  Also, the factor $z^{-p/2}$ is there to cancel the extra factor of $z$ on the
left-hand side of \eqref{CMVtransfer}.  From a strictly
OPUC point of view, it is more natural to omit this factor (as in \cite{OPUC1,OPUC2}); however, as
the magic formula is an operator identity, we have elected to use the definition best
adapted to this perspective.  The only negative side effect of this choice is that our Lyapunov
exponent (defined below) differs by $-\tfrac12\log|z|$ from that in \cite{OPUC1,OPUC2}.

For OPRL, the discriminant is a real polynomial of degree $p$ with leading behavior
\begin{align}
\lb{2.5xnew}
\Delta(x) &= (a_1 \dots a_p)^{-1} \biggl[\, \prod_{j=1}^p (x-b_j) + O(x^{p-2})\biggr] \\
\lb{2.4new}
          &= (a_1 \cdots a_p )^{-1} \biggl[ x^p - \biggl(\, \sum_{j=1}^p b_j\biggr) x^{p-1} + O(x^{p-2})\biggr]
\end{align}

For OPUC, it is a Laurent polynomial of total degree $p$ with
\begin{equation} \lb{1.13}
\Delta(\bar z) = \ol{\Delta(1/z)}
\end{equation}
so $\Delta$ is real-valued on $\partial\bbD$.  Moreover,
\begin{equation} \lb{1.13a}
\Delta(z) =(\rho_0 \rho_1 \cdots \rho_{p-1})^{-1} (z^{p/2} + \cdots + z^{-p/2})
\end{equation}

\subsection{The Lyapunov Exponent}
On an exponential scale, the behavior of formal eigenfunctions is determined by
the Lyapunov exponent
\begin{align} \lb{1.16}
\gamma(z) &= \lim_{n\to\infty}\, \f{1}{np} \, \log \| T^n(z) \| \\
&= \f{1}{p} \, \log \text{(spectral radius of $T(z)$)} \notag \\
&= \f{1}{p}\, \log |\lambda_+(z)| \lb{1.18x}
\end{align}
where $\lambda_\pm$ are the eigenvalues of $T(z)$ with the convention $\abs{\lambda_+} \geq \abs{\lambda_-}$.

As $\det(T(z))\equiv 1$, these eigenvalues are the roots of
\begin{equation} \lb{1.22}
\lambda^2 -\Delta(z)\lambda + 1=0
\end{equation}
which implies
\begin{equation} \lb{1.23}
\lambda_\pm(z)  = \f{\Delta(z)}{2} \pm \f{\sqrt{\Delta^2(z) -4}}{2}
\end{equation}
and so
\begin{equation} \lb{1.24}
\gamma(z) = \tfrac{1}{p} \log \left| \tfrac12 \Delta(z) + \tfrac12 \sqrt{\Delta(z)^2-4}\right|
\end{equation}

\subsection{Gaps and Bands.}
Our recurrence relations admit bounded solutions for a given $z$ if and only if
$\Delta(z)\in[-2,2]$.  In the Jacobi/OPRL case, this is a collection of intervals in $\bbR$.
For CMV/OPUC, it is a collection of arcs in $\partial \bbD$.
In either case, one may partition this set into $p$ bands.  These are the closures of the (disjoint)
regions where $\Delta(z)\in (-2,2)$.  These can only intersect at the `band edges', $\Delta^{-1} (\{-2,2\})$.

The open gaps are the intervals/arcs that are complementary to the bands---excluding the two semi-infinite
intervals in the OPRL case.   When two bands touch, we refer to the common band edge as a closed gap.

$\Delta^2-4$ has simple zeros at the edges of the open gaps and double zeros at the closed gaps; indeed,
this is a complete list of its zeros.  It is possible to distinguish whether these zeros correspond to
$\Delta(z)=\pm 2$ from the fact that there must be two zeros of $\Delta\pm2$ between consecutive zeros of $\Delta\mp2$
and the fact that $\Delta$ has positive leading coefficient.

\subsection{Spectrum.} In both cases, the spectrum of the two-sided operator (acting on $\ell^2(\bbZ)$)
is the union of the bands: $\sigma = \Delta^{-1} ([-2,2])$.  It is purely absolutely continuous and of
multiplicity two.

The spectrum of a two-sided $p$-periodic operator uniquely determines its discriminant; see Lemma~\ref{L2.4}.
One consequence of this was noted already in the introduction: isospectral tori are the
classes of $p$-periodic recurrence coefficients that lead to the same discriminant.

In the case of a one-sided operator, the essential spectrum remains $\Delta^{-1} ([-2,2])$; it is absolutely
continuous with multiplicity one.  In addition, up to one eigenvalue may appear in each open gap.

\subsection{Potential Theory.} From the way they are defined, one can see
that $\gamma(z)$ vanishes on the bands and is both positive and harmonic in the
complement (in the OPUC case one must also exclude $z=0$).  This leads to the solution of
the Dirichlet problem for a charge at infinity,
\begin{equation}\label{1.24new}
g_{\bbC\setminus \sigma} (z;\infty) = \begin{cases}
\ \f{1}{p} \log \left| \f{\Delta}{2} + \sqrt{\f{\Delta^2}{4}-1}\right| & \text{(OPRL)} \\[0.25ex]
\ \f{1}{2} \log \abs{z} + \f{1}{p} \log \left| \f{\Delta}{2} + \sqrt{\f{\Delta^2}{4}-1}\right| & \text{(OPUC)}
\end{cases}
\end{equation}
and so to the logarithmic capacity of the spectrum,
\begin{equation}\label{1.24newA}
\Capac(\sigma) = \begin{cases}
\bigl(\prod_{j=1}^p a_j\bigr)^{1/p}           & \text{(OPRL)} \\[1.25ex]
\bigl(\prod_{j=0}^{p-1} \rho_j\bigr)^{1/p}    & \text{(OPUC)}
\end{cases}
\end{equation}

\subsection{Harmonic Measure.} Taking normal derivatives in \eqref{1.24new} leads to a formula
for harmonic measure on $\sigma$ (aka equilibrium measure for the logarithmic potential),
\begin{equation} \lb{1.25}
d\nu = \begin{cases}
\f{2}{p} \, \f{\abs{\Delta'(x)}}{\sqrt{4-\Delta^2(x)}}\, \f{dx}{2\pi} & \text{(OPRL)} \\[1.25ex]
\f{2}{p} \, \f{\abs{\Delta'(e^{i\theta})}}{\sqrt{4-\Delta^2(e^{i\theta})}}\, \f{d\theta}{2\pi}    & \text{(OPUC)}
\end{cases}
\end{equation}
where $\supp(d\nu)= \sigma=\{ z :\abs{\Delta(z)}\leq 2\}$.
({\it Note}: \cite{OPUC2} has $1/p$ rather than $2/p$, but that is an error.)

Recognizing
\begin{equation} \lb{1.26}
\f{\Delta'(x)}{\sqrt{4-\Delta^2 (x)}} = \f{d}{dx}\, \arccos \biggl( \f{\Delta(x)}{2}\biggr)
\end{equation}
we see that the harmonic measure of each band is exactly $1/p$.
In particular, the connected components of the union of the bands all have
rational harmonic measure. This gives strong restrictions on sets that can be bands.
In the OPRL case, rational harmonic measure of connected components is also sufficient for a set to
be the spectrum of a periodic Jacobi matrix. In the OPUC case, there is an additional condition
needed: after breaking the bands into arcs of harmonic measure $1/p$, the
harmonic midpoints $\{\zeta_j\}_{j=1}^p$ of these intervals must obey $\prod_{j=1}^p \zeta_j =1$.
Clearly, the condition on the harmonic midpoints can be achieved by simply rotating $\sigma$.
Discarding this condition gives rise to Verblunsky coefficients that are
$p$-automorphic, $\alpha_{n+p} = e^{i\phi} \alpha_n$, rather than $p$-periodic.

\subsection{Thouless Formula}\label{ssTF}
Harmonic measure appears naturally in the theory in several other ways.  It is the density of states measure:
\begin{equation}\label{DOS}
 \lim_{N\to\infty} \frac{1}{2N} \sum_{n=-N}^N f(J)_{nn} = \int f(x) \,d\nu(x)
\end{equation}
for every polynomial (or continuous function) $f$.  The same formula holds with $\calC$ replacing $J$.
This connection, or more precisely the resulting expression for $\gamma(z)$ in terms of $\Capac(\sigma)$ and
the logarithmic potential of $d\nu$, is known as the Thouless formula.

Two further characterizations of $d\nu$ involve the orthogonal polynomials.  $d\nu$ is the weak limit of
$\f{1}{N} \sum_{n=0}^{N-1} p_n^2 (x)\, d\mol(x)$ (resp. $ \f{1}{N} \sum_{n=0}^{N-1}
\abs{\varphi_n (e^{i\theta})}^2 d\mu(\theta)$).  It is also the limiting density of zeros, that is, the
weak limit of the probability measures, $d\nu_n$, which give weight $1/n$ to each of the zeros of $p_n$ (resp.\
$\varphi_n$).  These two characterizations are closely linked to \eqref{DOS}; however, in the OPUC case
one should keep in mind that for each $n$, the zeros of $\varphi_n$ lie strictly inside $\bbD$.

\subsection{Floquet Theory.}
Looking at the eigenvalues of $T$, one sees that when $\lambda_+\neq \lambda_-$, there is
a basis of formal (i.e., non-$\ell^2$) eigenfunctions obeying
\begin{equation} \lb{1.29}
u_{m+kp} = \lambda_\pm^k u_m
\end{equation}
If $\lambda_+ =\lambda_-$, which happens precisely at the band edges, then both are $\pm1$.  If the edge
abuts an open gap, there is only one eigenfunction obeying \eqref{1.29} since $T$ has a Jordan block structure.
At closed gaps, $T=\pm \bdone$ and so all solutions obey \eqref{1.29}.

Solutions obeying
\begin{equation} \lb{1.30}
u_{m+kp} = e^{ik\theta} u_m
\end{equation}
are called Floquet solutions and $e^{i\theta}$ is called the Floquet index;
they have much the same role as plane waves in Fourier analysis. Since $\lambda_-=\lambda_+^{-1}$,
if \eqref{1.30} has a solution, then $e^{-i\theta}$ is also a Floquet index.

In the OPRL case,
\begin{equation} \lb{1.31}
\text{\eqref{1.30} holds} \iff \Delta(x) =2\cos\theta
\end{equation}
Thus, by the discussion above, for each $\theta\in (0,\pi)$, \eqref{1.30} or \eqref{1.31}
holds for exactly $p$ values of $x$: $x_1(\theta) <x_2(\theta) <\cdots < x_p(\theta)$.
These $x_j(\theta)$ are known as the band functions.

The changes in the OPUC case are purely notational:
\begin{equation} \lb{1.34}
\text{\eqref{1.30} holds} \iff \Delta(z) = 2\cos\theta
\end{equation}
For $\theta\in (0,\pi)$, this has $p$ solutions all of which lie in $\partial\bbD$.

\subsection{Direct Integrals.}
Let $S$ denote the right shift,
\begin{equation} \lb{1.35}
(Su)_n = u_{n-1}
\end{equation}
If the sequences of coefficients are $p$-periodic, then $J$ (or $\calC$) commutes with $S^p$,
which means that the two operators can be `simultaneously diagonalized'.  We elaborate this
point in the OPRL case; the OPUC is almost identical.

Let us write
$$
\calH_p := \int^\oplus \ell^2_\theta \, \f{d\theta}{2\pi} = L^2\bigl([0,2\pi), \tfrac{d\theta}{2\pi}; \bbC^p\bigr)
$$
where $\ell^2_\theta$ is the $p$-dimensional Hilbert space
$$
\ell^2_\theta = \{ u | u_{n+p}=e^{i\theta} u_n\}
\qquad
\langle u |  v \rangle_\theta = \sum_{n=1}^p \bar{u}_nv_n
$$

{}From Fourier analysis, there is a unitary operator $\calF \colon \ell^2 (\bbZ)\to \calH_p$
so that $\calF S^p\calF^{-1}$ is multiplication by $e^{i\theta}\bdone$ and $\calF J\calF^{-1}$
acts fiber-wise (i.e., on each $\ell^2_\theta$) as a $p\times p$ matrix, $J(\theta)$.
In particular, the eigenvalues of $J(\theta)$ are the solutions of \eqref{1.31}, that is,
they are the band functions $x_j(\theta)$.

\subsection{Hyperelliptic Riemann Surfaces.} As $\sqrt{\Delta^2 -4}$ appears repeatedly in the
theory, it is natural that the associated Riemann surface should enter the analysis.
$\Delta^2 -4$ has simple zeros at the edges of open gaps and at $\inf \sigma(J)$ and $\sup \sigma(J)$.
It has double zeros at the closed gaps.  Let $\ell$ denote the number of open gaps, then
$\sqrt{\Delta^2 -4}$ has square root singularities at
$2(\ell+1)$ points, and so its natural analyticity domain is the genus $\ell$ Riemann
surface, $\calS$, obtained by taking two copies of $\bbC\setminus\sigma(J)$, gluing
at the bands and adding points at $\infty$. There is a natural projection $\pi \colon
\calS \to \bbC \cup\{\infty\}$ which is $2$ to $1$ except at the branch points of
$\sqrt{\Delta^2 -4}$. A similar analysis works for OPUC, but now there are $\ell$ gaps
and the genus is $\ell-1$.

\subsection{Minimal Herglotz and Carath\'eodory Functions.} For a half-line
periodic Jacobi matrix, the $m$ function is defined by
\begin{equation} \lb{1.40}
m(z) = \langle\delta_0, (J-z)^{-1}\delta_0\rangle\qquad\forall z\in\bbC\setminus\sigma(J)
\end{equation}
This can be shown to obey a quadratic equation with polynomial coefficients
\begin{equation} \lb{1.41}
A(z) m(z)^2 + B(z) m(z) + C(z) =0
\end{equation}
Moreover, these coefficients can be chosen to obey
\begin{equation} \lb{1.42}
B^2 -4AC=\Delta^2 - 4
\end{equation}
This implies that $m(z)$ has meromorphic continuation to $\calS$. Indeed, $m$ has minimal
degree (i.e., degree $\ell+1$ in the OPRL case and $\ell$ in the OPUC case) among all
meromorphic functions on $\calS$ that are not of the form $g\circ\pi$ with $g$
meromorphic on the Riemann sphere. It can be shown that there is a one-one
correspondence between minimal meromorphic functions obeying $\Ima m(z) >0$
if $\Ima z >0$ and $m(z) =-z^{-1} + O(z^{-2})$ on the top sheet of $\calS$
and all periodic Jacobi parameters with the same $\Delta$. (We call these
minimal Herglotz functions.)

There is a similar description for OPUC, but now one uses
\begin{equation} \lb{1.43x}
F(z) =\langle\delta_0, (\calC+z) (\calC-z)^{-1}\delta_0\rangle
\end{equation}
which obeys $\Real F(z) >0$ if $\abs{z}<1$ and $F(0)=1$. Again $F$ obeys a quadratic
equation, showing that $F$ has a meromorphic continuation to $\calS$ of minimal degree,
and again there is a one-one correspondence between all $\{\alpha_n\}_{n=0}^{p-1}$
with the same $\Delta$ and all minimal Carath\'eodory functions.

\subsection{Dirichlet Data.} One can describe the set of minimal
Herglotz functions in terms of their poles. For each open gap, $\{G_j\}_{j=1}^\ell$,
$\pi^{-1} (\ol{G}_j)\equiv T_j$ is a circle since $\pi$ is $2$ to $1$ on $G_j$ and
one-one on $\ol{G}_j\setminus G_j$. A meromorphic Herglotz function has
$\ell+1$ simple poles, one at $\infty$ on the second sheet and the other $\ell$,
one in each $T_j$. Thus, the set of meromorphic Herglotz functions is homeomorphic to
$\bigtimes_{j=1}^\ell T_j$ under the bijective map from such functions to its poles.
A similar analysis holds for OPUC but now there is no pole at infinity, there are
$\ell$ gaps, and $\bigtimes_{j=1}^\ell T_j$ describes the possible poles. The difference
is that for OPRL, the dimension of the torus is $\ell$, and for OPUC it is $\ell-1$.

\subsection{Isospectral Tori.}\label{ssIT} By combining the bijective maps from periodic
OPRL to minimal Herglotz functions and of such functions to Dirichlet data, we see
for a $\Delta$ of period $p$ with $\ell$ gaps,
$$
\{(a_n, b_n)_{n=1}^p : \text{the discriminant is }\Delta\}
$$
is an $\ell$-dimensional torus in $\bbR^{2p}$. Generically, $\ell=p-1$.
In the OPUC case, generically $\ell=p$ and the torus is naturally embedded in $\bbC^p$.
This torus is the isospectral torus which we will denote by $\calT$ or $\calT_{J_0}$
if a given periodic $J_0$ underlies our construction.  For clarity of exposition, we will
typically blur the distinction between $p$-tuples $(a_n, b_n)_{n=1}^p$ and the corresponding
infinite sequences $\{a_n,b_n\}_{n\in\bbZ}$ of period $p$. Because of our perturbation theory
viewpoint, we use $J_0$ to label the torus, but we emphasize that from another point of view,
the torus is associated to the set $\sigma_\ess(J_0)$ and not to $J_0$.

\subsection{Isospectral Flows.} The fact that spaces of $p$-periodic coefficients foliate into tori
suggests that there is some kind of completely integrable system in the background.
That is true: it is the Toda flow in the OPRL case and the defocusing Ablowitz--Ladik
flow in the OPUC case. Since we will not need these below, we say no more about them,
but see Chapter~6 of \cite{SimonNew} for the OPRL case and Section~11.11 of \cite{OPUC2}
for OPUC.

\section{The Magic Formula for Jacobi Matrices} \lb{s2}

Our goal in this section is to prove

\begin{theorem}\lb{T2.1} Let $J_0$ be a two-sided $p$-periodic Jacobi
matrix with discriminant $\Delta_{J_0}(x)$ and isospectral torus $\calT_{J_0}$.
Let $J$ be a two-sided {\rm{(}}not a priori periodic{\rm{)}} Jacobi matrix. Then
\begin{equation} \lb{2.1}
\Delta_{J_0}(J) = S^p + S^{-p} \iff J\in\calT_{J_0}
\end{equation}
where $S$ is the right shift, \eqref{1.35}, on $\ell^2(\bbZ)$.
\end{theorem}

We provide two proofs of the `harder' direction $\Rightarrow$ or rather of
\begin{equation}\label{2.11a}
\Delta_{J_0}(J) =S^p + S^{-p}\Rightarrow J\text{ is periodic}
\end{equation}
which is the key step. Our first proof is immediately below;
the second, suggested to us by L.~Golinskii, appears after Lemma~\ref{L2.4A}.

\begin{lemma}\lb{L2.2} Let $\ell=1,2,\dots$. Then
\begin{equation} \lb{2.2}
(J^\ell)_{m,m+k} = \begin{cases}
0 & k>\ell \\
a_m a_{m+1} \cdots a_{m+k-1} & k=\ell \\
a_m a_{m+1} \cdots a_{m+\ell-2} (b_m + b_{m+1} + \cdots + b_{m+\ell -1}) & k=\ell-1
\end{cases}
\end{equation}
\end{lemma}

\begin{proof} Writing
\begin{equation} \lb{2.3}
(J^\ell)_{m,m+k} = \sum_{i_1, \dots, i_{\ell-1}} J_{m, i_1} J_{i_1, i_2} \cdots
J_{i_{\ell-1}, m+k}
\end{equation}
We see that since $J$ is tridiagonal, all terms are zero if $k>\ell$, that we must
have (with $i_0\equiv m$, $i_\ell = m+k$) that $i_q - i_{q-1}=1$ for $q=1, \dots,
\ell$ if $k=\ell$, and that if $k=\ell-1$, $i_q -i_{q-1}=1$ for all but one $q\in
\{1,\dots,\ell\}$ and it is zero for that $q$.
\end{proof}

\begin{lemma}\lb{L2.4} If $J$ and $J_0$ are periodic, then
$\sigma(J)=\sigma(J_0)$ if and only if $\Delta_J=\Delta_{J_0}$.
\end{lemma}

\begin{remark}
This lemma says that the spectrum determines the discriminant and
vice versa.  That the discriminant determines the spectrum is elementary:
$\sigma = \{ x : \Delta(x)\in[-2,2]\}$.  Therefore we only prove the other direction---indeed,
we give two proofs.
\end{remark}

\begin{proof}[First Proof] Harmonic measure $d\nu$ is intrinsic to the set $\sigma$; it is
the solution of an electrostatic problem there.  But then $d\nu$ determines $\Delta$ via \eqref{1.25}.
\end{proof}

\begin{proof}[Second Proof] $\sigma$ determines the gaps---even closed gaps---via
harmonic measure. The gap edges determine the zeros of $\Delta -2$ and so $\Delta -2$
up to a constant. The zeros of $\Delta +2$ then determine the constant.
\end{proof}

\begin{proof}[Proof of Theorem~\ref{T2.1}] For all $\theta\in[0,2\pi)$, $J(\theta)$ is self-adjoint
and so diagonalizable.  Moreover, the eigenvalues of $J_0(\theta)$ are precisely the roots of
$\Delta(x)=2\cos(\theta)$. Thus%
\begin{equation} \lb{2.5}
\Delta_{J_0}(J_0 (\theta)) = (2\cos\theta)\bdone
\end{equation}
But then $\Delta_{J_0}(J_0)$ and $S^p + S^{-p}$ both have direct integral decomposition with fibers $(2\cos\theta)\bdone$,
so
\begin{equation} \lb{2.6}
\Delta_{J_0}(J_0) = S^p + S^{-p}
\end{equation}
Since $J\in \calT_{J_0}\Rightarrow \Delta_J=\Delta_{J_0}$, this proves $\Leftarrow$ in \eqref{2.1}.

Now suppose LHS of \eqref{2.1} holds. By \eqref{2.4new}, \eqref{2.2}, and
\begin{equation} \lb{2.7}
(S^p + S^{-p})_{m, m+p}=1 \qquad
(S^p + S^{-p})_{m, m+p-1} =0
\end{equation}
this implies
\begin{equation} \lb{2.8}
a_m \cdots a_{m+p-1} = a_1^{(0)} \cdots a_p^{(0)}
\end{equation}
and
\begin{equation} \lb{2.9}
\sum_{j=0}^{p-1} b_{m+j} =\sum_{j=0}^{p-1} b_{j+1}^{(0)}
\end{equation}
In particular,
\[
a_m\cdots a_{m+p-1} = a_{m+1} \cdots a_{m+p} \qquad
\sum_{j=0}^{p-1} (b_{m+j+1}-b_{m+j}) =0
\]
which lead to
\begin{equation} \lb{2.10}
a_m = a_{m+p} \qquad b_m = b_{m+p}
\end{equation}
so $J$ is periodic.

Since $J$ is periodic, $\Delta_J(J)=S^p + S^{-p}$; moreover, $\Delta_{J_0}(J)=S^p + S^{-p}$ by hypothesis.
Thus we learn that applying the polynomial $\Delta_{J}-\Delta_{J_0}$ to $J$ gives zero.  By the $k=\ell$ case
of \eqref{2.2}, it must therefore be the zero polynomial, that is, $\Delta_{J}=\Delta_{J_0}$.
Lemma~\ref{L2.4} now completes the proof.
\end{proof}

\begin{remarks}
1. Showing that $J$ was periodic only required equality in $\Delta_{J_0}(\calT)=S^p + S^{-p}$,
for the two most extreme upper (or lower) diagonals.  Nevertheless, $J\in\calT_{J_0}$ requires equality everywhere.

2. We need not suppose a priori that each $a_n >0$ and
can allow some $a_n =0$ ($J$ can still be defined on $\ell^2 (\bbZ)$), for \eqref{2.8}
implies that if LHS of \eqref{2.1} holds, then each $a_n >0$.
\end{remarks}

We now turn to our second proof of \eqref{2.11a}.

\begin{lemma}[Na\u\i man \cite{Nai2}]\lb{L2.4A} Let $A$ be a two-sided {\rm{(}}bounded{\rm{)}}
infinite matrix of finite width {\rm{(}}i.e., for some $w$, we have that $\abs{k-\ell} > w
\Rightarrow A_{k\ell} =0${\rm{)}}. Suppose
\begin{equation} \lb{2.11b}
[A, S^p + S^{-p}] =0
\end{equation}
for some $p$, then
\begin{equation} \lb{2.11c}
[A,S^p] =0
\end{equation}
\end{lemma}

\begin{remarks} 1. This is Lemma~2 in \cite{Nai2}; no proof is given.

2. $[A,B]\equiv AB-BA$

3. \eqref{2.11c} has an equivalent form:
\begin{equation} \lb{2.11d}
[A,S^p]=0 \iff A_{k+p,\ell+p}=A_{k,\ell} \qquad \text{for all $k,\ell$}
\end{equation}

4. As \eqref{2.11d} shows,  $[J,S^p] =0$ for a Jacobi matrix if and only if $a_k$ and $b_k$
are $p$-periodic.
\end{remarks}

\begin{proof} Since $A$ has finite width, we can find diagonal matrices $D_{k_1},D_{k_1+1}, \ldots,
D_{k_2}$ with $D_{k_1}\neq 0\neq D_{k_2}$, so that
\begin{equation} \lb{2.11e}
A=\sum_{j=k_1}^{k_2} D_j S^j
\end{equation}
Since $D_j$ is diagonal, so is $S^p D_j S^{-p}$. Thus
\begin{gather*}
(S^p + S^{-p})A = \sum_{j=k_1}^{k_2} (S^p D_j S^{-p}) S^{j+p} + \sum_{j=k_1}^{k_2}
(S^{-p} D_j S^p) S^{j-p} \\
A (S^p + S^{-p}) =\sum_{j=k_1}^{k_2} D_j S^{j+p} + \sum_{j=k_1}^{k_2} D_j S^{j-p}
\end{gather*}
Since the composition \eqref{2.11e} uniquely determines each $D_j$, \eqref{2.11b} implies
\[
S^p D_{k_2} S^{-p} = D_{k_2}
\]
that is, $D_{k_2}$ is periodic. Thus, $D_{k_2} S^p$ commutes with $S^p + S^{-p}$, so
we can remove it from \eqref{2.11e} without losing \eqref{2.11b}. This shows inductively
that each $D_j$ is periodic.
\end{proof}

\begin{proof}[Second proof of \eqref{2.11a}] $J$ commutes with $\Delta (J)$, so \eqref{2.11c}
holds.
\end{proof}

Our next goal is to compare $\ti d_m ((a,b),\calT)$ given by \eqref{1.54} and $d_m
((a,b),\calT)$ given by \eqref{1.43}.  As well as satisfying natural curiosity, this
relation also plays an important role (via Theorem~\ref{T9.11}) in the proofs of Theorems~\ref{T1.3} and~\ref{T1.4}.

To capture the essence of what follows, let us pause to ponder the following:
suppose $\tilde d_m((a,b),\calT)=0$ for all $m$, does this mean that $(a,b)\in\calT$?
The hypothesis tells us that each length-$p$ block belongs to the
isospectral torus; it does not a priori even guarantee that the coefficients are periodic.
Example~\ref{E4.5} shows that periodicity can fail in the OPUC case.  However, such problems do not
arise for OPRL. The reason is simple: within the isospectral torus, $a_1,\ldots,a_{p-1}$ determines $a_p$ and
$b_1,\ldots,b_{p-1}$ determines $b_p$.


\begin{proposition}\lb{P2.5} Given a $p$-periodic Jacobi matrix $J_0$, $1\leq q \leq \infty$, and $\veps>0$,
there is a constant $C$ so that
\begin{equation} \lb{2.13}
e^{1-p} \bigl\| \tilde d_m \bigl((a,b),\calT_{J_0}\bigr) \bigl\|_{\ell^q}
\leq \bigl\| d_m \bigl((a,b),\calT_{J_0}\bigr) \bigl\|_{\ell^q}
    \leq C \bigl\| \tilde d_m \bigl((a,b),\calT_{J_0}\bigr) \bigl\|_{\ell^q}
\end{equation}
for all sequences $\{(a_n,b_n)\}$ obeying $\veps^{-1} > a_n > \veps >0$.  All $\ell^q$ norms are over $m\in\{1,2,3,\ldots\}$.
\end{proposition}

The key input is

\begin{lemma}\lb{L2.6}
Given $\{(a_n,b_n)\}$ obeying $\veps^{-1} > a_n > \veps >0$,
\begin{equation*}
\bigl| a_n^{ } - a^{(0)}_n \bigr| + \bigl|b_n^{ } - b^{(0)}_n \bigr|
\leq \tilde d_{m}\bigl( (a,b), (a^{(0)},b^{(0)}) \bigr) +
    C \sum_{r=m}^{n-p+1} \tilde d_{r}\bigl( (a,b), \calT_{J_0} \bigr)
\end{equation*}
for all $n\geq m$.  The constant $C$ depends only on $\epsilon$.
\end{lemma}

\begin{proof}
The proof is by induction on $n$.  For $m\leq n\leq m+p-1$, the
result is immediate from the definition of $\tilde d_m$.

For $n>m+p-1$, we consider the functions
\begin{align*}
f(a_1,\ldots,a_p) := \sum_{j=1}^p \bigl[ \log(a_j) - \log( a_j^{(0)} ) \bigr]
\qquad
g(b_1,\ldots,b_p) := \sum_{j=1}^p \bigl[ b_j - b_j^{(0)} \bigr]
\end{align*}
These vanish on $\calT_{J_0}$, as explained in the proof of Theorem~\ref{T2.1}.

As $g$ is Lipschitz (with constant $1$),
\begin{align*}
\bigl|b_{n} - b_{n-p} \bigr|
&= \bigl| g(b_{n},\ldots,b_{n-p+1}) - g(b_{n-1},\ldots,b_{n-p}) \bigr| \\
&\leq \bigl| g(b_{n},\ldots,b_{n-p+1}) \bigr|  + \bigl| g(b_{n-1},\ldots,b_{n-p}) \bigr| \\
&\leq \tilde d_{n-p+1}\bigl( (a,b), \calT_{J_0} \bigr)
    + \tilde d_{n-p}\bigl( (a,b), \calT_{J_0} \bigr)
\end{align*}
In a similar way,
\begin{align*}
\bigl|\log[a_n] - \log[a_{n-p}] \bigr|
&= \bigl| f(a_{n},\ldots,a_{n-p+1}) - f(a_{n-1},\ldots,a_{n-p}) \bigr|
\end{align*}
leads to
\begin{align*}
\bigl| a_n - a_{n-p} \bigr|
&\leq C_\epsilon \bigl[ \tilde d_{n-p+1}\bigl( (a,b), \calT_{J_0} \bigr)
    + \tilde d_{n-p}\bigl( (a,b), \calT_{J_0} \bigr) \bigr]
\end{align*}

Combining these two inequalities gives
\begin{align*}
\bigl| a_n^{ } - a^{(0)}_n \bigr| + \bigl|b_n^{ } - b^{(0)}_n \bigr|
&\leq \bigl| a_{n-p}^{ } - a^{(0)}_{n-p} \bigr| + \bigl|b_{n-p}^{ } - b_{n-p}^{(0)} \bigr| \\
&\qquad+ (1+C_\veps) \bigl[ \tilde d_{n-p+1}\bigl( (a,b),
\calT_{J_0} \bigr)
    + \tilde d_{n-p}\bigl( (a,b), \calT_{J_0} \bigr) \bigr]
\end{align*}
which completes the proof of the inductive step.
\end{proof}

\begin{proof}[Proof of Proposition~\ref{P2.5}]  The left-hand inequality in \eqref{2.13} follows immediately
from the definitions of $d_m$ and $\tilde d_m$; we focus on the second inequality.

Choose $(a^{(0)},b^{(0)})$ minimizing $d_m((a,b),\calT_{J_0})$; strictly, this amounts to a (inconsequential) change
in $J_0$.  Applying Lemma~\ref{L2.6} in the definition of $d_m$ gives
\begin{align*}
d_m\bigl((a,b),\calT_{J_0}\bigr) &\leq \tfrac{e}{e-1}\, \tilde d_m\bigl((a,b),\calT_{J_0}\bigr)
    + C \sum_{k=0}^\infty \sum_{r=m}^{m+k}  e^{-k} \, \tilde d_{r} \bigl((a,b),\calT_{J_0}\bigr) \\
&\leq C' \sum_{j=0}^\infty e^{-j} \, \tilde d_{m+j} \bigl((a,b),\calT_{J_0}\bigr)
\end{align*}
The proposition follows because convolution with $e^{-j}\chi_{[0,\infty)}(j)$ is a bounded operator on all $\ell^q$ spaces.
\end{proof}

\section{The Magic Formula for CMV Matrices} \lb{s3}

Our goal in this section is to prove

\begin{theorem}\lb{T3.1} Let $p$ be even and let $\calC_0$ be a two-sided
$p$-periodic CMV matrix with discriminant $\Delta_{C_0}(z)$ and
isospectral torus $\calT_{\calC_0}$. Given a two-sided {\rm{(}}not
a priori periodic{\rm{)}} CMV matrix, $\calC$,
\begin{equation} \lb{3.1}
\Delta_{\calC_0}(\calC) =S^p + S^{-p} \iff \calC\in\calT_{\calC_0}
\end{equation}
\end{theorem}

\begin{remarks} 1. Notice that since $\calC$ is unitary and $\Delta (e^{i\theta})$
is real, $\Delta_{\calC_0}(\calC)$ is self-adjoint.

2. By \eqref{1.13a} and the fact that $\calC$ is five-diagonal, $\Delta_{\calC_0}
(\calC)$ has $2(p/2)$ diagonals above/below the main diagonal.

3. As in Section~\ref{s2}, we will first present our initial proof that
\begin{equation} \lb{3.1a}
\Delta_{\calC_0}(\calC)=S^p + S^{-p} \Rightarrow \{\alpha_n\} \text{ is periodic}
\end{equation}
and then a proof based on Golinskii's suggestion.
\end{remarks}

\begin{lemma}\lb{L3.2} We have:
\begin{align}
&(\calC^\ell)_{m,m+k}= (\calC^{-\ell})_{m,m+k} = 0 \qquad \text{ if }k>2\ell \lb{3.2} \\
&(\calC^\ell)_{2m,2m+2\ell} = \rho_{2m} \rho_{2m+1} \dots \rho_{2m+2\ell-1} \lb{3.3} \\
&(\calC^\ell)_{2m+1, 2m+2\ell+1} = 0 \lb{3.4} \\
&(\calC^{-\ell})_{2m, 2m+2\ell} = 0 \lb{3.5} \\
&(\calC^{-\ell})_{2m+1, 2m+\ell +1} = \rho_{2m+1} \rho_{2m+2} \dots \rho_{2m+2\ell} \lb{3.6} \\
&(\calC^\ell)_{2m, 2m+2\ell -1} = \rho_{2m} \rho_{2m+1} \dots \rho_{2m+2\ell-2}
\bar\alpha_{2m+2\ell-1} \lb{3.7} \\
&(\calC^\ell)_{2m+1, 2m+2\ell} = -\alpha_{2m} \rho_{2m+1} \dots \rho_{2m+2\ell-1} \lb{3.8} \\
&(\calC^{-\ell})_{2m, 2m+2\ell-1} = -\bar\alpha_{2m-1} \rho_{2m} \rho_{2m+1} \dots
\rho_{2m+2\ell-2} \lb{3.9} \\
&(\calC^{-\ell})_{2m+1, 2m+2\ell} = \rho_{2m+1} \dots \rho_{2m+2\ell-1}
\alpha_{2m+2\ell} \lb{3.10}
\end{align}
\end{lemma}

\begin{proof} As $\calL$ and $\calM$ are tridagonal,
$\calC^\ell$ is a product of $2\ell$ tridiagonal matrices, so \eqref{3.2} is immediate.

We will prove the results for $\calC^\ell$. The results for $\calC^{-\ell}$ are similar
if we note
\begin{equation} \lb{3.15}
\Theta (\alpha)^{-1} = \Theta (\bar\alpha)
\end{equation}
since $\Theta$ is unitary and symmetric.

Equation \eqref{3.3} follows from
\begin{equation} \lb{3.16}
\calL_{2m, 2m+1} = \rho_{2m} \qquad
\calM_{2m+1, 2m+2} = \rho_{2m+1}
\end{equation}
and \eqref{3.4} from
\begin{equation} \lb{3.17}
\calL_{2m+1, 2m+2} =0
\end{equation}
Because of \eqref{3.17}, the only way for $\calC^\ell$ to get from $2m$ to $2m+2\ell-1$ is to
increase index in the first $2\ell-1$ factors, which leads to \eqref{3.7}. For the same reason,
to get from $2m+1$ to $2m+2\ell$, the last $2\ell -1$ factor must increase index, leading to
\eqref{3.8}.
\end{proof}

\begin{lemma}\lb{L3.3} If $\calC$ and $\calC_0$ are $p$-periodic, then
$\sigma (\calC)=\sigma (\calC_0)$ if and only if $\Delta_\calC=\Delta_{\calC_0}$.
\end{lemma}

\begin{proof} Either proof of Lemma~\ref{L2.4} carries over with no change.
\end{proof}

\begin{proof}[Proof of Theorem~\ref{T3.1}] The proof that
\begin{equation} \lb{3.18}
\Delta_{\calC_0} (\calC_0) =S^p + S^{-p}
\end{equation}
is identical to the proof of \eqref{2.5}.

For the converse, suppose
\begin{equation} \lb{3.19}
\Delta_{\calC_0}(\calC) = S^p + S^{-p}
\end{equation}
In particular,
\begin{equation} \lb{3.20}
\Delta_{\calC_0} (\calC)_{2m,2m+p-1} =0
\end{equation}

By \eqref{1.13a} and Lemma~\ref{L3.2}, this implies (recall $p$ is even)
\[
(\rho_0^{(0)}\dots\rho_{p-1}^{(0)})^{-1} (\rho_{2m} \rho_{m+1} \dots
\rho_{2m + p-2}) (\bar\alpha_{2m+p-1} - \bar\alpha_{2m-1}) = 0
\]
so
\begin{equation} \lb{3.21}
\alpha_{2m+p-1} = \alpha_{2m-1}
\end{equation}

Similarly, since
\[
\Delta_{\calC_0} (\calC)_{2m+1, 2m+p} =0
\]
we get
\[
(\rho_0^{(0)}\dots\rho_{p-1}^{(0)})^{-1} (\rho_{2m+1} \dots \rho_{2m+p-1})
(\alpha_{2m+2}-\alpha_{2m+p-1}) =0
\]
which leads to
\begin{equation} \lb{3.22}
\alpha_{2m+p} = \alpha_{2m}
\end{equation}

Thus, $\alpha$ has period $p$.  That $\calC\in\calT_{\calC_0}$ follows from
Lemma~\ref{L3.3} and the same argument used in the OPRL case.
\end{proof}

Next, we give a proof using Na\u\i man's lemma. We will need

\begin{lemma}\lb{L3.3A} Let $\calC$ be the extended CMV matrix associated to
$\{\alpha_n\}_{n=-\infty}^\infty$. Let $p$ be even. If $S^p\calC =\calC S^p$,
then
\begin{equation} \lb{3.22a}
\alpha_{n+p} =\alpha_n
\end{equation}
for all $n$.
\end{lemma}

\begin{proof} We have that $\calC_{2j\,\, 2j+1}^2 + \calC_{2j\,\, 2j+2}^2 =
\rho_{2j}^2$ (see (4.2.14) of \cite{OPUC1}), so $\rho_{2j}$ is periodic. Thus,
$\calC_{2j\,\, 2j+2}/\rho_{2j}=\rho_{2j+1}$ is also periodic. So $\bar\alpha_{2j+1}
=\calC_{2j\,\, 2j+1}/\rho_{2j}$ is periodic as is $\alpha_{2j} =\calC_{2j+1\,\,2j+2}/
(-\rho_{2j+1})$.
\end{proof}

\begin{proof}[Second proof that \eqref{3.1a} holds] $\calC$ commutes with $S^p + S^{-p}$,
so by Na\u\i man's lemma (Lemma~\ref{L2.4A}), which did not require that $A$ be self-adjoint,
$S^p\calC=\calC S^p$, which implies $\alpha$ is periodic by Lemma~\ref{L3.3A}.
\end{proof}

%

We now turn to the OPUC version of Proposition~\ref{P2.5}.  As noted in the introduction, it is not
sufficient to sum over exactly one period:

\begin{example}\label{E4.5} $(0,\f12, 0, \f12, 0, \dots)$ and $(0,-\f12, 0, -\f12, 0, \dots)$ are
in the same isospectral torus, namely, the one with $p=2$ and
\[
\Delta(z) = \sqrt{\tfrac43}\, (z+z^{-1})
\]
Now consider $\alpha =(0,\f12, 0, -\f12, 0, \f12, 0, -\f12, \dots)$. If $\ti d_m (\alpha,
\calT_{\calC_0})$ were defined as sum from $k=0$ to $p-1$, it would be zero for all $m$,
but $d_m (\alpha, \calT_{\calC_0})$ is not small.
\qed
\end{example}

The problem, as this example shows, is that for sequences in $\calT_{\calC_0}$,
$(\alpha_0,\dots, \alpha_{p-2})$ does not determine $\alpha_{p-1}$. But by periodicity,
$\alpha_0, \dots, \alpha_{p-1}$ determines $\alpha_p$.  Thus, if we define
\begin{equation} \lb{3.23}
\ti d_m (\alpha,\alpha'):=\sum_{k=0}^p\, \abs{\alpha_{m+k}-\alpha'_{m+k}}
\end{equation}
then
$$
| \alpha_{m+p} - \alpha_m | \leq \tilde d_m(\alpha,\calT_{\calC_0})
$$
Plugging this into the proofs of Lemma~\ref{L2.6} and Proposition~\ref{P2.5} leads quickly to

\begin{proposition}\lb{P3.5} Let $\calC_0$ be a fixed periodic CMV matrix, then
\begin{equation} \lb{3.24}
e^{-p} \bigl\| \tilde d_m (\alpha,\calT_{\calC_0}) \bigl\|_{\ell^q}
\leq \bigl\| d_m (\alpha,\calT_{\calC_0}) \bigl\|_{\ell^q}
    \leq C \bigl\| \tilde d_m (\alpha,\calT_{\calC_0}) \bigl\|_{\ell^q}
\end{equation}
for any sequence of Verblunsky coefficients $\{\alpha_n\}$.
\end{proposition}

\section{The Magic Formula for Schr\"odinger Operators} \lb{s4}

In this section, we want to illustrate the potential applicability of our
central idea to the spectral theory of one-dimensional Schr\"odinger
operators,
\begin{equation}\label{SO}
 H = -\f{d^2}{dx} + V(x)
\end{equation}
However, we will not pursue the applications in this paper.

We will suppose $V\in L_\unif^1$, that is, $\sup_x \int_{x-1}^{x+1}
\abs{V(y)}\, dy <\infty$. In that case, $V$ is a form bounded perturbation
of $-\f{d^2}{dx^2}$ on $L^2 (\bbR, dx)$ with relative bound zero, so
$H$ is a self-adjoint operator.  Its form domain is the Sobolev space $H^1(\bbR)$.

We need to say something about periodic Schr\"odinger operators. Suppose $V_0$
has period $L$, that is,
\begin{equation} \lb{4.1}
V_0 (x+L) = V_0(x)
\end{equation}

For arbitrary $V$ in $L_\loc^1$ and $E\in\bbC$, let $u_D (x,E;V)$ and $u_N
(x,E;V)$ (we will often drop the $V$ if it is clear which $V$ is intended) be
the solutions of
\begin{equation} \lb{4.2}
-u'' + Vu =Eu
\end{equation}
obeying the boundary conditions
\begin{equation} \lb{4.3}
u_D(0)=0 \qquad u'_D(0)=1 \qquad u_N(0)=1 \qquad u'_N(0)=0
\end{equation}
There is a unique solution of \eqref{4.2} in distributional sense which is
absolutely continuous.

The transfer matrix that updates solutions of \eqref{4.2} (with data written as
$\binom{u}{u'}$) is
\begin{equation} \lb{4.4}
T(x,E;V) = \begin{pmatrix}
u_N(x,E) & u_D(x,E) \\
u'_N(x,E) & u'_D(x,E)
\end{pmatrix}
\end{equation}
$\det(T)=1$ by constancy of the Wronskian. For periodic $V_0$, we define the
discriminant by
\begin{align}
\Delta_{V_0}(E) &= \tr (T(L,E; V_0)) \notag \\
&= u_N (L,E) + u'_D(L,E) \lb{4.5}
\end{align}

As in the OPRL and OPUC cases, it is easy to see for the whole-line
operator that
\begin{equation} \lb{4.6}
\sigma \bigl( -\tfrac{d^2}{dx^2} + V_0 \bigr) = \Delta_{V_0}^{-1} ([-2,2])
\end{equation}
and is purely absolutely continuous. Moreover (see, e.g., \cite{RS4}), if
\begin{equation} \lb{4.7}
(S_y u)(x) = u(x-y)
\end{equation}
then $H=-\f{d^2}{dx^2}+V_0$ commutes with $S_L$ and so has a direct integral decomposition,
\begin{equation} \lb{4.8}
H=\int^\oplus H(\theta)\, \f{d\theta}{2\pi}
\end{equation}
whose fibers, $H(\theta)$, are the operator \eqref{SO} on $[0,L]$ with
\begin{equation} \lb{4.9}
u(L)=e^{i\theta} u(0) \qquad u'(L) =e^{i\theta} u'(0)
\end{equation}
boundary conditions. $H(\theta)$ has purely discrete spectrum (i.e., $(H(\theta)
+i)^{-1}$ is compact); the eigenvalues are precisely the solutions of
\begin{equation} \lb{4.10}
\Delta (E) =2\cos (\theta)
\end{equation}

Two periodic potentials of period $L$ are called isospectral if and only if they have
the same $\Delta$.  As in the Jacobi and CMV cases, the spectrum determines $\Delta$, but
this is more difficult to prove in the Schr\"odinger case.
It is also known (\cite{McKT} for nice $V_0$'s; \cite{BuFi,FIT,Iwa,RT})
that the set of $V$'s isospectral to $V_0$ is a torus of dimension equal to the
number of gaps which is typically infinite, so we will refer to an isospectral
torus, $\calT_{V_0}$. We can now state the main result in this section:

\begin{theorem}\lb{T4.1} Let $V_0$ be periodic obeying \eqref{4.1} and let $\Delta_{V_0}$
be its discriminant and $\calT_{V_0}$ its isospectral torus.
Let $V$ be in $L_{\loc,\unif}^1$ on $\bbR$ and $H=-\f{d^2}{dx^2} +V$\!. Then
\begin{equation} \lb{4.11}
\Delta_{V_0}(H) = S_L + S_{-L}\iff V\in\calT_{V_0}
\end{equation}
Here $S_{\pm L}$ denotes the shift operator, as in \eqref{4.7}.
\end{theorem}

\begin{remarks} 1. $\Delta_{V_0}(H)$ is defined by the functional calculus.

2. As in the last two sections, we will provide our initial proof that
\begin{equation}\lb{4.11a}
\Delta_{V_0}(H) = S_L + S_{-L}\Rightarrow V \text{ periodic}
\end{equation}
and then a simpler proof using an analog of Na\u\i man's lemma. This argument does not
require Theorems~\ref{T4.2} and \ref{T4.3} and the considerable machinery their
proofs entail. That said, to show $\Delta_{V_0}(H)=S_L + S_{-L}$ plus $V$ periodic
implies $V\in\calT_{V_0}$ {\it does\/} require Theorem~\ref{T4.3}, but it should be
noted that one can prove Theorem~\ref{T4.3} fairly easily without needing transformation
formulae of Delsarte, Levitan, Gel'fand, Marchenko type.
\end{remarks}

We need two preliminaries whose proofs we defer to later in the section. We first
make a definition:

\begin{definition} For any $y>0$, $\calR_y$ consists of operators on $L^2 (\bbR)$
of the form
\begin{equation} \lb{4.12}
(Af)(x) = \tfrac12\, f(x+y) + \tfrac12\, f(x-y) + \int_{x-y}^{x+y}
K(x,z)f(z)\, dz
\end{equation}
where $K$ is continuous and uniformly bounded on $\{(x,z) : \abs{x-z} \leq y\}$.
\end{definition}

{\it Note.} It can happen that $K(x,x\pm y)\neq 0$, so if we think of $K$ as an integral
kernel on $\bbR\times\bbR$, it can be discontinuous at $\abs{x-z}=y$.

\begin{theorem} \lb{T4.2} If $V_0$ is $L$-periodic and $V$ in $L_{\loc,\unif}^1$,
then $\f12 \Delta(H)\in\calR_L$ and
\begin{equation} \lb{4.13}
K(x,x+L)=-\tfrac14 \int_x^{x+L} (V(z) - V_0(z))\, dz
\end{equation}
\end{theorem}

Note that \eqref{4.13} describes the `matrix elements' of $\Delta_{V_0} (H)-(S_L + S_{-L})$
that are farthest from the diagonal.  Indeed, just as in the other cases, one does not need the full statement
$\Delta_{V_0} (H)=S_L + S_{-L}$ to see that $V$ is periodic, only that $\langle f, (\Delta (H)-S_L - S_{-L})g\rangle
=0$ for $f$ supported near $x_0$ and $g$ near $x_0 +L$ (for all $x_0$).

\begin{theorem}\lb{T4.3} $\Delta(E)$ is an entire function which obeys
\begin{alignat}{2}
& \text{\rm{(i)}} \qquad && \abs{\Delta(E)} \leq C\exp \bigl( L\sqrt{\abs{E}}\bigr) \lb{4.14} \\
& \text{\rm{(ii)}} \qquad && \lim_{\substack{E\to -\infty \\ E\text{ real}}}
\, \f{\Delta(E)}{\exp (L\sqrt{\abs{E}})}=1 \lb{4.15}
\end{alignat}
\end{theorem}

\begin{proof}[Proof of Theorem~\ref{T4.1}] If $V\in\calT_{V_0}$, then $\Delta_V =
\Delta_{V_0}$, so for the $\Leftarrow$ direction we need only prove
\begin{equation} 
\Delta_{V_0}\bigl( \tfrac{d^2}{dx^2} + V_0) = S_L + S_{-L}
\end{equation}
As before, this is equivalent to $\Delta_{V_0}(H(\theta))=2\cos\theta$ which follows
from \eqref{4.10}.

Conversely, if $\Delta_{V_0}(H)=S_L + S_{-L}$, then from Theorem~\ref{T4.2} and the periodicity
of $V_0$, we see
\begin{equation} \lb{4.16}
\int_x^{x+L} V(z)\, dz = \text{constant}
\end{equation}
This implies that $V(x+L)-V(x)=0$ for a.e.\ $x$, that is, $V$ is periodic.

If $H(\theta)$ are the fibers of $H$ in the direct integral decomposition,
$\Delta_{V_0}(H)=S_L + S_{-L}$ implies
\begin{equation} \lb{4.17}
\Delta_{V_0}(H(\theta))=2\cos\theta
\end{equation}
so, if $\Delta$ is the discriminant for $V$\!, we have $\Delta
(z)=\pm 2\Rightarrow \Delta_{V_0}(z) =\pm 2$. Moreover, \eqref{4.17} implies $\sigma(H)
\subseteq \sigma (-\f{d^2}{dx^2} + V_0)$, so any double zero of $\Delta\pm 2$ is
a double zero of $\Delta_{V_0}\pm 2$. It follows that
\begin{equation} \lb{4.18}
g(z) = \f{\Delta_{V_0}^2 (z) -4}{\Delta^2(z) -4}
\end{equation}
is analytic.

Since $\Delta_{V_0}$ and $\Delta$ are entire functions of order $\f12$ (by Theorem~\ref{T4.3}),
$g(z)$ is of the form
\begin{equation} \lb{4.19}
g(z) = C\prod_{j=1}^J \biggl( 1-\f{z}{z_j}\biggr)
\end{equation}
where $z_1 < z_2 < \cdots$ are bounded from below. By \eqref{4.15}, $\lim_{E\to -\infty}
g(E)=1$, which implies $g\equiv 1$, that is, $\Delta = \Delta_{V_0}$.
\end{proof}

The argument used at the end of the proof to conclude that missing zeros cannot
occur is reminiscent of ideas connected with the Hochstadt--Lieberman \cite{HL}
and related theorems \cite{GShpa,GStams}.

We now turn to the proofs of Theorems~\ref{T4.2} and \ref{T4.3}. A critical role
is played by the wave equation and the transformation operator formalism of Gel'fand--Levitan,
further important work is due to Delsarte, Levin, and Marchenko;
see the book of Marchenko \cite{Marc} for references and history.

Define for $s>0$,
\begin{equation} \lb{4.20}
C_s(z) =\cos \bigl(s\sqrt{z}\bigr) \qquad
S_s(z) = z^{-1/2} \sin \bigl(s\sqrt{z}\bigr)
\end{equation}
which are entire functions of $z$ bounded on $(a,\infty)$ for any $a\in\bbR$. Thus
$C_s(H)$ and $S_s(H)$ are bounded operators for any $H$ that is bounded from below. We will
need to study the form of $C_s (-\f{d^2}{dx^2}+V)$. For bounded continuous $V$\!,
this is discussed in Marchenko \cite{Marc}. While his proofs extend to the
$L_\loc^1$ case, it seems simpler to sketch the ideas:

\begin{proposition}\lb{P4.4} $C_s^0:=C_s (-\f{d^2}{dx^2})\in\calR_s$; indeed,
\begin{equation} \lb{4.21}
(C_s^0 f)(x) =\tfrac12\, [f(x+s) + f(x-s)]
\end{equation}
If $S_s^0:=S_s (-\f{d^2}{dx^2})$, then
\begin{equation} \lb{4.22}
(S_s^0 f)(x) = \tfrac12 \int_{x-s}^{x+s} f(y)\, dy
\end{equation}
\end{proposition}

\begin{remark} If $w(x,s):=(C_s^0 f)(x)+(S_sg)(x)$, then $w$ obeys the wave equation
$(\f{\partial^2}{\partial s^2} -\f{\partial^2}{\partial x^2}) w=0$ with
initial data $w(x,0)=f$ and $\partial_s w(x,0)=g(x)$.  Thus the proposition basically encodes d'Alembert's
solution of the wave equation. From this point of view, Theorem~\ref{T4.2} is connected to
finite propagation speed for the wave equation.
\end{remark}

\begin{proof} Since $\cos$ is even,
\begin{equation} \lb{4.23}
\cos (s\abs{k})=\cos (sk) = \tfrac12\, (e^{iks} + e^{-iks})
\end{equation}
\eqref{4.21} is just the Fourier transform of this. \eqref{4.22} follows from
\begin{equation} \lb{4.24}
S_s(z) = \int_0^s C_t(z)\, dt
\end{equation}
and \eqref{4.21}.
\end{proof}

We are heading towards

\begin{theorem}\lb{T4.5} Let $V\in L_{\loc,\unif}^1 (\bbR)$ and let
$H=-\f{d^2}{dx^2}+V$\!. Then $C_s(H)\in\calR_s$ and the associated kernel
$K_s$ of \eqref{4.12} obeys
\begin{equation} \lb{4.25}
K_s (x,x+s) = -\tfrac14\int_x^{x+s} V(u)\, du
\end{equation}
and for each $t\in (0,\infty)$,
\begin{equation} \lb{4.26}
\sup_{x,y, \abs{s}\leq t}\, \abs{K_s (x,y)} <\infty
\end{equation}

In addition,
\begin{equation} \lb{4.27}
(S_s(H)f)(x) = \int_{x-s}^{x+s} L_s (x,y) f(y)\, dy
\end{equation}
where
\begin{equation} \lb{4.28}
L_s (x,x+s) =\tfrac12
\end{equation}
\end{theorem}

\begin{lemma}\lb{L4.6} It suffices to prove Theorem~\ref{T4.5} for $s$ small.
\end{lemma}

\begin{proof} Since $\cos (2u)=2\cos^2 (u)-1$, one sees
\begin{equation} \lb{4.29}
C_{2s}(A) = 2C_s(A)^2 -\bdone
\end{equation}
Thus, if $C_s\in\calR_s$, one sees $C_{2s}\in\calR_{2s}$ and
\begin{equation} \lb{4.30}
K_{2s}(x,y) = K_s (x,y+s) + K_s (x,y-s) + K_s (x+s,y) +K_s (x-s,y)
\end{equation}
where $K(x,y)=0$ if $\abs{x-y}>s$. Thus
\begin{equation} \lb{4.31}
K_{2s} (x,x+2s) = K_s (x,y+s) + K_s (x+s, y+2s)
\end{equation}

This shows that if the formula is known for $\abs{s}\leq T$, one gets
it successively for $2T, 4T, 8T, \dots$.

Using \eqref{4.24}, one sees that the result for $C_s(H)$ implies \eqref{4.27}
and \eqref{4.28}.
\end{proof}

\begin{proof}[Proof of Theorem~\ref{T4.5}] If $A$ is a bounded self-adjoint
operator on $\calH$ which is bounded from below, and $B$ is the operator on
$\calH\oplus\calH$ given by
\begin{equation} \lb{4.32}
B = \left(
\begin{array}{rc}
\bdzero & \bdone \\
-A & \bdzero
\end{array}\right)
\end{equation}
then
\begin{equation} \lb{4.33}
e^{sB} = \left(
\begin{array}{rc}
C_t(A) & S_t(A) \\
-AS_t(A) & C_t(A)
\end{array}\right)
\end{equation}

This formula can be checked by showing that the right side of \eqref{4.33} is
a bounded semigroup whose derivative at $t=0$ is $B$. DuHamel's formula for
$A,\ti A$ bounded says that
\begin{align}
e^{t \ti B} &= e^{tB} + \int_0^t  e^{sB} (\ti B-B) e^{(t-s)\ti B}\, ds \lb{4.34} \\
&= e^{tB} + \int_0^t e^{s\ti B} (\ti B-B) e^{(t-s)B}\, ds \lb{4.35}
\end{align}
Using \eqref{4.33}, we obtain
\begin{align}
C_t(\ti A) &= C_t(A) - \int_0^t S_s (\ti A) (\ti A-A) C_{t-s} (A)\, ds \lb{4.36} \\
&= C_t(A) - \int_0^t S_s (A) (\ti A-A) C_{t-s} (\ti A)\, ds \lb{4.37}
\end{align}

By taking limits, it is easy to obtain these formulae for $A=-\f{d^2}{dx^2}$,
$\ti A=-\f{d^2}{dx^2} +V$ with $V$ bounded. By obtaining a priori bounds below
depending only on certain $L^1$ norms of $V$\!, we get estimates for $V$ in
$L^1$ and so, using the lemma, prove the theorem.

By iterating \eqref{4.37}, one gets an expansion (which converges if $V$ is
bounded and whose estimates then extend),
\begin{gather}
C_t \bigl( -\tfrac{d^2}{dx^2} + V(x) \bigr) = C_t^{(0)} +
\sum_{n=1}^\infty C_t^{(n)} \lb{4.38} \\
C_t^{(n)} = (-1)^n \int_{0\leq s_1 + \cdots + s_n\leq t} S_{s_1}^{(0)} V
S_{s_2}^{(0)} \dots V S_{s_n}^{(0)} V C_{t-s_1- \cdots -s_n}^{(0)}\, ds_1, \dots ds_n \lb{4.39}
\end{gather}

Apply the integrand in $C_t^{(n)}$ to a function $f$ and evaluate at $x$ for
fixed $s_1, \dots, s_n$. Each $S_{s_j}^{(0)}V$ evaluates $V$ at points and integrals
using \eqref{4.22}. The integrands in $V$ are in the interval $(x-t, x+t)$, so
if we take absolute values, we see this integrand is bounded by
\[
\biggl(\tfrac12 \int_{x-t}^{x+t} \abs{V(y)}\, dy \biggr)^n
[\tfrac12\, f(x+t-s_1 - \cdots -s_n) + \tfrac12\, f(x-t+s_1 + \cdots + s_n)]
\]

Now we can do the integral over $s_1, \dots, s_n$. For $t-s_1 -\cdots - s_n$ fixed,
the new integrand is independent of $s_1, \dots, s_{n-1}$ and is bounded by $t^{n-1}$.
We find
\begin{equation} \lb{4.40}
\abs{(C_t^{(n)}f)(x)}\leq t^{n-1} \biggl( \tfrac12\int_{x-t}^{x+t} \abs{V(y)}\, dy
\biggr)^n \int_{x-t}^{x+t} \abs{f(y)}\, dy
\end{equation}

Moreover, $C_t^{(n)}$ has a continuous integral kernel $K_t^{(n)}(x,y)$ supported
in $\abs{x-y}\leq t$. Since $V$ is uniformly locally $L^1$, by taking $t$ small, we
can be sure $\sup_x \f12 \int_{x-t}^{x+t} \abs{V(y)}\, dy <1$, which yields
uniform convergence of $K_t^{(n)}$ to a uniformly bounded kernel.

By \eqref{4.24}, we get \eqref{4.27} from $C_s(H)\in\calR_s$, and \eqref{4.28}
comes from noting that
\begin{equation} \lb{4.41}
\abs{L_s (x,y)-\tfrac12}\leq (s-\abs{x-y})\, \sup_{x,y,u\leq s}\,
\abs{K_u (x,y)}
\end{equation}

Finally, using \eqref{4.36}, we see that
\[
K_t (x,x+t) = -\tfrac12 \int_0^t L_s (x ,x+s) V(x+s)\, ds
\]
proving \eqref{4.25}.
\end{proof}

To complete the proofs of Theorems~\ref{T4.2} and \ref{T4.3} (and so Theorem~\ref{T4.1}),
we need the transformation formulae of Delsarte, Levitan, Gel'fand, and Marchenko
\cite{Marc}:

\begin{theorem}\lb{T4.7} If $V\in L^1 ([0,R])$ for $R<\infty$, then there exist
functions $K_N,K_D$ $C^1$ in $\{(y,x) : 0\leq y\leq  x\leq R\}$ so that for
$0\leq x\leq R$,
\begin{align}
u_N(x,E) &=C_x(E) + \int_0^x K_N (y,x) C_y (E)\, dy  \lb{4.42} \\
u_D (x,E) &= S_x (E) + \int_0^x K_D (y,x) S_y (E)\, dy \lb{4.43}
\end{align}
Moreover,
\begin{equation} \lb{4.44}
K_D(x,x) = K_N(x,x) = \tfrac12 \int_0^x V(t)\, dt
\end{equation}
\end{theorem}

\begin{remarks} 1. These formulae are in Marchenko \cite[p.~9 and (1.2.28)]{Marc}.
He supposes $V$ is continuous, but his proof works if $V$ is $L^1$; indeed, see
Remark~2.

2. Defining $\ti u_X (x,k) = u_X (x,k^2)$ for $X=D,N$ and $Q_X (x,y)$ as the
Fourier transform of $\ti u_X$ in $k$, we see \eqref{4.2} becomes
\begin{equation} \lb{4.45}
\f{\partial^2 Q}{\partial x^2} - \f{\partial^2 Q}{\partial y^2} = VQ(x,y)
\end{equation}
with initial conditions
\begin{alignat*}{2}
Q_N(x=0,y)&=\delta (y) \qquad & Q_n' (x=0,y) &=0 \\
Q_D(x=0,y) &= 0 \qquad & Q'_D (x=0,y) &= \delta(y)
\end{alignat*}
Thus, Theorem~\ref{T4.7} is essentially Theorem~\ref{T4.5} with a time-dependent
$V$ used.

3. By \eqref{4.5}, \eqref{4.42}, and \eqref{4.43}, we obtain a critical representation
for $\Delta$:
\begin{equation} \lb{4.46}
\begin{split}
\Delta(E) = 2C_x (E) &+ \int_0^L L_1 (t) C_t (E)\, dt \\
& \qquad + \int_0^L L_2 (t) S_t(E)\, dt + K_D (L,L) S_L(E)
\end{split}
\end{equation}
where $L_1,L_2$ are continuous in $[0,L]$. Indeed,
\[
L_1(t) = K_N (t,L) \qquad L_2 (t) = \left.\f{\partial}{\partial  x}\, K_D (t,x)\right|_{x=L}
\]
\end{remarks}

\begin{proof}[Proof of Theorem~\ref{T4.3}] The analyticity is immediate from \eqref{4.46}
as is \eqref{4.14} given
\[
\abs{C_x (E)} + \abs{S_x(E)} \leq C\exp \bigl( x\sqrt{\abs{E}}\bigr)
\]
Moreover, since for $t<L$,
\[
\lim_{E\to -\infty}\, \f{C_t(E)}{C_L(E)} =0 \qquad\text{and}\qquad
\lim_{E\to -\infty}\, \f{S_L(E)}{C_L(E)} =0
\]
we have \eqref{4.15}.
\end{proof}

\begin{proof}[Proof of Theorem~\ref{T4.2}] By \eqref{4.46} and Theorem~\ref{T4.5},
$\f12 \Delta (-\f{d^2}{dx^2}+V)$ is in $\calR_L$. Moreover, the only terms contributing
to $K(x,x+L)$ come from $C_L (-\f{d^2}{dx^2}+V)$ and $K_D (L,L) S_L(-\f{d^2}{dx^2}+V)$.
By \eqref{4.25}, \eqref{4.28}, and \eqref{4.44},
\[
K(x,x+L) = -\tfrac14 \int_x^{x+L} V(y)\, dy + \tfrac12 \biggl( \tfrac12 \int_0^L
V_0(y)\, dy \biggr)
\]
which, given the periodicity of $V_0$, is \eqref{4.13}.
\end{proof}

There is a second proof of \eqref{4.11a}. It depends on this analog of
Na\u\i man's lemma:

\begin{lemma} \lb{L4.8} If $V$ is $L_{\loc,\unif}^1$ and $-\f{d^2}{dx^2}+V$
commutes with $S_L+S_{-L}$, then
\begin{equation}\lb{4.48a}
V(x+L)=V(x)
\end{equation}
\end{lemma}

\begin{proof} Suppose first that $V$ is bounded. Then $S_L + S_{-L}$ leaves
$D(-\f{d^2}{dx^2})$ invariant and commutes with it, so $S_L + S_{-L}$ commutes
with $V$\!. If $f$ is supported in a small neighborhood of $x_0$, $(x_0-\delta,
x_0 +\delta)$ with $\abs{\delta} < L/2$, then $(S_L+S_{-L})(Vf)$ is two separate
pieces $V(x-L) f(x-L)$ supported near $x_0+L$ and $V(x+L)f(x+L)$ supported near
$x_0-L$, while $V(S_L+S_{-L})f$ is two pieces $V(x) f(x-L)$ and $V(x)f(x+L)$.
Since the pieces are disjoint,
\[
V(x) f(x-L) =V(x+L)f(x-L)
\]
which implies \eqref{4.48a}.

For general $V$\!, take $g\in \calC_0^\infty(\bbR)$ with $\int g(x)\, dx =1$
and note that $\int g(x) S_x (-\f{d^2}{dx^2}+V) S_{-x}\, dx$ is $-\f{d^2}{dx^2}
+g*V$ and it commutes with $S_L+S_{-L}$ also. But $g*V$ is bounded, so it is periodic.
\eqref{4.48a} follows by using an approximate $\delta$-function.
\end{proof}

\section{Block Jacobi Matrices and Matrix Orthogonal Polynomials} \lb{s5}

What the magic formula suggests is that the Jacobi matrix $J$ has parameters that approach
an isospectral torus if and only if $\Delta(J)$ approaches $S^p + S^{-p}$. $\Delta(J)$
is a matrix of width $2p+1$ (i.e., $\Delta(J)_{k\ell}=0$ if $k-\ell\notin \{0,\pm 1, \dots,
\pm p\}$) and $S^p + S^{-p}$ is a matrix with $1$'s at the extremes.

A matrix of width $2\ell+1$ has the structure of a tridiagonal matrix if rewritten in terms
of $\ell\times\ell$ blocks and $S^\ell + S^{-\ell}$ corresponds to $B_k=\bdzero$, $A_k=\bdone$,
the identity matrix, so $\Delta(J)\sim S^\ell + S^{-\ell}$, at the matrix level,
approaches the `free case.'  This will allow us to reduce our main theorems to
matrix analogs of the theorems on perturbations of the free case.

Of course, the association of the block matrix to orthogonal polynomials is critical---the
orthogonality will be with respect to a matrix-valued measure. There is a huge literature
on MOPRL (see, e.g.,
\cite{CFMV03,DJLP,DS2002,DS02,Dur99,DD02,DLR96,DLR01,DLR04,DLRS99,JDP96,MS93,MY03,YM2001,YMP01})
and MOPUC (see, e.g., \cite{BakCon,DG92,DGK78,DGK3,DGK2,DGK79,DGK81,Ger81,Lev,MRG89,Rod90,newVA,YK}).
In this section, our main purpose is to set notation and discuss the important notion
of equivalent families of block Jacobi matrices, a notion discussed more explicitly
in \cite{DPSprep}.

Given a semi-infinite complex matrix $M=\{m_{ij}\}_{1\leq i,j\leq\infty}$ and
$\ell=1,2,\dots$, we define the $\ell\times\ell$ block decomposition as the family
of $\ell\times\ell$ matrices $\{M_{qr}\}_{1\leq q,r\leq \infty}$ by
\begin{equation} \lb{5.1}
(M_{qr})_{ij} = m_{\ell(q-1)+i,\, \ell(r-1)+j} \qquad i,j=1,\dots\ell
\end{equation}

\begin{definition} A block Jacobi matrix is an $M$ where
\begin{equation} \lb{5.2}
M_{qr} = \begin{cases}
B_q &\text{if } r=q\geq 1 \\
A_q & \text{if } r=q+1,\, q\geq 1 \\
A_{q-1}^\dag &\text{if } r=q-1, \, q\geq 2 \\
0 & \abs{q-r}\geq 2
\end{cases}
\end{equation}
with each $A_q$ invertible and each $B_q$ Hermitian.  Here, following \cite[Section~2.13]{OPUC1},
we use $^\dagger$ for Hermitian adjoint; this is to avoid confusion with the Szeg\H{o} dual $\Phi_n^*$ appearing in OPUC.
\end{definition}

We will start writing $\calJ$ for such matrices.

In analogy with the scalar case, one may be tempted to require
\begin{equation} \lb{5.3}
A_q >0
\end{equation}
but to include $\Delta(J)$, we do not want to do that exclusively. If \eqref{5.3} holds,
we say that $\calJ$ is of type $1$.

If instead
\begin{equation} \lb{5.4}
A_1 \dots A_n >0
\end{equation}
for all $n$, we say $\calJ$ is of type $2$.

An $\ell\times\ell$ matrix, $K$, is said to be in $\calL$ if it is lower triangular with
strictly positive diagonal elements, that is,
\begin{equation} \lb{5.5}
K_{ij} = \begin{cases}
0 & \text{if } i<j \\
>0 & \text{if } i=j
\end{cases}
\end{equation}
If each $A_q\in\calL$, we say that $\calJ$ is of type $3$.

The calculations in Section~\ref{s2} and \ref{s3} show:

\begin{proposition}\lb{P5.1}
\begin{SL}
\item[{\rm{(i)}}] If $\Delta$ is the discriminant of a periodic Jacobi matrix, $J_0$, of period $\ell$,
then for any Jacobi matrix, $J$, $\Delta(J)=\calJ$ is a block Jacobi matrix of type 3.

\item[{\rm{(ii)}}] If $\Delta$ is the discriminant of a periodic CMV matrix, $\calC_0$,
of even period $\ell$, then for any CMV matrix, $\calC$, of $\Delta(\calC)=\calJ$ is
a block Jacobi matrix of type 3.
\end{SL}
\end{proposition}

We will see that distinct $\calJ$'s may correspond to the same measure. Indeed, in the scalar
case, the $b_n$'s and $\abs{a_n}$'s are fixed by the measure, but the $\arg (a_n)$'s are
arbitrary. Thus, we define

\begin{definition} Two block Jacobi matrices, $\calJ$ and $\wti\calJ$, are called equivalent
if and only if there is an $\ell\times\ell$ block diagonal unitary $\calU =\bdone\oplus U_2
\oplus U_3\oplus\cdots$ (we will use $U_1$ for $\bdone$) so that
\begin{equation} \lb{5.6}
\ti\calJ=\calU\calJ\calU^{-1}
\end{equation}
This is equivalent to
\begin{subequations}\lb{5.7}
\begin{align}
\wti B_n &= U_n B_n U_n^{-1} \lb{5.7a} \\
\ti A_n &= U_n A_n U_{n+1}^{-1} \lb{5.7b}
\end{align}
\end{subequations}
\end{definition}

We will be interested in equivalence classes of $\calJ$'s.

\begin{proposition}[\cite{DPSprep}] Each equivalence class of $\calJ$'s has exactly one element
each of type 1, type 2, and type 3.
\end{proposition}

\begin{definition} The Nevai class is the set of $\calJ$'s for which
\begin{equation} \lb{5.8}
B_n\to 0 \qquad A_n^\dagger A_n^{ }\to\bdone
\end{equation}
\end{definition}

The following is immediate from \eqref{5.7}:

\begin{proposition}\lb{P5.3} If some $\calJ$ is in the Nevai class, so are all
equivalent $\calJ$'s.
\end{proposition}

Damanik, Pushnitski, and Simon \cite{DPSprep} prove that

\begin{proposition}[\cite{DPSprep}]\lb{P5.4} If $\calJ$ is in the Nevai class and is type 1,
type 2, or type 3, then
\begin{equation} \lb{5.9}
A_n\to \bdone
\end{equation}
\end{proposition}

We will sometimes need the MOPRL, the matrix orthogonal polynomials.  What we describe here are
the left OPs. There are also right OPs (see \cite{DPSprep}), which we do not need here. An $\ell$-dimensional
matrix-valued measure is a positive scalar measure $d\mol_t(x)$ and a nonnegative $\ell\times \ell$ matrix-valued function $M(x)$.
The matrix-valued measure
\begin{equation} \lb{5.10}
d\mol(x) = M(x)\, d\mol_t(x)
\end{equation}
can always be normalized by
\begin{equation} \lb{5.11}
\tr(M(x))=\ell
\end{equation}

We will always assume $d\mol$ is normalized, that is,
\begin{equation} \lb{5.12}
\int d\mol(x)=\bdone
\end{equation}
The proper notion of nontriviality is a little subtle; it is discussed in detail in \cite{DPSprep}.
For our purpose here, it is sufficient that $\llangle \cdot,\cdot\rrangle_\Lt$ defined below
is nondegenerate on polynomials.

If $f,g$ are two $\ell$-dimensional matrix-valued functions, we define
\[
\llangle f,g\rrangle_\Lt=\int g(x) M(x) f(x)^\dagger\, d\mol_t(x)
\]
Note that this `inner product' returns matrix values.  Recall also that ${}^\dagger$ denotes the Hermitian conjugate of
a matrix.  The subscript `$L$' is for `left' and reflects the fact that if $C$ is an $\ell\times\ell$ matrix, then
\begin{align}
\llangle f,Cg\rrangle_\Lt &= C\llangle f,g\rrangle_\Lt \lb{5.13} \\
\llangle Cf,g\rrangle_\Lt &= \llangle f,g\rrangle_L C^\dagger \lb{5.14}
\end{align}
We will normally just write $\llangle \cdot,\cdot\rrangle$ from now on.

Left orthonormal polynomials are of the form
\[
p_n(x) = \kappa_n x^n + \text{lower order}
\]
with matrix coefficients, defined by
\begin{equation} \lb{5.15}
\llangle p_n,p_m\rrangle =\delta_{nm} \bdone
\end{equation}
So long as $d\mol$ is nontrivial, the $p_n$ exist. They are not unique since if $\{U_n\}_{n=1}^\infty$
are unitary $\ell\times\ell$ matrices,
\begin{equation} \lb{5.16}
\ti p_n(x) = U_{n+1} p_n(x)
\end{equation}
are also MOPRL. We demand $\kappa_0=\bdone$, that is, $p_0(x)=\bdone$, and so $U_1 =\bdone$.

$\{p_j\}_{j=0}^n$ are a left module basis for matrix polynomials of degree $n$, that is,
if $f$ is any polynomial of degree $n$, then there are unique $\ell\times\ell$ matrices $f_0, \ldots, f_n$ so that
\[
f(x)=\sum_{j=0}^n f_j p_j(x)
\]
Indeed,
\begin{equation} \lb{5.16a}
f_j =\llangle p_j,f\rrangle
\end{equation}

For $n=1,2,\ldots$, define
\begin{equation} \lb{5.17}
B_n=\llangle p_{n-1}, xp_{n-1}\rrangle \qquad A_n = \llangle p_n, xp_{n-1}\rrangle
\end{equation}
Then, since $xp_j =\sum_{\ell=0}^{j+1} C_\ell p_\ell$ implies $\llangle p_j, xp_n\rrangle
=\llangle xp_j, p_n\rrangle=0$ if $j\leq n-2$, we have
\begin{equation} \lb{5.18}
xp_n(x) =A_{n+1} p_{n+1}(x) + B_{n+1} p_n(x) + A_n^\dagger p_{n-1}(x)
\end{equation}

This implies $A_{n+1} \kappa_{n+1} = \kappa_n$ so
\begin{equation} \lb{5.19}
\kappa_n = (A_1 \dots A_n)^{-1}
\end{equation}
and the type $2$ condition is equivalent to $\kappa_n >0$.

Looking at \eqref{5.18}, we see that \eqref{5.16} holds for $\ti p_n, p_n$ if and only if
$\ti A_n, \wti B_n$ are related to $A_n,B_n$ by \eqref{5.7}. Jacobi matrix equivalence is
just a `change of phase' in the MOPRL.

Given a block Jacobi matrix, we can view it as acting on the Hilbert space $\ell^2
(\{1,2,\dots\},\bbC^\ell)$ with inner product
\begin{equation} \lb{5.20}
\langle f,g\rangle = \sum_{n=1}^\infty \langle f_n, g_n\rangle_{\bbC^\ell}
\end{equation}
If $\{e_j\}_{j=1}^\ell$ is the standard basis of $\bbC^\ell$, then
$\{\delta_{k;j}\}_{k=1}^\infty {}_{j=1}^{\ell}$, defined by
\begin{equation} \lb{5.21}
(\delta_{k;j})_n = \delta_{kn} e_j
\end{equation}
is a basis. $\calJ$ acts on $\ell^2(\{1,2,\dots\},\bbC^\ell)$ via
\begin{equation} \lb{5.22}
(\calJ f)_n = A_{n-1}^\dagger f_{n-1} + B_n f_n + A_n f_{n+1}
\end{equation}
(with $A_0=0$).

The spectral measure for $\calJ$ is the $\ell\times\ell$ matrix-valued measure with
\begin{equation} \lb{5.23}
\langle \delta_{0;j}, f(\calJ) \delta_{0;k}\rangle =
\int f(x)\, d\mol_{jk}(x)
\end{equation}
for any scalar-valued function $f$. It is easy to see (e.g., \cite{DPSprep}) that this map from
$\calJ$ to $\mol$ inverts the one given by forming the MOPRL and defining $\calJ$
by \eqref{5.17}. Moreover, $\calJ$ and $\wti\calJ$ are equivalent if and only if
$d\ti\mol \equiv d\mol$.

The $m$-function is defined by
\begin{align}
m(E) &= \int \f{1}{x-E}\, d\mol(x) \lb{5.24} \\
&= \langle \delta_{0;\,\bddot}\, , (\calJ-E)^{-1} \delta_{0;\, \bddot}\,\rangle \lb{5.25}
\end{align}
It is an $\ell\times\ell$ matrix-valued Herglotz function:
\begin{equation} \lb{5.26}
\Ima E>0\Rightarrow \Ima m(E)>0
\end{equation}
that is, $\tfrac{1}{2i}(m-m^\dagger)$ is positive definite in the upper half-plane.
For information on matrix Herglotz functions, see
\cite{AG94,AG95,BT92a,BT92b,Fu76,Ger77,GT2000,GKS98,HS81,Jo87,Kr89,KO78,Sak92,Wei87}.
Obviously, by \eqref{5.24}, $m$ is constant over equivalence classes.

As in the scalar case (see \cite[Section~1.2]{OPUC1}), one has that for a.e.\ $x\in\bbR$,
$\lim_{\veps\downarrow 0} m(x+i\veps)\equiv m(x+i0)$ exists and
\begin{equation} \lb{5.27a}
d\mol_\ac (x) = \pi^{-1} \Ima m (x+i0)\, dx
\end{equation}
Here
\begin{equation} \lb{5.27b}
d\mol_\ac (x) = M(x) \, d\mol_{t;\ac}(x)
\end{equation}
where $d\mol_{t;\ac}$ is the a.c.\ part of $d\mol_t$. Alternatively, $d\mol_\ac$ is the
unique matrix-valued measure which is a.c.\ (i.e., $\mol_\ac (I) =0$ for any set with
$\abs{I}=0$) and where there is a set $K$ with $\abs{K}=0$ so $(\mol-\mol_\ac) (\bbR\setminus K)=0$.

Given a block Jacobi matrix, $\calJ$, by $\calJ^{(n)}$ we mean the matrix with the top $n$
(block matrix) rows and leftmost $n$ columns removed, that is,
\begin{equation} \lb{5.27x}
B_k^{(n)} = B_{k+n}^{} \qquad A_k^{(n)} = A_{k+n}^{ }
\end{equation}
We write $m^{(n)}(z)$ for the $m$-function associated to $\calJ^{(n)}$. Equivalent $\calJ$'s do not have the same
$m^{(n)}$ for $n\geq 1$ (although they are unitarily related). We see $m^{(0)}\equiv m$.

We will need the following result of Aptekarev--Nikishin \cite{AN84} (see also \cite{DPSprep}),
a matrix analog of the well-known Jacobi--Stieltjes recursion for OPRL:

\begin{theorem}[\cite{AN84,DPSprep}]\lb{T5.5} We have that
\begin{equation} \lb{5.27}
m^{(n)} (E)^{-1} = E  -B_{n+1} -A_{n+1} m^{(n+1)} (E) A_{n+1}^\dagger
\end{equation}
for $n=0,1,2,\dots$.
\end{theorem}

Next, we need to note the following analog of a well-known scalar result (see, e.g.,
\cite{Sim_Sturm}) proven in \cite{DPSprep}:

\begin{theorem}\lb{T5.5A} Let $\calJ$ be a block Jacobi matrix with $\sigma_\ess (\calJ)
\subset [a,b]$.  Then, for any $\veps$, there is a $K$ so that for $k\geq K$,
\begin{equation} \lb{5.32}
\sigma(\calJ^{(k)})\subset [a-\veps, b+\veps]
\end{equation}
\end{theorem}

Finally, we need to look at poles and zeros of $\det (m(z))$. In the scalar $(\ell=1)$
case, poles occur precisely at eigenvalues of $\calJ$ and zeros at eigenvalues of
$\calJ^{(1)}$, the once stripped Jacobi matrix. In that scalar case, these eigenvalues
are distinct.

In the matrix case, $\calJ$ and $\calJ^{(1)}$ can have eigenvalues in common (as can
be easily arranged by taking a direct sum of suitable scalar $J$'s) so there can
be cancellations. We say a scalar meromorphic function, $f(z)$, has a zero/pole of
order $k\in\bbZ$ at $z_0\in\bbC$ if $(z-z_0)^{-k} f(z)$ has a finite nonzero limit
as $z\to z_0$. We will need the following result from \cite{DPSprep}:

\begin{theorem}[\cite{DPSprep}]\lb{T5.6} Let $x_0\in\bbR$. Let $q_0$ be the multiplicity
of $x_0$ as an eigenvalue of $\calJ$, and $q_1$ its multiplicity as an eigenvalue of $\calJ_1$.
Then
\begin{SL}
\item[\rm{{(a)}}] $q_0 + q_1 \leq \ell$
\item[{\rm{(b)}}] $\det(m(z))$ has a zero/pole of order $q_1-q_0$.
\end{SL}
\end{theorem}

We will also need the following result from Aptekarev--Nikishin \cite{AN84}:

\begin{theorem}\lb{T5.8} Let $\calJ$ be a block Jacobi matrix with $\sigma_\ess (\calJ)=
[-2,2]$, $\sigma(\calJ)\setminus\sigma_\ess (\calJ)$ a finite set and with $g(x)=d\mol_\ac
(x)/dx$ we have
\begin{equation} \lb{5.33}
\int (4-x^2)^{-1/2} \log (\det(g(x))\, dx > -\infty
\end{equation}
Suppose $\calJ$ is type 2. Then
\[
\lim_{n\to\infty}\, A_1 \dots A_n
\]
exists and is a strictly positive matrix.
\end{theorem}

\section{A Denisov--Rakhmanov Theorem for MOPRL} \lb{s6}

As preparation for proving Theorem~\ref{T1.2} in Section~\ref{s7}, in
this section we will prove

\begin{theorem}\lb{T6.1} Let $d\mol$ be a nondegenerate $\ell\times\ell$ matrix-valued
measure on $\bbR$ with associated block Jacobi matrix $\calJ$ of type 3 so that
\begin{alignat}{2}
&\text{\rm{(i)}} \qquad && \sigma_\ess (\calJ) = [-2,2] \lb{6.1} \\
&\text{\rm{(ii)}} \qquad && d\mol = f(x)\, dx + d\mol_\s \lb{6.2}
\end{alignat}
with $d\mol_\s$ singular and
\begin{equation} \lb{6.3}
\det(f(x)) >0
\end{equation}
a.e.\ on $[-2,2]$. Then
\begin{equation} \lb{6.4}
B_n\to \bdzero \qquad A_n\to \bdone
\end{equation}
\end{theorem}

\begin{remark} \eqref{6.3} says the a.c.\ spectrum has multiplicity $\ell$.
\end{remark}

If \eqref{6.1} is replaced by the stronger $\sigma(J)=[-2,2]$ and type $3$ by type $2$,
this is a theorem of Yakhlef--Marcell\'an \cite{YM2001}. We will prove Theorem~\ref{T6.1}
by modifying their proof.

The shift from type $2$ to  $3$ is easy on account of Proposition~\ref{P5.4}. By applying
the argument of \cite{YM2001}, we get $\ti A_n\to\bdone$ for the equivalent $\ti\calJ$ of type $2$, conclude
the whole equivalence class is in the Nevai class, and see $A_n\to\bdone$. So we will only
worry about the changes needed to go from $\sigma(J) =[-2,2]$ to $\sigma_\ess (J)=
[-2,2]$, where we follow Denisov's approach for the scalar case \cite{DenPAMS}.

\cite{YM2001} relies on a matrix version of Rakhmanov's theorem proven by van Assche
\cite{newVA}. We need to extend it slightly to allow a.c.\ spectrum on a large subset of
$\partial\bbD$ rather than all of $\partial\bbD$:

\begin{theorem}\lb{T6.2} Let $d\mu$ be an $\ell\times\ell$ matrix-valued measure on
$\partial\bbD$ and let $\{\alpha_n\}_{n=0}^\infty$ denote its matrix Verblunsky coefficients .
Suppose
\begin{equation} \lb{6.5}
d\mu = w(\theta)\, \f{d\theta}{2w} + d\mu_\s
\end{equation}
where $d\mu_\s$ is singular, and let
\begin{equation} \lb{6.6}
\Omega =\{\theta : \det (w(\theta)) >0\}
\end{equation}
Then
\begin{equation} \lb{6.7x}
\limsup_{n\to\infty}\, \|\alpha_n\| \leq 2 \sqrt{2\ell} \, \biggl( 1 -
\biggl( \f{\abs{\Omega}}{2\pi}\biggr)^3\biggr)^{1/2}
\end{equation}
\end{theorem}

\begin{remarks} 1. For notation on MOPUC, see \cite{DPSprep}.

2. Where we use $\{\alpha_n\}_{n=0}^\infty$, van Assche \cite{newVA} uses $\{H_n\}_{n=1}^\infty$
related to $\alpha_n$ by
\begin{equation} \lb{6.7}
H_n = -\alpha_{n-1}^\dagger
\end{equation}

3. We follow notation from \cite{newVA} and the variant of the scalar proof as in
\cite[Section~9.1]{OPUC2} where $a_n, b_n, c_n, d_n$ below all appear.

We define
\begin{align*}
a_n & = \| \alpha_n \| \\
b_{n,q} & = \frac{1}{2\pi} \int_0^{2\pi} \| [ \varphi_n^\Lt (e^{i\theta})
\varphi_{n+q}^\Lt (e^{i\theta})^{-1} ] \, [ \varphi_n^\Lt (e^{i\theta})
\varphi_{n+q}^\Lt(e^{i\theta})^{-1} ]^\dagger - I \| \, d\theta \\
c_{n,q} & = \frac{1}{2 \pi \ell} \int_0^{2\pi} \tr ( \varphi_n^\Lt(e^{i\theta})
\varphi_{n+q}^\Lt(e^{i\theta})^{-1} [ \varphi_{n+q}^\Lt(e^{i\theta})^\dagger ]^{-1}
\varphi_n^\Lt(e^{i\theta})^\dagger )^{1/2} \, d\theta \\
d_n & = \frac{1}{2 \pi \ell} \int_0^{2\pi} \tr ( [ \varphi_n^\Lt
(e^{i\theta}) w(\theta) \varphi_n^\Lt(e^{i\theta})^\dagger ]^{1/2} ) \, d\theta
\end{align*}
\end{remarks}

\begin{proposition}\lb{P6.3} For every $n \geq 0$, we have that
\begin{alignat}{2}
 a_n & \leq b_{n,q} && \qquad \text{ for every } q \geq 1  \label{a} \\
 b_{n,q}^2 & \leq 8 \ell^2 (1 - c_{n,q}) && \qquad \text{ for every } q \geq 1  \label{b} \\
 d_n^2 & \leq \inf_{q \geq 1} c_{n,q} \label{c}
\end{alignat}
Moreover, we have that
\begin{equation}\label{d}
\biggl( \frac{\abs{\Omega}}{2\pi} \biggr)^{3/2} \leq \liminf_{n \to \infty} \, d_n
\end{equation}
Consequently,
\begin{equation} \label{e}
\limsup_{n \to \infty} \, a_n \leq 2 \sqrt{2}\, \ell
\biggl( 1 - \biggl( \frac{\abs{\Omega}}{2\pi} \biggr)^3 \biggr)^{1/2}
\end{equation}
\end{proposition}

\begin{proof} The second to last displayed formula on \cite[p.~7]{newVA} is \eqref{a}.
The estimates on the bottom half of \cite[p.~12]{newVA} show that
\begin{align*}
b_{n,q}^2 & = \frac{1}{4 \pi^2} \biggl( \int_0^{2\pi} \| [
\varphi_n^\Lt(e^{i\theta}) \varphi_{n+q}^\Lt(e^{i\theta})^{-1} ] \,
[ \varphi_n^\Lt(e^{i\theta})
\varphi_{n+q}^\Lt(e^{i\theta})^{-1} ]^\dagger - I \| \, d\theta \biggr)^2 \\
& \leq \frac{2 \ell}{\pi} \int_0^{2\pi} \| ( [ \varphi_n^\Lt(e^{i\theta})
\varphi_{n+q}^\Lt(e^{i\theta})^{-1} ]\, [ \varphi_n^\Lt(e^{i\theta})
\varphi_{n+q}^\Lt(e^{i\theta})^{-1} ]^\dagger )^{1/2} - I\|^2 \, d\theta \\
& \leq \frac{2\ell}{\pi}\, ( 4 \pi \ell - 4 \pi \ell \, c_{n,q} ) \\
& = 8 \ell^2 (1 - c_{n,q})
\end{align*}
which is \eqref{b}. The third displayed formula on \cite[p.~14]{newVA} is \eqref{c}.

Now, mimicking the estimates on the bottom half of \cite[p.~14]{newVA},
\[
\begin{split}
\int_\Omega \tr &( [ f(\theta) w(\theta) f(\theta)^\dagger ]^{1/4}
) \, d\theta \\
& \leq \biggl( 2 \pi \ell \int_\Omega \tr ( f(\theta)
\varphi_n^\Lt(e^{i\theta})^{-1} ( \varphi_n^\Lt(e^{i\theta})^\dagger )^{-1} f(\theta)^\dagger )\,
 \frac{d\theta}{2\pi} \biggr)^{1/4} ( 2 \pi \ell \, d_n )^{1/2}
 \end{split}
\]
Taking $n \to \infty$, we see that
\[
\begin{split}
\int_\Omega \tr & ([ f(\theta) w(\theta) f(\theta)^\dagger ]^{1/4})
\, d\theta \\
& \leq \biggl( 2 \pi \ell \, \tr \int_\Omega f(\theta) d\mu(\theta)
f(\theta)^\dagger )^{1/4} \bigl( 2 \pi \ell \, \liminf_{n \to \infty}\, d_n \bigr)^{1/2}
\end{split}
\]
Removing the singular part as in \cite{newVA}, we obtain
\[
\begin{split}
\int_\Omega \tr &( [ f(\theta) w(\theta) f(\theta)^\dagger ]^{1/4}) \,
d\theta \\
&  \leq \biggl( 2 \pi \ell \, \tr \int_\Omega f(\theta) w(\theta)
f(\theta)^\dagger d\theta \biggr)^{1/4} \bigl( 2 \pi \ell \, \liminf_{n \to \infty}\,
d_n \bigr)^{1/2}
\end{split}
\]
Proceeding as in \cite[pp.~15--16]{newVA}, it then follows that
\[
\abs{\Omega} \ell \leq ( 2 \pi \ell \,\abs{\Omega} \ell )^{1/4} \bigl( 2 \pi \ell \, \liminf_{n
\to \infty}\,  d_n\bigr)^{1/2}
\]
which implies \eqref{d}.

Putting these estimates together,
\[
a_n \leq b_{n,1} \leq 2 \sqrt{2} \, \ell (1 - c_{n,1})^{1/2} \leq 2 \sqrt{2}\, \ell (1 - d_n^2)^{1/2}
\]
and hence
\[
\limsup_{n \to \infty}\, a_n \leq 2 \sqrt{2}\, \ell
\biggl( 1 - \biggl( \frac{\abs{\Omega}}{2\pi}\biggr)^3 \biggr)^{1/2}
\]
which is \eqref{e}.
\end{proof}

In particular, for $2\pi - \abs{\Omega}$ small,
\[
\limsup_{n \to \infty}\, a_n = O \biggl( \biggl( 1 - \biggl( \f{\abs{\Omega}}{2\pi} \biggr)^3
\biggr)^{1/2} \biggr) = O (( 2\pi - \abs{\Omega} )^{1/2})
\]
as in the scalar case.

We have thus proven Theorem~\ref{T6.2}. To get Theorem~\ref{T6.1}, we follow \cite{YM2001} using
the analog of Denisov's arguments for the case $\ell=1$.

\begin{proof}[Proof of Theorem~\ref{T6.1}] By Proposition~\ref{P5.4}, we need only prove
for the type $2$ choice, for any $\veps >0$, we have
\begin{equation} \lb{6.8}
\limsup \, (\|\ti A_n-1\| + \|\ti B_n\|) \leq \veps
\end{equation}
By the Szeg\H{o} mapping and Geronimus connection formulae in \cite{YM2001}, this holds by
Theorem~\ref{T6.2} so long as for any $\veps' >0$, we can find $k$ so $\sigma (\calJ^{(k)})
\subset [-2-\veps', 2+\veps']$, and this is true by Theorem~\ref{T5.5A}.
\end{proof}

\section{A Denisov--Rakhmanov Theorem for Periodic OPRL} \lb{s7}

Our main goal in this section is to prove Theorem~\ref{T1.2}. We will also prove the
`hard' half of Theorem~\ref{T1.1}. The simplicity of the proof shows the magic in the
magic formula!

\begin{proof}[Proof of Theorem~\ref{T1.2}] By a right limit of $J$, we mean a two-sided
Jacobi matrix, $J_r$, (but with some $a$'s allowed to vanish) so that for some subsequence
$n_j\to\infty$ and any $k\in\bbZ$,
\begin{equation}\lb{7.1}
a_{n_j+k}\to (a_r)_k \qquad b_{n_j+k}\to (b_r)_k
\end{equation}
By our standing convention, Jacobi parameters are uniformly bounded, so by compactness,
if $d_m ((a,b),\calT_{J_0})\nrightarrow 0$, there exists a right limit $J_r\notin\calT_{J_0}$.
Thus, it suffices to show that any right limit $J_r$ has $J_r\in\calT_{J_0}$.

By the hypotheses of Theorem~\ref{T1.2}, the spectral mapping theorem, and the fact that
$\Delta$ maps $\sigma_\ess (J_0)$ to $[-2,2]$ with a $p$-fold cover on $(-2,2)$, we see
that
\[
\Delta(J)_\ess = [-2,2]
\]
and $\Delta (J)$ has a.c.\ spectrum of multiplicity $p$. So thinking of $\calJ\equiv\Delta(J)$
as a block Jacobi matrix, $\calJ$ is of type $3$ and the hypotheses of Theorem~\ref{T6.1} apply.
It follows that $A_n\to \bdone$, $B_n\to \bdzero$. This means that
\[
\Delta (J_r) =S^p + S^{-p}
\]
so by the magic formula (Theorem~\ref{T2.1}), $J_r\in\calT_{J_0}$.
\end{proof}

Rakhmanov's theorem is often related to issues of $\wlim \abs{p_n}^2\, d\mu$ and to the
density of zeros. We note that there are also results of that genre here:

\begin{theorem}\lb{T7.1A} If $J_0$ is a periodic Jacobi matrix of period $p$ and $J$
is a Jacobi matrix with bounded Jacobi parameters whose right limits all lie in
$\calT_{J_0}$ {\rm{(}}in particular, if the hypotheses of Theorem~\ref{T1.2} hold{\rm{)}},
then {\rm{(}}with $d\mu$ the measure for $J${\rm{)}}
\begin{SL}
\item[{\rm{(a)}}]
\begin{equation}\lb{7.2}
\wlim_{n\to\infty}\, \f{1}{p} \sum_{j=1}^p\, \abs{p_{j+n}(x)}^2\, d\mu(x)=d\nu
\end{equation}
the density of zeros for $J_0$.
\item[{\rm{(b)}}] The density of zeros of $p_n(x)$ converges to $d\nu$.
\end{SL}
\end{theorem}

\begin{proof} If $J_1\in\calT_{J_0}$, then the spectral measure $d\mu_k^{(J_1)}$ associated
to $\delta_k$ has period $p$ in $k$ since $J_1$ is periodic. Thus $\lim_{N\to\infty}\f{1}{2N+1}
\sum_{\abs{j}\leq N} d\mu_j^{(J_1)}=\f{1}{p}\sum_{j=1}^p d\mu_j^{(J_1)}$, but the
limit is $d\nu$ by the discussion in Subsection~\ref{ssTF}. Since
\begin{equation} \lb{7.3}
\int x^\ell \abs{p_j^{(J)}(x)}^2\, d\mu(x) = \langle \delta_j, J^\ell \delta_j\rangle
\end{equation}
and $J_{j,j}, J_{j,j\pm 1}$ is very close to some $(J_1)_{j,j}$, $(J_1)_{j,j+1}$
for $\abs{j-j_0}\leq M$ for fixed $M$ and $j_0\to\infty$, we see that moments of LHS
of \eqref{7.2} are close to moments of $d\nu$. This proves (a).

If $J^{(n)}$ denotes the top left $n\times n$ submatrix of $J$, then
\[
\int x^\ell\, d\nu_n =\f{1}{n}\, \tr ((J^{(n)})^\ell)
\]
so
\[
\lim_{n\to\infty}\, \int x^\ell\, d\nu_n = \lim_{n\to\infty}\, \int x^\ell
\biggl[ \, \f{1}{n}\sum_{j=0}^{n-1} p_j(x)^2\, d\mu(x)\biggr]
\]
and thus (a) implies (b).
\end{proof}

We also have

\begin{theorem}\lb{T7.1} If $d_m ((a,b),\calT_{J_0})\to 0$, then
\[
\sigma_\ess (J) \subset \sigma_\ess (J_0)
\]
\end{theorem}

\begin{proof} By the magic formula, compactness, and the fact that every right limit of $J$ is
in $\calT_{J_0}$, we see that every right limit of $\Delta(J)$ is $S^p + S^{-p}$, that is,
$A_n\to\bdone$, $B_n\to \bdzero$. Thus, by Weyl's theorem, $\sigma_\ess (\Delta(J))=[-2,2]$. Since
\[
\sigma_\ess (\Delta(J)) = \Delta (\sigma_\ess (J))
\]
we see $\sigma_\ess (J)\subset \Delta^{-1} ([-2,2])=\sigma_\ess (J_0)$.
\end{proof}

\begin{remarks} 1. Since $\Delta$ is $p$ to $1$, we cannot conclude that $\sigma_\ess (J_0)
=\sigma_\ess (J)$ from $\Delta (\sigma_\ess (J_0))=\Delta (\sigma_\ess (J))$.

2. That $\sigma_\ess (J_0)\subset\sigma_\ess (J)$ is a simple trial function argument given that
$J$ must have some right limits; see \cite{LastS,LS_jdam}.
\end{remarks}

\section{Denisov--Rakhmanov Sets} \lb{new-s9}

In this section, we want to show how one can take suitable limits of Theorem~\ref{T1.2}
to get a `cheap' proof of similar theorems in other nonperiodic cases. We will also
present an insight into the proper general form of Denisov--Rakhmanov-type theorems.

The right limits we have discussed so far involve the weak product topology on the Jacobi
parameters, so we will emphasize this fact by using the phrase `weak right limits' in this section.
We are also interested in limits in the $\ell^\infty$-topology for two-sided sequences, that is,
$\{c_n^{(k)}\}_{n=-\infty}^\infty \to \{c_n^{(\infty)}\}_{n=-\infty}^\infty$ in this
topology if and only if, as $k\to\infty$,
\[
\sup_n\, \abs{c_n^{(k)}-c_n^{(\infty)}}\to 0
\]

In terms of weak limits, we note the following:

\begin{proposition}\lb{Pn9.1} Let $\calE$ be a closed set. Let $J$ be a half-line Jacobi matrix with
\begin{equation} \lb{n9.1}
\Sigma_\ac  (J)=\sigma_\ess(J)=\calE
\end{equation}
Let $J_r$ be a weak right limit of $J$. Then
\begin{equation} \lb{n9.2}
\Sigma_\ac (J_r) =\sigma(J_r) =\calE
\end{equation}
\end{proposition}

\begin{remark} Note that \eqref{n9.2} has $\sigma(J_r)$, not merely $\sigma_\ess (J_r)$.
\end{remark}

\begin{proof} By results in \cite{LastS},
\[
\sigma (J_r)\subset \calE\subset\Sigma_\ac (J_r)
\]
Since $\Sigma_\ac (J_r) \subset \sigma (J_r)$ trivially, \eqref{n9.2} holds.
\end{proof}

Recall that a sequence $\{c_n\}_{n=-\infty}^\infty$ is called uniformly almost periodic
(in the general theory of almost periodic functions, this defines `almost periodic'---we
add `uniformly' because the term is sometimes used in a weaker sense in the spectral theory
literature) if and only if $\{c^{(\ell)}\}_{\ell=-\infty}^\infty$ given by $(c^{(\ell)})_n=
c_{n+\ell}$ has compact closure in the $\ell^\infty$-topology.

\begin{definition} A set $\calE$ is called {\it essentially perfect\/} if and only if $\calE$ is closed,
and for all $E\in \calE$ and $\delta >0$, $\abs{(E-\delta, E+\delta) \cap \calE}>0$.
\end{definition}

\begin{remark} Essentially perfect sets are precisely the sets, $\calE$, for which there is a purely
a.c.\ measure $d\mol$ with $\supp(d\mol)=\calE$.
\end{remark}

\begin{definition} A set $\calE$ is said to be a {\it Denisov--Rakhmanov set\/} if and only
if
\begin{SL}
\item[(i)] $\calE$ is essentially perfect and bounded.
\item[(ii)] There is a set $\calT_{\calE}$ compact in the uniform topology so that for any
bounded Jacobi matrix, $J$, for which \eqref{n9.1} holds, the set of right limits of
$J$ lies in $\calT_{\calE}$.
\end{SL}
\end{definition}

The definition says nothing explicit about $\calT_{\calE}$ being a torus, but by Proposition~\ref{Pn9.1},
if $J_r\in\calT_{\calE}$, then \eqref{n9.2} holds, and since $\calT_{\calE}$ is closed under translations,
each $J_r$ in $\calT_{\calE}$ is almost periodic. By Kotani theory (see \cite{Kot,Kot97,S168,OPUC2}),
$\langle\delta_n, (J_r-E-i\veps)^{-1}\delta_n\rangle$ has real boundary values for a.e.\ $E$. In many
cases and, in particular, if $\calE$ is a finite union of closed intervals, Sodin--Yuditskii \cite{SY}
(see also \cite{AntK,BGHT}) proved there is a natural torus so that any almost periodic $J_r$ with real boundary
values lies in this torus. Thus for such cases, that $\calE$ is a Denisov--Rakhmanov set can be
connected to approach to an isospectral torus. In particular, our Theorem~\ref{T1.2} implies
the statement that $\sigma_\ess (J_0)$ is a Denisov--Rakhmanov set.

Given an essentially perfect set, $\calE$, we define $\calD(\calE)$ to be the set of Jacobi matrices
obeying \eqref{n9.1}.

The following two simple results will be the basis of our approximation theorems:

\begin{proposition}\lb{Pn9.2} Let $\calE$ be an essentially perfect set. Suppose there are
uniformly compact sets $\{\calT^{(n)}\}_{n=1}^\infty$ and $\calT^{(\infty)}$ of
two-sided Jacobi matrices so that
\begin{SL}
\item[{\rm{(1)}}] If $J_n\in\calT^{(n)}$ and $J_n\to J_\infty$ weakly, then $J_\infty
\in\calT^{(\infty)}$.
\item[{\rm{(2)}}] For any weak right limit point $J_r$ of some $J\in\calD(\calE)$, there is
$\ti J\in\calT^{(n)}$ so
\begin{equation} \lb{n9.3}
\sup_{\abs{j}\leq n}\, \abs{a_j^{(r)} -\ti a_j} + \abs{b_j^{(r)}-\ti b_j} \leq \f{1}{n}
\end{equation}
\end{SL}
Then $\calE$ is a Denisov--Rakhmanov set.
\end{proposition}

\begin{proof} Let $J_n$ be the $\ti J$ guaranteed by \eqref{n9.3}. Then clearly, $J_n$
converges weakly to $J_r$ so, by (1), $J_r\in\calT^{(\infty)}$. Since $\calT^{(\infty)}$
is uniformly compact, $\calE$ is a Denisov--Rakhmanov set.
\end{proof}

\begin{proposition}\lb{Pn9.3} Let $J_0$ be a fixed periodic Jacobi matrix with essential
spectrum $\calE_0$. Then for any $n$, there is a $\delta >0$ so that for any set $\calE$ with
\begin{alignat}{2}
&\text{\rm{(a)}} \qquad && \calE \subset \{E : \dist (E,\calE_0)<\delta\} \lb{n9.4} \\
&\text{\rm{(b)}} \qquad && \abs{\calE} > (1-\delta) \abs{\calE_0} \lb{n9.5}
\end{alignat}
and any $J\in\calD(\calE)$, we have that any right limit, $J_r$, obeys \eqref{n9.3}
for some $\ti J\in\calT_{J_0}$.

Moreover, if $p$ is fixed and $C$ is a compact subset of $[(0,\infty)\times
\bbR]^p$, then $\delta$ can be picked to work for all $J_0 =\{(a_n,b_n)\}_{n=1}^p\in C$.
\end{proposition}

\begin{proof} The uniformity claimed in the last statement comes from noting that choices
can be made uniformly in the proof below.

Let $p$ be the period of $J_0$. We first claim that given $\delta_1$, we can find
$\delta$ so if $\calE$ obeys \eqref{n9.4}--\eqref{n9.5}, then
\begin{gather}
\dist (\Delta(\calE), [-2,2]) < \delta_1 \lb{n9.6} \\
\abs{\{x\in (-2,2) : \text{all $p$ solutions of $\Delta(E)=x$ lie in }\calE\}} >
4-\delta_1 \lb{n9.7x}
\end{gather}
This is immediate from the continuity of $\Delta$ and its derivatives.

Next, we note that given $\veps_1$, we can find $\delta_1$ so that if $\calJ$ is
a $p\times p$ block Jacobi matrix so that
\begin{gather*}
\sigma_\ess (\calJ)\subset [-2-\delta_1, 2+\delta_1] \\
\abs{\{E\in [-2,2] : \calJ\text{ has a.c. spectrum at $E$ of multiplicity $p$}\}}
>4-\delta_1
\end{gather*}
then
\[
\limsup_{k,m\to\infty}\, \abs{\calJ_{km} - (S^p + S^{-p})_{km}} <\veps_1
\]
The proof of this is identical to the proof of the matrix Denisov--Rakhmanov theorem.

Combining these steps, we are reduced to showing for any $n$ and $\veps$, there
is $\veps_1$ so for all two-sided $J_r$ with $\dist(\sigma(J_r), \calE) <\veps$, we have that
\begin{equation} \lb{n9.7}
\sup_{k,m}\, \abs{\Delta(J_r)-(S^p+S^{-p})_{km}} <\veps_1
\end{equation}
implies there is a $\ti J\in\calT_{J_0}$ so that \eqref{n9.3} holds. To do this, we
first follow the proof of Theorem~\ref{T2.1} to note that for $n$, $\veps_2$, and
$\veps_3$ fixed, we can find $\veps_1$ so \eqref{n9.7} implies there is a $p$-periodic $J^\sharp$ such that
\begin{equation} \lb{n9.8}
\|\Delta (J^\sharp)-\Delta(J)\| <\veps_2
\end{equation}
and
\[
\sup_{\abs{j}\leq n}\, \abs{a_j^{(r)}-a_j^\sharp} +\abs{b_j^{(r)}-\ti b_j}\leq\veps_3
\]

Finally, a compactness argument shows that for any $n$, we can find $\veps_4$ so
for any periodic $J^\sharp$ with
\[
\|\Delta(J^\sharp)-(S^p +S^{-p})\| <\veps_4
\]
there is a $\ti J\in\calT_{J_0}$ so that
\[
\|J^\sharp-\ti J\| \leq \f{1}{2n}
\]
Putting these together implies \eqref{n9.3}.
\end{proof}

\begin{theorem}\lb{Tn9.4} Let $\ell_1,\ell_2, \dots$ be an arbitrary sequence in
$(2,3,4,\dots)$. For any $\ell_1$-periodic Jacobi matrix $J^{(0)}$, there exist
$k_2, k_3, \dots$ so that for any limit periodic $J$ with Jacobi coefficients
\begin{align}
a_n &=a_n^{(0)} +\sum_{m=2}^\infty \Real [A_m e^{2\pi in/\ell_1\ell_2\dots\ell_m}] \lb{n9.9} \\
b_n &=b_n^{(0)} +\sum_{m=2}^\infty \Real [B_m e^{2\pi in/\ell_1\dots\ell_m}] \lb{9.10}
\end{align}
obeying
\begin{equation} \lb{n9.11x}
\abs{A_m} + \abs{B_m} \leq k_m
\end{equation}
we have that $\sigma(J)$ is a Denisov--Rakhmanov set.
\end{theorem}

\begin{remark} The study of limit periodic discrete Schr\"odinger operators with small
tails was initiated by Avron--Simon \cite{S147} and Chulaevsky \cite{Chu}. They prove purely
a.c.\ spectrum.
\end{remark}

\begin{proof} As in \cite{S147,Chu}, one can pick the $k_m$'s so the
spectrum is purely a.c.\ and so that the union of all isospectral tori for the periodic
approximates lie in a fixed $\ell^\infty$ compact sets. This implies the limit periodic
potentials also have compact isospectral sets, and within this compact set, weak
convergence implies norm convergence so hypothesis (1) of Proposition~\ref{Pn9.2} holds.
By decreasing the $k_m$'s if necessary, Proposition~\ref{Pn9.3}, continuity of the
spectrum in $\ell^\infty$ norm, and absolute continuity of periodic spectrum imply
we can be sure that \eqref{n9.3} holds. Thus Proposition~\ref{Pn9.2} implies
this theorem.
\end{proof}

Our final theorem in this section is the following:

\begin{theorem}\lb{Tn9.5} Fix $\ell$. Let $\calG=\{(\alpha_1, \beta_1, \alpha_2, \dots,
\alpha_{\ell+1}, \beta_{\ell+1})\in\bbR^{2\ell+2} :  \alpha_1 <\beta_1 < \alpha_2 < \beta_2 <
\cdots < \beta_{\ell+1}\}$. For $(\vec{\alpha}, \vec{\beta})\in\calG$, define
\begin{equation} \lb{n9.11}
\calE(\vec{\alpha},\vec{\beta})=\bigcup_{j=1}^{\ell+1}\, [\alpha_j,\beta_j]
\end{equation}
Then $\{(\vec{\alpha},\vec{\beta}) : \calE(\vec{\alpha},\vec{\beta})\text{ is a Denisov--Rahkmanov
set}\}$ contains a dense $G_\delta$.
\end{theorem}

\begin{remarks} 1. As we have seen, the $\calE(\vec{\alpha},\vec{\beta})$ which arise from
periodic problems are precisely those where the harmonic measure of each $e_j=[\alpha_j, \beta_j]$
is rational. In particular, if we fix $\vec{\alpha}$ and $\beta_{\ell+1}$, the set of $\beta$'s
that are periodic is countable, and so certainly not a $G_\delta$. We show that the family that
leads to Denisov--Rahkmanov sets is uncountable.

\smallskip
2. It is a reasonable conjecture that every $\calE$ is a Denisov--Rakhmanov set, so this result
is weak. We include it because it is such a `cheap' way to go beyond the periodic case
using only that case.
\end{remarks}

\begin{proof} For each $(\vec{\alpha},\vec{\beta})$, it is known \cite{SY} that there is
an isospectral torus $\calT$ of almost periodic $J$'s where (whole line) spectrum
is precisely $\calE(\vec{\alpha},\vec{\beta})$. It follows from the construction in \cite{SY}
that if $(\vec{\alpha}^{(n)}, \vec{\beta}^{(n)})\in\calG$ converge to $(\vec{\alpha}^{(\infty)},
\vec{\beta}^{(\infty)})\in\calG$, then condition (1) of Proposition~\ref{Pn9.2} holds.

Let $\calG_p$ be the subset of $\calG$ coming from periodic problems---this is dense in $\calG$.
For $(\vec{\alpha}^{(0)},\vec{\beta}^{(0)})\in\calG_p$, pick $J(\vec{\alpha}^{(0)},
\vec{\beta}^{(0)})$ periodic with $\calE(\vec{\alpha}^{(0)},\vec{\beta}^{(0)})$ as spectrum and
pick $\delta_n (\vec{\alpha}^{(0)},\vec{\beta}^{(0)})$ via Proposition~\ref{Pn9.3} requiring
$\delta_n <\f12\min(\abs{\beta_j^{(0)}-\alpha_j^{(0)}}, \abs{\alpha_{j+1}^{(0)} -\beta_j^{(0)}})$.
Let $U^{(n)}(\vec{\alpha}_0,\vec{\beta}_0)=\{(\vec{\alpha},\vec{\beta}) : E(\vec{\alpha},
\vec{\beta})$ obeys \eqref{n9.4}/\eqref{n9.5} for $\calE = \calE(\vec{\alpha},\vec{\beta})$, $\calE_0=
\calE(\vec{\alpha}^{(0)},\vec{\beta}^{(0)})$ and $\delta=\delta_n\}$, and let
\[
U^{(n)}= \bigcup_{\calG_p} U^{(n)} (\alpha^{(0)},\beta^{(0)})
\]
This is dense and open. Then $\cap_n U^{(n)}$ is a dense $G_\delta$ whose points, by
construction and Proposition~\ref{Pn9.2}, correspond to Denisov--Rakhmanov sets.
\end{proof}

\section{Sum Rules for MOPRL} \lb{s8}

In this section, our main goal is to prove the following two theorems about block Jacobi
matrices:

\begin{theorem}[$P_2$ Sum Rule for MOPRL]\lb{T8.1} Let $\calJ$ be a block Jacobi
matrix with $\ell\times\ell$ Jacobi parameters $\{A_n\}_{n=1}^\infty, \{B_n\}_{n=1}^\infty$
and matrix measure $d\mol$. Let $m(E)$ be given by \eqref{5.24} and suppose $\sigma_\ess
(\calJ)=[-2,2]$. Define for $z\in\bbD\setminus \{z=E+E^{-1} : E\in\sigma(\calJ)
\setminus [-2,2]\}$
\begin{equation} \lb{8.1}
M(z) =-m(z+z^{-1})
\end{equation}
Let $F,G$ be the functions
\begin{equation} \lb{8.2}
F(\beta + \beta^{-1}) = \tfrac14\, [\beta^2 -\beta^{-2}-\log \beta^2]
\end{equation}
for $\beta\in\bbR\setminus [-1,1]$, that is, $E=\beta + \beta^{-1}\in\bbR
\setminus [-2,2]$ and
\begin{equation} \lb{8.3}
G(a)=a^2 -1 - \log(a^2) \qquad a\in (0,\infty)
\end{equation}
Then $\lim_{r\uparrow 1} M(re^{i\theta})$ exists for a.e.\ $\theta$ and
\begin{equation} \lb{8.4}
\begin{split}
\f{1}{2\pi}\int \log &\biggl( \f{\sin^\ell \theta}{\det(\Ima M(e^{i\theta}))}\biggr)
\sin^2 \theta\, d\theta \\
&+ \sum_{E\in\sigma(\calJ)\setminus [-2,2]}\, F(E)
= \sum_{n=1}^\infty \tr (\tfrac14\, B_n^2 + \tfrac12\, G(|A_n|))
\end{split}
\end{equation}
\end{theorem}

\begin{remarks} 1. All terms are positive (since $F$ and $G$ are positive, this is
evident for two terms; positivity of the integral will be seen below), so this
sum rule always makes sense, although some terms may be $+\infty$.

2. Recall that $|A_n|=\sqrt{A_n^\dagger A_n^{ }}$; although since the formula for $G(a)$ only involves $a^2$,
one does not need to take a square-root.

3. Because of the trace and absolute value, $\tr (\f14 B_n^2 + \f12 G(\abs{A_n}))$
is constant over equivalence classes of Jacobi matrix parameters.

4. In the type 1 case, the RHS of \eqref{8.4} is finite if and only if $\calJ-S^p-S^{-p}$ is Hilbert--Schmidt.
This is also true when $\calJ$ is of type 3; see Proposition~\ref{P9.12}.
\end{remarks}

\begin{theorem}[Sharp Case $C_0$ Sum Rule for MOPRL]\lb{T8.2} Consider the
three quantities:
\begin{align}
Z(\calJ) &= \f{1}{4\pi}\int_0^{2\pi} \log
\biggl( \f{\sin^\ell\theta}{\det(\Ima M(e^{i\theta}))}\biggr)\, d\theta \lb{8.5} \\
\calE_0(\calJ) &=\sum_{E\notin\sigma(\calJ)\setminus [-2,2]}\, \log (\abs{\beta}) \lb{8.6}
\end{align}
where $\beta$ is related to $E$ by
\begin{equation} \lb{8.4x}
\beta\in\bbR\setminus [-1,1] \qquad E=\beta+\beta^{-1}
\end{equation}
and
\begin{equation} \lb{8.5x}
A_0 (\calJ) = \lim_{N\to\infty} -\sum_{n=1}^N \log (\det(\abs{A_n}))
\end{equation}
which we suppose exists but it may be $+\infty$ or $-\infty$. Then
\begin{SL}
\item[{\rm{(i)}}] If any two of $Z,\calE_0, A_0$ are finite, then so is the
third.

\item[{\rm{(ii)}}] If all are finite, then
\begin{equation} \lb{8.6x}
Z(J)=A_0(J) + \calE_0(J)
\end{equation}

\item[{\rm{(iii)}}] If all are finite, then
\begin{equation} \lb{8.7}
\lim_{N\to\infty} \, \sum_{n=1}^N \, \tr (B_n)
\end{equation}
exists.
\end{SL}
\end{theorem}

\begin{remark}  We will prove (and actually use it to prove Theorem~\ref{T1.3}) that if
$\calE_0(J)<\infty$, then $Z(J)<\infty$ so long as
\begin{equation} \lb{8.8}
\underline{A}_0 (\calJ) =\liminf_N \biggl( -\sum_{n=1}^N \log (\det(\abs{A_n}))\biggr) <\infty
\end{equation}
\end{remark}

Theorem~\ref{T8.1} is a matrix-valued analog of the OPRL $P_2$ sum rule of
Killip--Simon \cite{KS}, and Theorem~\ref{T8.2} of the OPRL Case $C_0$ sum rule by
Simon--Zlato\v{s} \cite{SZ}. Both were refinements of sum rules of Case
\cite{Case1,Case2} who in turn was motivated by earlier KdV and Toda sum rules.
Case only considered short-range $\abs{a_n-1} + \abs{b_n}$, while \cite{KS,SZ}
considered the necessary techniques to go up to the borderline of validity.
\cite{SZ} had some simplifications of \cite{KS}, and \cite{Sim288} further
simplified, although each of the later two proofs depends heavily on the earlier
ones. Here, following Simon \cite{Sim288}, we will prove a nonlocal step-by-step sum
rule. As there, the key is a suitable representation theorem for meromorphic Herglotz
functions---in this case, extended to matrix-valued functions.

For $a\in (-1,1)$, we define Blaschke factors as usual by
\begin{equation} \lb{8.9}
b(z,a) = \begin{cases} \f{a-z}{1-az} & 0<a<1 \\
\f{z-a}{1-az} & -1 < a\leq 0 \end{cases}
\end{equation}

\begin{proposition}\lb{P8.3} Let $f(z)$ be an $\ell\times\ell$ matrix-valued meromorphic
function on $\bbD$ so that
\begin{alignat}{2}
& \text{\rm{(i)}} \qquad && \pm \Ima f(z) >0 \text{ when } \pm\Ima z >0 \lb{8.10} \\
& \text{\rm{(ii)}} \qquad && \lim_{z\to  0}\, f(z) z^{-1} =\bdone \lb{8.11}
\end{alignat}
where $\Ima f \equiv \frac1{2i}(f-f^\dagger)$.  Then
\begin{SL}
\item[{\rm{(a)}}] For a.e.\ $\theta$, $\lim_{r\uparrow 1} f(re^{i\theta})\equiv
f(e^{i\theta})$ exists.

\item[{\rm{(b)}}] $\log \abs{\det(f(e^{i\theta}))} \in \cap_{1\leq p <\infty}
L^p (\partial\bbD, d\theta/2\pi)$

\item[{\rm{(c)}}] All the zeros and poles of $\det (f(z))$ lie on $(-1,1)$ and are
of finite order. Let $\{z_j\}_{j=1}^\infty$ and $\{p_j\}_{j=1}^\infty$ be those
zeros and poles of $\det (f(z))$ repeated up to multiplicity {\rm{(}}it can also happen
that both sets are finite{\rm{)}}. $z=0$ is not included in $\{z_j\}$. Then
\begin{equation} \lb{8.12}
B_\infty (z) =\lim_{r\uparrow 1} \, \f{\prod_{\abs{z_j}<r} b(z,z_j)}
{\prod_{\abs{p_j}<r} b(z, p_j)}
\end{equation}
exists and obeys:
\begin{SL}
\item[{\rm{(i)}}] $B_\infty$ is analytic and nonvanishing on $\bbC\setminus \{z_j\}
\cup \{p_j\}\cup\{z_j^{-1}\}\cup \{p_j^{-1}\}\cup \{\pm 1\}$
\item[{\rm{(ii)}}] $\abs{B_\infty (e^{i\theta})}=1$ on $\partial\bbD\setminus\{\pm 1\}$
\item[{\rm{(iii)}}]
\begin{equation} \lb{8.12a}
\abs{\arg (B_\infty(z))}\leq 2\pi \ell
\end{equation}
for $\abs{z}<1$ with $\arg$ normalized by $\arg B_\infty (0)=0$.
\end{SL}

\item[{\rm{(d)}}] We have the representation
\begin{equation} \lb{8.12b}
\det(f(z)) =z^\ell B_\infty (z) \exp\biggl( \int \f{e^{i\theta}+z}{{e^{i\theta}-z}}\,
\log \abs{\det (f(e^{i\theta}))}\, \f{d\theta}{\pi}\biggr)
\end{equation}
\end{SL}
\end{proposition}

\begin{remarks} 1. It should be possible to prove that $0<\arg (B_\infty (z))
<\pi \ell$ for $\Ima z>0$; we settle for the weaker result.

\smallskip
2. \eqref{8.11} is not central for a result of this type, but it is true in applications
and simplifies the notation.

\smallskip
3. This result for $\ell=1$ is in \cite{Sim288}. $\ell >1$ has some subtleties, but
the basic strategy we use is that of \cite{Sim288}.
\end{remarks}

We will prove this result in a sequence of lemmas:

\begin{lemma}\lb{L8.4} $\det (f(z))$ is analytic and nonvanishing in $\Omega\equiv
\{z : z\in\bbD, \, \Ima z>0\}$, and $\arg (\det(f(z)))$ can be chosen in that
region to be continuous so that
\begin{equation} \lb{8.13}
0 < \arg (\det(f(z))) <\pi \ell
\end{equation}
\end{lemma}

\begin{proof} In $\Omega$, all matrix elements $\langle\varphi, f(z)\varphi\rangle$
are analytic and have a.e.\ boundary values (since they are scalar Herglotz functions),
so by polarization, $f(z)$ is analytic on $\Omega$ and has a.e.\ boundary values. Thus
$\det (f(z))$ as a polynomial in matrix elements is analytic on $\Omega$.

Consider
\begin{equation} \lb{8.14}
P(\lambda, z) = \det(\lambda\bdone -f(z))
\end{equation}
which is a polynomial in $\lambda$ with analytic coefficients away from the poles of $f$.
It follows, for $z$ near any $z_0$ about which $f$ is analytic, that the roots
$P(\lambda,z)=0$ written as a function of $z$ are analytic functions in $(z-z_0)^{1/k}$
for some $k$ depending on $z_0$. It then follows that near any fixed $z_0$, all roots are
analytic, that is, singularities are isolated.

Pick $x_0\in (0,\veps)$ so that $x_0 f(x_0) >0$, so all eigenvalues $\lambda_1(x_0),
\dots, \lambda_\ell(x_0)$ are in $(0,\infty)$.  Let $z\in\Omega$ be a point about which
all eigenvalues are analytic, and let $\gamma(z)$ be a simple closed path from $x_0$
to $z$ which avoids the discrete set where eigenvalues are not analytic and lie
in $\Omega$ except for $x_0$ with, say, $\gamma(0)=x_0$, $\gamma(1)=z$. By
analytically continuing eigenvalues, we get function $\{\lambda_j(z)\}_{j=1}^\ell$,
so $\lambda_j (z)$ are all the eigenvalues of $f(\gamma(t))$ and $\lambda_j(0)\in
(0,\infty)$. By $\Ima f>0$, $\Ima \lambda_j(z) >0$, so if we define $\arg (\lambda_j
(z))$ with $\arg (\lambda_j(0))=0$, we have
\[
0 < \arg (\lambda_j(z)) <\pi
\]
Thus
\[
\arg (\det(f(z))) = \sum_{j=1}^\ell \arg (\lambda_j(z))
\]
normalized by $\arg (\det (f(x_0)))=0$ obeys \eqref{8.13}.

By analyticity of $\det (f(z))$ and the fact that it is nonvanishing, $\arg(\det(f(z)))$ is
uniquely defined as a continuous function on $\Omega$ with $\lim_{\veps\downarrow 0}
\arg(\det(f(x_0 + i\veps)))=0$. By the above, \eqref{8.13} holds at all points $z$ in
$\Omega$ where all eigenvalues are analytic and so, by continuity and the open mapping
theorem for analytic functions, all points.
\end{proof}

\begin{lemma}\lb{L8.5} Let $a<b$ lie in $(-1,1)$ so that both $a$ and $b$ are neither
a zero nor a pole of $\det (f(z))$. Let $Z(a,b), P(a,b)$ be the number of zeros, poles
of $\det(f(z))$ in $(a,b)$ counting multiplicity. Then
\begin{equation} \lb{8.15}
\abs{Z(a,b)-P(a,b)} \leq\ell
\end{equation}
\end{lemma}

\begin{proof} By the argument principle, $2\pi (Z-P)$ is the change of $\arg (\det(f(z))$
along the circle through $a$ and $b$ centered at $\f12 (a+b)$. By Lemma~\ref{L8.4}, this
is at most $2(\ell\pi)$.
\end{proof}

\begin{lemma}\lb{L8.6} The sets of zeros and poles {\rm{(}}with multiplicity{\rm{)}}
of $\det(f(z))$, including the $\ell$-fold zero at $z=0$, can be written as
$\ell$ subsets $z_j^{(k)}, p_j^{(k)}$ with $k=1, \dots, \ell$ and $-\ti N_k < j < N_k$
{\rm{(}}with $N_k$ and $\ti N_k$ among $1,2,\dots, \infty${\rm{)}} so that
$z_0^{(k)}=0$ and
\begin{equation} \lb{8.16}
z_j^{(k)} < p_j^{(k)} < z_{j+1}^{(k)}
\end{equation}
for all allowed values of $j$.
\end{lemma}

\begin{remarks} 1. If there are infinitely many $z$ in $(-1,0)$ and in $(0,1)$, then
$\ti N_k=N_k=\infty$ for all $k$. The awkwardness requiring $N_k,\ti N_k$ is
only needed if there are finitely many zeros.

\smallskip
2. To avoid notational complexity, we slightly lied if $N_k$ or $\ti N_k$ is finite.
If $N_k$ is finite, $z_j^{(k)}$ runs to $j=N_k$. $p_j^{(k)}$ can then run to either
$N_k$ or $N_k-1$.
\end{remarks}

\begin{proof} Construct $S_1, S_2, \dots, S_\ell$ as follows: Set $z_0^{(1)}=0$. Let
$p_0^{(1)}$ be the first pole larger than $z_0^{(1)}$, $z_1^{(1)}$ the first zero
larger than $p_0^{(1)}$, $p_1^{(1)}$ the next pole, etc. This either continues
indefinitely, in which case we set $N_1=\infty$, or stops because there is no next
zero or pole. Then do the same to the left of $0$, that is, $p_{-1}^{(1)}$ is the
first pole smaller than $z_1^{(1)}$, etc. Clearly, the points in $S_1$ obey \eqref{8.16}.
Now remove the points of $S_1$ (or decrease their multiplicity by $1$) and repeat the
construction (starting with $z_0^{(2)}=0$) to make $S_2, S_3, \dots, S_\ell$.

We claim that after we construct $\ell$ $S_j$'s, we have exhausted all the poles and
zeros. Let us show this is true for $(0,1)$; the argument for $(-1,0)$ is similar
(and since $0$ has multiplicity $\ell$, it is removed after $\ell$ steps).

Suppose $\ti z$ is a zero that is left and it is closer to zero than any leftover zero
or pole. If $\ti z$ lies in some $(p_j^{(k)}, z_{j+1}^{(k)})$, $j=0,1,\dots$, we could
have used it as $z_{j+1}^{(k)}$ so it cannot lie in any such interval. Put differently,
there are only matched zeros and poles in $(0,\ti z)\cap \cup_{j=1}^\ell S_j$.
By the choice of $\ti z$, there are no other poles in $(0,\ti z)$. Thus, for small
$\delta$, the interval $(-\delta, \ti z+\delta)$ has $\ell+1$ extra zeros over poles,
violating Lemma~\ref{L8.5}. So the closest leftover point is not a zero.

Suppose $\ti p$ is a pole that is left and it is closer to zero than any other
leftover zero or pole. As above, $\ti p$ cannot lie in any $(z_j^{(k)}, p_j^{(k)})$,
$j=0,1,\dots$, so there are only matched zeros and poles in $[0,\ti p)\cap\cup_{j=1}^\ell
S_j$. But then, for small $\delta$, $(\delta, \ti p+ \delta)$ has $\ell+1$ extra poles,
violating Lemma~\ref{L8.5}. Thus $\cup_{j=1}^\ell S_j$ includes all zeros and poles.
\end{proof}

\begin{lemma}\lb{L8.7} The limit $B_\infty(z)$ of \eqref{8.12} exists and obeys
conditions {\rm{(i)--(iii)}} of Proposition~\ref{P8.3}{\rm{(c)}}.
\end{lemma}

\begin{proof} Renumber the $p_j^{(k)}$ into a single sequence $p_1, p_2, \dots$, so
$\abs{p_1} \leq \abs{p_2} \leq \cdots$ and let $z_m$ be the correspondingly paired
$z_{j+1}^{(k)}$ (paired to the $p_j^{(k)}$ that is $p$). Since $\{(p_j^{(k)},
z_{j+1}^{(k)})\}_{j=1}^{N_j}$ are disjoint subsets of $(0,1)$ for each fixed $k$,
\[
\sum_{j=1}^\infty \, \abs{z_{j+1}^{(k)} - p_j^{(k)}} =
\sum_{j=1}^\infty z_{j+1}^{(k)} - p_j^{(k)} <1
\]
so we see that
\[
\sum_{j=1}^\infty \abs{z_j-p_j} \leq 2\ell
\]
The existence of $B_\infty$ then follows by Proposition~13.8.2 of \cite{OPUC2}, as do
(i) and (ii).

To get (iii), we note that just taking the zeros and poles in a single $S_j$ yields
a set obeying (13.8.5) and (13.8.6) of \cite{OPUC2}. So, by (13.8.10), that product
has $\arg$ bounded by $2\pi$. The $\ell$-fold product thus obeys \eqref{8.12a}.
\end{proof}

\begin{proof}[Proof of Proposition~\ref{P8.3}] Given Lemma~\ref{L8.7}, the proof
is essentially that of Theorem~13.8.3 of \cite{OPUC2}.
\begin{equation} \lb{8.17}
g(z)\equiv \f{\det(f(z))}{z^\ell B_\infty (z)}
\end{equation}
is analytic and nonvanishing on $\bbD$ with $g(0)>0$ (since $B_\infty (0) >0$).
Moreover, by \eqref{8.13} and \eqref{8.12a},
\begin{equation} \lb{8.18}
\abs{\arg g(z)} \leq 4\pi \ell
\end{equation}
so, by M.~Riesz' theorem, $\log (g(z))\in\cap_{p<\infty} H^p (\partial\bbD)$
from which (a), (b) of the theorem are immediate and (d) follows from the Poisson
representation for $\log (g(z))$ since $\log (\abs{g(e^{i\theta})})=\log
(\abs{\det (f(e^{i\theta}))}$.
\end{proof}

Now we turn to block Jacobi matrices where we obtain:

\begin{theorem}[Nonlocal Step-by-Step Sum Rule for Block Jacobi Matrices]\lb{T8.8}\leavevmode\\
Let $\calJ$ be a block Jacobi matrix with $\sigma_\ess (\calJ)\subset [-2,2]$ and Jacobi parameters
$\{A_n, B_n\}_{n=1}^\infty$. Let $\calJ^{(1)}$ denote this Jacobi matrix with the top
row of blocks and left-most column of blocks removed. Let $m(E), m^{(1)}(E)$ be the $m$-functions
given by \eqref{5.24}. Let $M,M^{(1)}$ be defined on $\bbD$ by
\begin{equation} \lb{8.19}
M(z) =-m(z+z^{-1})
\end{equation}
with poles at $\{p_i\}_{i=1}^N$ where $p_i + p_i^{-1}$ are eigenvalues of $\calJ$.
We repeat each $p_i$ a number of times equal to the multiplicity of the eigenvalues
{\rm{(}}equivalently, the rank of the residue{\rm{)}}. Let $\{z_i\}_{i=1}^{N-1}$
be the corresponding points for $\calJ^{(1)}$. Then
\begin{SL}
\item[{\rm{(a)}}] The Blaschke product, $B_\infty (z)$, defined by the $\{z_i\}
\cup\{p_i\}$ via \eqref{8.12} exists and obeys {\rm{(i)--(iii)}} of
Proposition~\ref{P8.3}{\rm{(c)}}.

\item[{\rm{(b)}}] $M(z)$ and $M^{(1)}(z)$ have limits $M(e^{i\theta})$ and
$M^{(1)}(e^{i\theta})$ as $r\uparrow 1$ for $z=re^{i\theta}$ for a.e.\ $\theta$
in $\partial\bbD$ and
\begin{equation} \lb{8.20}
\log\biggl[ \f{\det (\Ima M)(e^{i\theta})}{\det(\Ima M^{(1)}(e^{i\theta}))}\biggr]
\in\bigcap_{1\leq p <\infty} L^p \biggl( \partial\bbD, \f{d\theta}{2\pi}\biggr)
\end{equation}

\item[{\rm{(c)}}]
\begin{equation} \lb{8.21}
\det\biggl( \f{\abs{A_1} M(z)}{z}\biggr)
 = B_\infty (z) \exp\biggl( \int \f{e^{i\theta}+z}{e^{i\theta}-z} \,
    \log\biggl[ \f{\det (\Ima M(e^{i\theta}))}{\det (\Ima M^{(1)}(e^{i\theta}))}\biggr]
    \frac{d\theta}{4\pi}\biggr)
\end{equation}
\end{SL}
\end{theorem}

\begin{remark} As in the case $\ell=1$, it can happen (although not in examples
where sum rules are finite) that $\det (\Ima M(e^{i\theta}))=\det (\Ima M^{(1)}
(e^{i\theta}))=0$ for $\theta$ in a set of positive measure. (b) and (c) are shorthand for
the more precise
\begin{SL}
\item[(i)] For a.e.\ $\theta$, $\det (\Ima M(e^{i\theta}))=0$ if and only if
$\det (\Ima M^{(1)}(e^{i\theta}))=0$.
\item[(ii)] There is an a.e.\ positive function $g(\theta)$ on $\partial\bbD$, equal
to $\det (\Ima M(e^{i\theta}))/\det (\Ima M^{(1)}(e^{i\theta}))$ when the
ratio is not $0/0$ so that \eqref{8.20} and \eqref{8.21} hold if the formal ratio
is replaced by $g(\theta)$.
\end{SL}
\end{remark}

\begin{proof} Given Proposition~\ref{P8.3}, this is essentially identical to the
proof of Theorem~13.8.4 of \cite{OPUC2} with care given to matrix issues. We begin
by noting that \eqref{5.27} for $n\to n+1$ first implies near $z=0$
\begin{equation} \lb{8.22}
M^{(n+1)}(z)^{-1} = z^{-1} + O(1)
\end{equation}
and then by \eqref{5.24} that
\begin{equation} \lb{8.23}
\biggl( \f{M^{(n)}(z)}{z}\biggr)^{-1} = 1-B_{n+1} z - (A_{n+1}^\dagger A_{n+1} -1)
z^2 + O(z^3)
\end{equation}
Since $M^{(n)}(z)/z$ is near $1$ for $z$ small, we can compute its determinant using
\begin{equation} \lb{8.24}
\det(C) = \exp (\tr(\log(C))
\end{equation}
which holds if $\|C-1\|<1$. Thus
\begin{equation} \lb{8.25}
\begin{split}
\log &\det \biggl( \f{M^{(n)}(z)}z \biggr) \\
&= \tr (B_{n+1})z +
[\tr\{ ([A_{n+1}^\dagger A_{n+1}-1]+ \tfrac12 B_{n+1}^2)\}]
z^2 + O(z^3)
\end{split}
\end{equation}

In addition, \eqref{5.27} implies
\begin{equation} \lb{8.26}
\Ima [M(z)^{-1}] = \Ima (z+z^{-1}) - A_1 \Ima M_1(z) A_1^\dagger
\end{equation}
so at points where $M(z)$ has radial limits (a.e.\ $\theta$, see below),
\begin{equation} \lb{8.27}
-[M(e^{i\theta})^\dagger]^{-1} \Ima M(e^{i\theta}) [M(e^{i\theta})]^{-1} =
-A_1 \Ima M(e^{i\theta}) A_1^\dagger
\end{equation}
which, using (on account of $\det (\abs{C})^2=\det (C^\dagger)\det (C)$)
\[
\abs{\det (A_1)} = \det (\abs{A_1})
\]
yields
\begin{equation} \lb{8.28}
\abs{\det (\abs{A_1} M(e^{i\theta}))}^2 =
\f{\det (\Ima M(e^{i\theta}))}{\det (\Ima M_1 (e^{i\theta}))}
\end{equation}

We now apply Proposition~\ref{P8.3} to $M(z)$ which obeys \eqref{8.10}
(since $\Ima (z+z^{-1})<0$ on $\bbD$ and \eqref{8.19} has a minus sign) and
\eqref{8.11} by \eqref{8.22}.

By Theorem~\ref{T5.6}, our $B_\infty (z)$ here (after perhaps canceling some
zeros and poles) is the $B_\infty (z)$ of Proposition~\ref{P8.3}. (a) and (b)
immediately follow from Proposition~\ref{P8.3}. We get \eqref{8.21} from \eqref{8.12b}
by using \eqref{8.28} (noting \eqref{8.21} has a $1/4\pi$ while \eqref{8.12b} a
$1/2\pi$ on account of the square on the left side of \eqref{8.28}). We also use that
if $c$ is a positive constant,
\begin{equation} \lb{8.29}
\exp\biggl( \int \f{e^{i\theta}+z}{e^{i\theta}-z}\,
\log (c^2) \f{d\theta}{4\pi}\biggr) =c
\end{equation}
\end{proof}

As in \cite{Sim288}, we can get step-by-step $P_2$ (originally in \cite{KS}),
$C_0$, $C_1$ (originally in \cite{SZ}) sum rules immediately from Taylor expansion
of the $\log$ of \eqref{8.21}. We let $\beta_j (\calJ)$ be the numbers in $(-1,1)
\setminus \{0\}$ for which $E_j\equiv \beta_j + \beta_j^{-1}$ are eigenvalues of
$\calJ$ counting multiplicities.

\begin{theorem}[$C_0$, $C_1$, $P_2$ Step-by-Step Sum Rules] \lb{T8.9} $\qquad$ \\
\noindent {\rm{(i)}}
\begin{equation} \lb{8.30}
\begin{split}
\f{1}{4\pi} &\int_0^{2\pi} \log \biggl( \f{\det (\Ima M^{(1)}(e^{i\theta}))}
{\det (\Ima M (e^{i\theta}))} \biggr) d\theta \\
& -\sum_j \log (\abs{\beta_j(\calJ)}) -
\log (\abs{\beta_j (\calJ^{(1)})}) = -\log (\det\abs{A_1})
\end{split}
\end{equation}

\noindent{\rm{(ii)}}
\begin{equation} \lb{8.31}
\begin{split}
-\f{1}{2\pi} &\int_0^{2\pi} \log \biggl( \f{\det (\Ima M^{(1)}(e^{i\theta}))}
{\det (\Ima M (e^{i\theta}))}\biggr) \cos\theta
\\
&+ \sum_j [\beta_j (\calJ)
- (\beta_j^{-1} (\calJ)^{-1})] - [\beta_j (\calJ^{(1)}) - (\beta_j (\calJ^{(1)}))^{-1}]
\equiv \tr (B_1)
\end{split}
\end{equation}

\noindent{\rm{(iii)}}
\begin{equation} \lb{8.32}
\begin{split}
\f{1}{4\pi} &\int \log \biggl( \f{\det (\Ima M^{(1)}(e^{i\theta}))}
{\det (\Ima M (e^{i\theta}))}\biggr) \sin^2(\theta) \phi(\theta) \\
&+ \sum F (E_j(\calJ))
-F(E_j (\calJ^{(1)})) =  \tr \tfrac14 \bigl( B_1^2 + \tfrac12\, G(\abs{A_1}) \bigr)
\end{split}
\end{equation}
where $F$ is given by \eqref{8.2} and $G$ by \eqref{8.3}.
\end{theorem}

\begin{remark} The $E_j(\calJ)$ in $(-\infty, -2)$ and $(2,\infty)$ and $E_j(\calJ^{(1)})$
interlace in the $\ell=1$ case. In the general $\ell$ case, we have at most $\ell$ fewer
eigenvalues of $\calJ^{(1)}$ on any $(-\infty, -E_0)$ or $(E_0,\infty)$ so, as in
Lemma~\ref{L8.6}, one can decompose into $\ell$ interlacing subsets. This and the monotonicity
of functions like $F$ show the eigenvalue sums in \eqref{8.30}--\eqref{8.32} are
conditionally convergent. Similarly, the integrals are always convergent.
\end{remark}

\begin{proof} Apply $\log$ to both sides of \eqref{8.21} and take Taylor coefficients.
The constant term is \eqref{8.30} and the first derivative is \eqref{8.31}. If $L(z)$ is
the $\log$ of the left side and $R(z)$ of the right, then
\[
L(0) + \tfrac12\, L''(0) = R(0) + \tfrac12\, R''(0)
\]
is \eqref{8.32}.
\end{proof}

The proofs of Theorems~\ref{T8.1}--\ref{T8.2} are now identical to those of the scalar
case; see, for example, the discussion of Theorems~13.8.6 and 13.8.8  of \cite{OPUC2}.
In particular, $Z(\calJ)$ and $Q(\calJ)$ (the integral on the right of \eqref{8.4})
are negatives of relative entropies, and so, lower semi-continuous.

\section{Szeg\H{o} and Killip--Simon Theorems When All Gaps Are Open} \lb{s9}

Our goal here is to prove Theorems~\ref{T1.3} and \ref{T1.4}. Our strategy, of course,
will be to translate Theorems~\ref{T8.1} and \ref{T8.2} for $\Delta (J)$ to statements
about $J$. Firstly, we need to relate the a.c.\ part of the matrix measure for $\Delta(J)$ to the
a.c.\ part of the (scalar) measure for $J$. And secondly, to relate $\ell^2$ norms of
coefficients of $\Delta(J)$ to the distance of $J$'s Jacobi parameters to the isospectral torus.
We begin with the first question. Thus we take
\begin{equation} \lb{9.1}
d\mol_J(x)=\omega(x)\, dx + d\mol_{J,\s}(x)
\end{equation}
with $d\mu_{J,\s}$ singular and $\omega$ supported precisely on $\sigma_\ess (J_0)$. By this
assumption and the spectral mapping theorem, $\Delta(J)$ has a.c.\ spectrum precisely
on $[-2,2]$ so the matrix measure for $\Delta(J)$ has the form
\begin{equation} \lb{9.2}
d\mol_{\Delta(J)}(E)=W(E)\, dE + d\mol_{\Delta(J),\s}(E)
\end{equation}

\begin{proposition}\lb{P9.1} Let $J_0$ be a periodic Jacobi matrix with period $p$
and $J$ a Jacobi matrix with Jacobi parameters $\{a_n, b_n\}_{n=1}^\infty$ and measure
$d\mol_J$ of the form \eqref{9.1} with $\omega$ supported on $\sigma_\ess (J_0)$. Let
$\Delta$ be the discriminant for $J_0$ and $W(E)$ the a.c.\ part of the $p\times p$
matrix-valued measure $d\mol_{\Delta(J)}$ associated to $\Delta(J)$ {\rm{(}}so $W$\! is
a $p\times p$ matrix{\rm{)}}. Then for $E\in (-2,2)$ and $\Delta^{-1} (E) =
\{x_1, \dots, x_p\}$,
\begin{equation} \lb{9.3}
\det (W(E)) = \biggl( \prod_{j=1}^{p} a_j^{p-j} \biggr)^{-2}
\biggl( \prod_{j=1}^{p} a^{(0)}_j \biggr)^p
    \biggl( \prod_{j=1}^p \omega(x_j) \biggr)
\end{equation}
\end{proposition}

\begin{proof} In the block Jacobi form, $d\mol_{\Delta(J)}$ has $jk$ matrix element equal
to the spectral measure of the operator $\Delta(J)$ associated to
$\delta_j,\delta_k$, that is, $\int F(x) (d\mol_{\Delta}(x))_{jk}
=\langle \delta_j, F(\Delta (J))\delta_k\rangle$. But $\delta_j=
p_{j-1}(J)\delta_1$. It follows that
\begin{align}
W_{kj}(E) &= \sum_{\ell=1}^p \omega (x_\ell) \bigl( |\Delta'(x_\ell)| \bigr)^{-1} p_{k-1} (x_\ell) p_{j-1} (x_\ell) \lb{9.4}
\end{align}
Note that the factors of $1/\Delta'$ arise from the Jacobian $\tfrac{dE}{dx} = \Delta'(x)$.
We can re-write \eqref{9.4} as $W_{kj}(E) = (M\!AM^t)_{kj}$ where $A$ is the diagonal matrix
\begin{equation} \lb{9.7}
A_{\ell m} = \delta_{\ell m} \omega (x_\ell) \bigl( |\Delta'(x_\ell)| \bigr)^{-1}
\end{equation}
and $M$ is the matrix
\begin{equation} \lb{9.6}
M_{k\ell} = p_{k-1}(x_\ell) \qquad k=1, \dots, p;\, \ell=1, \dots, p
\end{equation}

Next we compute $\det(M)$; $\det(A)$ is easy.  Note that
\begin{equation} \lb{9.9}
p_{k-1} (x_\ell) = \biggl(\, \prod_{j=1}^{k-1} a_j\biggr)^{-1} x_\ell^{k-1} +
\text{lower order}
\end{equation}
Moreover, inductively one sees that the lower order terms can be neglected in the
determinant---they can be removed by subtracting a multiple of rows above (i.e.,
smaller values of $k$). Thus,
\begin{equation} \lb{9.10x}
\det(M) = \biggl[\, \prod_{k=1}^p \biggl(\, \prod_{j=1}^{k-1} a_j\biggr)^{-1}\biggr]
    \det (x_\ell^{k-1})
= \biggl( \prod_{j=1}^{p} a_j^{p-j} \biggr)^{-1} \biggl( \prod_{j
> k} (x_j-x_k) \biggr)
\end{equation}
by the well-known formula for Vandermonde determinants.  This can be simplified further.  The
points $x_j$ are precisely the zeros of the polynomial $\Delta(x)-E$; hence,
invoking \eqref{2.4new},
$$
\Delta(x) - E =  \biggl(\prod_{j=1}^p a^{(0)}_j \biggr)^{-1} \biggl(\prod_{k=1}^p (x-x_k) \biggr)
$$
In this way we discover that
\begin{equation} \lb{9.10new}
\det(M)^2 = \biggl( \prod_{j=1}^{p} a_j^{p-j} \biggr)^{-2} \biggl(
\prod_{j=1}^{p} a^{(0)}_j \biggr)^p \biggl(\prod_{k=1}^p
|\Delta'(x_k)| \biggr)
\end{equation}
Multiplying this by $\det(A)$ gives \eqref{9.3}.
\end{proof}

\begin{corollary}\lb{C9.2&3} If $J_0$ has all gaps open and $\alpha >-1$, then
\begin{equation} \lb{9.11}
\int_{-2}^2 (4-E^2)^\alpha \abs{\log \det (W(E))}\, dE <\infty
\end{equation}
if and only if
\begin{equation} \lb{9.12}
\int_{\sigma_\ess(J_0)}\dist (x,\bbR\setminus\sigma_\ess (J_0))^\alpha
\abs{\log \omega(x)}\, dx <\infty
\end{equation}
When $\alpha=-\frac12$, the same conclusion holds even if some gaps are closed.
\end{corollary}

\begin{remark} Since $\alpha >-1$, $(4-E^2)^\alpha$ (resp., $\dist (\dots)^\alpha$) are in
$L^p$ for some $p>1$, so the $\log_+ (\quad)$ is always integrable and these conditions are
equivalent to the integral without $\abs{\cdot}$ being larger than $-\infty$.
\end{remark}

\begin{proof}  Changing variables via $E=\Delta(x)$ and applying Proposition~\ref{P9.1} shows
that \eqref{9.11} holds if and only if
\begin{equation}
  \int_{\sigma_\ess(J_0)} |\log \omega(x) | \,(4-\Delta(x)^2)^\alpha \,|\Delta'(x)|\,dx  < \infty
\end{equation}

If all gaps are open then $|\Delta'(x)|$ is strictly positive on $\sigma_\ess(J_0)$, while $4-\Delta(x)^2$
is a polynomial with a simple zero at each band edge (and no others).  This proves the first claim.

At a closed gap, $4-\Delta(x)^2$ has a double zero and $\Delta'(x)$ a simple zero.
When $\alpha=-\frac12$,  these cancel exactly.
\end{proof}

Next we turn to the $\ell^2$ issue. Given any two-sided periodic
matrix $\ti J$ with Jacobi parameters $\{a_n, b_n\}_{n=1}^p$ and
fixed periodic $J_0$, let $B_{J_0}(\ti J), A_{J_0}(\ti J)$ be the
constant $p\times p$ blocks in $\Delta_{J_0} (\ti J)$. We are
heading towards showing that $\|B_{J_0}(\ti J)\|_2^2
+\|A_{J_0}(\ti J)-\bdone\|_2^2$ is comparable to $\dist
((a_n,b_n)_{n=1}^p,\calT_{J_0})$. This will be the key to showing
$\ell^2$ tails in the matrix pieces of
$\Delta_{J_0}(J)-S^p-S^{-p}$ for general $J$ is equivalent to
\eqref{1.52}.\marginpar{changes}

The crucial fact will be that the polynomial coefficients of
$\Delta_{J_0}({\ti J})-\Delta_{J_0}(J_0)$ are comparable to $\dist
((a_n,b_n)_{n=1}^p, \calT_{J_0})$.\marginpar{changes} For this we
need the following, which is a simple application of the implicit
function theorem and compactness:

\begin{lemma}\lb{L9.3} Let $F$ be a $C^\infty$ map of an open set $U\subset\bbR^n$
to $\bbR^\ell$ with $\ell <n$. Suppose $\calT =F^{-1}(y_0)$ is a smooth manifold of
dimension $n-\ell$ and compact for some $y_0\in\bbR^\ell$, and
\begin{equation} \lb{9.13}
\rank ((\nabla F)(x_0))=\ell
\end{equation}
for all $x_0\in\calT$. Then for any compact neighborhood, $K$, of
$\calT$, there are $c_K,d_K \in (0,\infty)$ so for all $x\in K$,
\begin{equation} \lb{9.14}
c_K\abs{F(x)-y_0} \leq \dist (x,\calT)\leq d_K\abs{F(x)-y_0}
\end{equation}
\end{lemma}

One can restate \eqref{9.13} in a more illuminated way in terms of the components
$F_1, \dots, F_\ell$ of $\calT$. Of course, $\nabla F_j(x_0)$ is orthogonal to
$\calT$ at $x_0$. The condition \eqref{9.13} is equivalent to saying that
$\{\nabla _j F(x_0)\}_{j=1}^\ell$ span the normal bundle to $\calT$. This is
equivalent to saying they are linearly independent. Notice that if $J_0$ has all
gaps open, $\calT_{J_0}$ is of dimension $p-1=2p-(p+1)$ and $\Delta_{J_0}$ is a
polynomial of degree $p$, hence with $p+1$ coefficients. Thus the following shows
we can use Lemma~\ref{L9.3}:

\begin{theorem}\lb{T9.4} Suppose all gaps are open for some periodic $J_0$. Then
at any point in $\calT_{J_0}$, the gradients of the derivatives of the coefficients
of $\Delta_J$ span the normal bundle of $\calT_{J_0}$ in $\bbR^{2p}$.
\end{theorem}

\begin{proof} $\Delta_{J_0}$ has the form
\[
\Delta_{J_0}(x)=(a_1 \dots a_p)^{-1} \prod_{j=1}^p (x-\lambda_j) =
\sum_{j=0}^p c_j x^j
\]
where $\lambda_j$ are the roots. The coefficients thus obey
\begin{align}
c_p^{-1} &= a_1 \dots a_p \lb{9.15} \\
c_\ell c_p^{-1} &= \sum_{1\leq k_1 \leq \cdots\leq k_{p-\ell}\leq p}
\lambda_{k_1} \dots \lambda_{k_{1-\ell}} \equiv s_{k-\ell} \qquad \ell <p \lb{9.16}
\end{align}
It is well known that if
\begin{equation} \lb{9.17}
t_\ell =\sum_{j=1}^p \lambda_j^\ell
\end{equation}
then $t_\ell$ is $\ell s_\ell$ plus a polynomial in $\{s_j\}_{j=1}^{\ell-1}$,
so $\{\nabla t_j\}_{j=1}^\ell$ and $\{\nabla s_j\}_{j=1}^\ell$ span
the same space. It follows that we need only show the gradients of $c_p^{-1}$
and $t_\ell$ span the normal bundle of $\calT_{J_0}$.

Let
\begin{equation} \lb{9.18}
\calM_0=\{(a_n,b_n)_{\ell=1}^p : a_1 a_2\dots a_p = a_1^{(0)} a_2^{(0)}\dots
a_p^{(0)}; \, b_1+\cdots + b_p =b_1^{(0)}+\cdots + b_p^{(0)}\}
\end{equation}
We know $\calM_0\supset\calT_{J_0}$. Clearly, $\nabla c_p^{-1}$ and $\nabla t_1$ span
the normal bundle to $\calM_0$ since $t_1=\sum_{j=1}^p b_j$ (see \eqref{2.4new}). Thus
we need only show the projections of $\{\nabla t_\ell\}_{\ell=2}^p$ into the
tangent space of $\calM_0$ span the normal bundle of $\calT_{J_0}$ in $\calM_0$.

Studies of the Toda flows show that $\calM_0$ is a symplectic manifold with
$\{t_\ell\}_{\ell=2}^p$ Poisson commuting. Since the symplectic form on $\calM_0$
is nondegenerate, to say $\{\nabla t_j\}_{j=2}^p$ span the normal bundle is
the same as saying that the Hamiltonian flows generated by $\{t_j\}_{j=2}^p$ span
the tangent bundle of $\calT_{J_0}$, or equivalently, given $\dim (\calT_{J_0})=p-1$,
that these Hamiltonian flows are independent.

This independence is a theorem of van Moerbeke \cite[Theorem~5.2]{vMoer} or \cite{SimonNew}.
\end{proof}

\begin{lemma}\lb{L9.5} Let $\chi_k$ be the projection onto the $k$-dimensional space
spanned by $\{\delta_j\}_{j=1}^k$. For any compact subset, $K$, of
period $p$ Jacobi matrices, there exist constants $c_K$ and $d_K$
in $(0,\infty)$ so for all $J\in K$,
\begin{equation} \lb{9.19}
c_K \biggl\|\, \sum_{\ell=0}^p \alpha_\ell J^\ell
\chi_{p+1}\biggr\|_2 \leq \biggl( \sum_{\ell=0}^p\,
\abs{\alpha_\ell}^2 \biggr)^{\frac12} \leq d_K \biggl\|\,
\sum_{\ell=0}^p \alpha_\ell J^\ell \chi_{p+1}\biggr\|_2
\end{equation}

\end{lemma}

\begin{proof} $\{J^\ell \chi_{p+1}\}_{\ell=0}^p$ are independent since $J^\ell$ has
strictly positive elements in the $\ell$-th diagonal and $\{J^k\}_{k<\ell}$ only
has zero elements there.
Hence, the matrix
\begin{equation} \lb{9.20b}
\left. \tr (\chi_{p+1} J^\ell J^k \chi_{p+1})\right|_{\ell,k=0,
\dots, p}
\end{equation}
is strictly positive so \eqref{9.19} holds for each fixed $J$. The optimal constants are
clearly continuous so uniformly bounded above and below on $K$.
\end{proof}

\begin{proposition}\lb{P9.6} Let $J_0$ be a periodic Jacobi matrix with all gaps open. For any
compact neighborhood $K$ of $\calT_{J_0}$ in $(0,\infty)^p\times\bbR^p$, there are constants
$c_K$ and $d_K$ in $(0,\infty)$ so that for all $J\in K$,
\begin{align*}
c_K (\|A_{J_0}(J)-\bdone\|^2 + \|B_{J_0}(J)\|^2)^{1/2}
&\leq  \dist(J,\calT_{J_0})  \\
&\leq d_K (\|A_{J_0}(J)-\bdone\|^2 + \|B_{J_0}(J)\|^2)^{1/2}
\end{align*}
\end{proposition}

\begin{proof} We have that
\begin{align} \lb{9.21}
2\|A_{J_0}(J)-\bdone\|^2 + \|B_{J_0}(J)\|^2 & \le
\|[\Delta_{J_0}(J)-(S^p + S^{-p})]\chi_p\|^2 \\
\nonumber & \le 4\|A_{J_0}(J)-\bdone\|^2 + 2\|B_{J_0}(J)\|^2
\end{align}
But by the magic formula,
\begin{equation} \lb{9.22}
\Delta_{J_0}(J_0)=S^p + S^{-p}
\end{equation}
so
\begin{equation} \lb{9.23}
[\Delta_{J_0}(J)-(S^p + S^{-p})]\chi_{p+1} = \sum_{\ell=0}^p
c_\ell J^\ell \chi_{p+1}
\end{equation}
where $c_\ell$ is the difference of coefficients for $J$ and $J_0$. By Lemma~\ref{L9.5},
\begin{equation} \lb{9.24}
\|[\Delta_{J_0}(J)-(S^p + S^{-p})]\chi_{p+1}\|^2
\sim\sum_{\ell=0}^p \, \abs{c_\ell}^2
\end{equation}
where $\sim$ means the ratio is bounded above and away from zero on compact subsets.

By Lemma~\ref{L9.3} and Theorem~\ref{T9.4},
\begin{equation} \lb{9.25}
\sum_{\ell=0}^p\, \abs{c_\ell}^2 \sim \dist (J,\calT_{J_0})^2
\end{equation}
Combining this with \eqref{9.24} proves the proposition.
\end{proof}

Now we take a general $J$ not periodic and form $\Delta_{J_0}(J)$ which is a one-sided
block Jacobi matrix with block elements $A_{n,J_0}(J), B_{n,J_0}(J)$.

\begin{lemma}\lb{L9.7} $\Delta_{J_0}(J)_{k\ell}$ for $k\leq \ell$ depends only on
$\{b_j\}_{j=k-\alpha}^{\ell+\alpha}$ and $\{a_j\}_{j=k-\alpha}^{\ell+\alpha -1}$
where $\alpha = \lfloor\f12 (p-(\ell-k))\rfloor$ is the
greatest integer less than or equal to $\f12[p-(\ell-k)]$.
\end{lemma}

\begin{proof} Each factor of $J$ changes index by at most one. In order to get from
$k$ to $\ell$, $\ell-k$ steps are needed. The remainder cannot go below $\ell-\alpha$
or above $k+\alpha$ and get back to $k$ in $p$ steps.
\end{proof}

\begin{lemma}\lb{L9.8} Let $J$ have Jacobi parameters $\{a_n,b_n\}_{n=1}^\infty$.
Let $\ti J$ be periodic with period $p$ and suppose $b_n=\ti b_n$ for $kp-p\leq n \leq
kp+2p$ and $a_n=\ti a_n$ for $kp-p\leq n\leq kp + 2p-1$. Then
\[
A_{k,J_0}(J)=A_{J_0}(\ti J) \qquad
B_{k,J_0}(J)=B_{J_0}(\ti J)
\]
\end{lemma}

\begin{proof} Immediate from Lemma~\ref{L9.7}.
\end{proof}

\begin{lemma}\lb{L9.9} Let $k\leq \ell$ and $\alpha = [\f12(p-(\ell-k))]$.
For any two $J$ and $\ti J$ and any $K$, there is $C_K$ so that
\begin{equation} \lb{9.26}
\abs{\Delta_{J_0}(J)_{k\ell} -\Delta_{J_0}(\ti J)_{k\ell}} \leq
C_K \sup_{k-\alpha\leq j\leq \ell+\alpha} [\abs{b_j-\ti b_j} +
\abs{a_j-\ti a_j}]
\end{equation}
so long as
\begin{equation} \lb{9.27}
\sup [\abs{b_j} + \abs{\ti b_j} + \abs{a_j} + \abs{\ti a_j}]\leq K
\end{equation}
\end{lemma}

\begin{proof} Immediate from Lemma~\ref{L9.7} and the fact that $\Delta_{J_0}$ has
matrix elements that are fixed (given $J_0$) polynomials in $a$'s and $b$'s.
\end{proof}

\begin{lemma}\lb{l.new1}
{\rm (a)} For any Jacobi matrix, $J$, and $\ell = 1,2,\ldots$, $m
= 1,2,\ldots$,
\begin{equation}\lb{f.new1a}
(J^\ell)_{m \, m+l} = a_m a_{m+1} \cdots a_{m+\ell-1}
\end{equation}
and for $\ell = 2,3,\ldots$, $m = 1,2,\ldots$,
\begin{equation}\lb{f.new1b}
(J^\ell)_{m \, m+\ell-1} = a_m \cdots a_{m+\ell-2} \Bigl(
\sum_{j=0}^{\ell-1} b_{m+j} \Bigr)
\end{equation}
{\rm (b)} For $J_0$ periodic of period $p \ge 2$ and $m =
1,2,\ldots$,
\begin{align}
\lb{f.new1c} \Delta_{J_0}(J)_{m \, m+p} & = \frac{a_m \cdots a_{m+p-1}}{a^{(0)}_m \cdots a^{(0)}_{m+p-1}} \\
\lb{f.new1d} \Delta_{J_0}(J)_{m \, m+p-1} & = \bigl( a^{(0)}_m
\cdots a^{(0)}_{m+p-1} \bigr)^{-1} (a_m \cdots a_{m+p-2}) \Bigl(
\sum_{j=0}^{p-1} (b_{m+j} - b^{(0)}_{m+j}) \Bigr)
\end{align}
\end{lemma}

\begin{proof}
(a) Since $J$ changes index by at most one,
$$
(J^\ell)_{m \, m+\ell} = (J_{m \, m+1}) \cdots (J_{m+\ell-1 \,
m+\ell})
$$
proving \eqref{f.new1a}, while
$$
(J^\ell)_{m \, m+\ell-1} = \sum_{j=0}^{\ell-1} (J^j)_{m \, m+j} \,
J_{m+j \, m+j} \, (J^{\ell+j-1})_{m+j \, m+\ell-1}
$$
which, given \eqref{f.new1a}, proves \eqref{f.new1b}.

(b) By \eqref{2.4new},
$$
\Delta_{J_0}(J) = \bigl( a^{(0)}_1 \cdots a^{(0)}_p \bigr)^{-1}
\left[ J^p - \sum_{j=0}^{p-1} b^{(0)}_{j+1} J^{p-1} + O(J^{p-2})
\right]
$$
which, given (a), $(J^{p-k})_{m \, m+p} = (J^{p-k})_{m \, m+p-1} =
0$ if $k = 2,3,\ldots$, and the periodicity of $a^{(0)}$ and
$b^{(0)}$ yields \eqref{f.new1c} and \eqref{f.new1d}.
\end{proof}

\begin{lemma}\lb{l.new2}
Suppose that $\Delta_{J_0}(J) -S^p - S^{-p}$ is a Hilbert--Schmidt
operator on $\ell^2(\{0,1,2,\ldots\})$. Then,
\begin{align}
\lb{f.new2a} \sum_n (a_n a_{n+1} \cdots a_{n+p-1} - a^{(0)}_n
a^{(0)}_{n+1} & \cdots a^{(0)}_{n+p-1})^2 < \infty \\
\lb{f.new2b} \sum_n \left( \sum_{j=0}^{p-1} (b_{n+j} - b^{(0)}_{n+j}) \right)^2 & < \infty \\
\lb{f.new2c} \sum_n (a_{n+p} - a_n)^2 < \infty & \\
\lb{f.new2d} \sum_n (b_{n+p} - b_n)^2 < \infty &
\end{align}
\end{lemma}

\begin{proof}
For a Hilbert--Schmidt operator, any subset of matrix elements
lies in $\ell^2$, so by \eqref{f.new1c},
$$
\sum_n \left| a_n \cdots a_{n+p-1} \bigl( a^{(0)}_n \cdots
a^{(0)}_{n+p-1} \bigr)^{-1} - 1 \right|^2 < \infty
$$
which, given that $a^{(0)}_n \cdots a^{(0)}_{n+p-1}$ is
$n$-independent, implies \eqref{f.new2a}.

Similarly, \eqref{f.new1d} implies \eqref{f.new2b} if we note that
$\{a_j\}$ bounded and $a_n \cdots a_{n+p-1} \to a^{(0)}_1 \cdots
a^{(0)}_p > 0$ implies $\inf a_j > 0$, so
$$
\inf_m  \bigl( a^{(0)}_m \cdots a^{(0)}_{m+p-1} \bigr)^{-1} (a_m
\cdots a_{m+p-2}) > 0
$$

Since the difference of $\ell^2$ sequences is $\ell^2$,
\eqref{f.new2a} implies (since $a^{(0)}_n$ is periodic)
$$
\sum_n (a_{n+p} - a_n)^2 (a_{n+1} \cdots a_{n+p-1})^2 < \infty
$$
which, given that $\inf a_j > 0$, implies \eqref{f.new2c}.

Similarly, since
$$
\sum_{j=0}^{p-1} (b_{n+1+j} - b_{n+j}) = b_{n+p} - b_n
$$
\eqref{f.new2b} implies \eqref{f.new2d}.
\end{proof}

Our next preliminary is to relate $A\in\calL$ to
\begin{equation} \lb{9.36}
\abs{A} = \sqrt{\displaystyle A^\dagger_{ }A}
\end{equation}

\begin{proposition}\lb{P9.12} The map $A\mapsto\abs{A}$ from $\calL$ to positive definite
matrices is a diffeomorphism.  In particular, for $A$'s in $\calL$ with
$\|A-1\| <\f12$, there exist constants $C_1$ and $C_2$ so that
\begin{equation} \lb{9.37}
C_1 \|A-1\|_2 \leq \|\,\abs{A}-1\|_2 \leq C_2 \|A-1\|_2
\end{equation}
\end{proposition}

\begin{proof} By \eqref{9.36}, $A\mapsto\abs{A}$ is a smooth map.

The inverse map (strictly $|A|^2\mapsto A$) is known as the Cholesky factorization; see \cite{GvL,Wat}.
Given $B>0$, apply the Gram--Schmidt procedure to the (linearly independent) columns of $B$ working from right
to left.  This gives a factorization $B=QA$ with $Q$ unitary and $A\in\calL$.
Note that $|A|=B$ and that because the columns of $B$ are linearly dependent, $B\mapsto A$ is also a smooth map.
\end{proof}

\begin{theorem}\lb{T9.11} Let $J_0$ be a two-sided $p$-periodic Jacobi matrix with all gaps open and let $\Delta_{J_0}$
denote its discriminant.  For a Jacobi matrix with parameters $(a_n,b_n)$, the following are equivalent:
\begin{SSL}
\item $\Delta_{J_0}(J) -S^p - S^{-p}$ is a Hilbert--Schmidt operator on $\ell^2(\{0,1,2,\ldots\})$.
\item $\sum_n \tr\{B_n^2 + |A_n-1|^2\} < \infty$.
\item $\sum_n \tr\{B_n^2 + (|A_n|-1)^2\} < \infty$.
\item $\sum_n \tr\{B_n^2 + G(|A_n|)\} < \infty$.
\item $\sum_m d_m((a,b),\calT_{J_0})^2 < \infty$.
\item $\sum_m \tilde d_m((a,b),\calT_{J_0})^2 < \infty$.
\end{SSL}
Here we have adopted the abbreviations $A_n:=A_{n,J_0}(J)$ and $B_n:=B_{n,J_0}(J)$
\end{theorem}

\begin{proof} (i)$\Leftrightarrow$(ii) amounts to the definition of the Hilbert--Schmidt norm.

(ii)$\Leftrightarrow$(iii) follows from Proposition~\ref{P9.12}.

(iii)$\Leftrightarrow$(iv) Notice that $G$, defined in \eqref{8.3}, obeys
$$
c_\veps (x-1)^2 \leq G(x) \leq c'_\veps (x-1)^2 \qquad \forall\ x\in (\veps, \veps^{-1})
$$
Applying this to the eigenvalues of $|A_n|$ yields this equivalence.

(v)$\Leftrightarrow$(vi) is the $q=2$ case of Proposition~\ref{P2.5}.

(vi)$\Rightarrow$(i) By Lemma~\ref{L9.7}, each matrix element of
$\Delta_{J_0}(J) -S^p - S^{-p}$ is a smooth function of $p$
consecutive pairs $(a_n,b_n)$; moreover, by the magic formula, all
of these smooth functions vanish if the corresponding $p$-tuple
belongs to $\calT_{J_0}$.  The implication now follows from the
fact that smooth functions are Lipschitz.

(i)$\Rightarrow$(vi) Define $J^{(k)}$ to be the $p$-periodic
Jacobi matrix that equals $J$ on block $k$, that is,
\begin{equation} \lb{9.31}
b_\ell^{(k)} =b_{kp+\ell} \qquad a_\ell^{(k)}=a_{kp+\ell}
\end{equation}
for $\ell=1,2,\dots, p$.  Obviously, $J^{(k)}=J$ on block $k$ and,
by \eqref{f.new2c} and \eqref{f.new2d}, the difference on blocks
$k-1$ and $k+1$ are in $\ell^2$, that is,
\begin{equation} \lb{9.32}
\sum_k \Bigl[\, \sup_{(k-1)p\leq j\leq (k+2)p-1}\, \abs{b_j-b_j^{(k)}} +
\abs{a_j -a_j^{(k)}}\;\Bigr]^2 <\infty
\end{equation}
Together with Lemma~\ref{L9.9} and Proposition~\ref{P9.6}, this
implies
\begin{equation}\lb{f.tnewa}
\sum_k \tilde d_{kp} (J^{(k)} , \calT_{J_0})^2 < \infty
\end{equation}
On the other hand, \eqref{9.32} implies that for $j = 1,\ldots,p$,
\begin{equation}\lb{f.tnewb}
\sum_k \tilde d_{kp+j} (J,J^{(k)})^2 < \infty
\end{equation}
By the triangle inequality,
$$
\tilde d_{kp+j} (J , \calT_{J_0})^2 \le 2 d_{kp+j} (J,J^{(k)})^2 +
2 \tilde d_{kp} (J^{(k)} , \calT_{J_0})^2
$$
so \eqref{f.tnewa} and \eqref{f.tnewb} imply (vi).
\end{proof}

\begin{proof}[Proof of Theorem~\ref{T1.4}]  We will refer to the three statements (i)--(iii) of Theorem~\ref{T1.4}
simply by their numbers.

Suppose first that
\begin{equation}\label{1.52again}
d_m\bigl( (a,b),\calT_{J_0})\bigr)\in\ell^2
\end{equation}
Then (i) holds by Theorem~\ref{T1.1}.  Moreover, by Theorem~\ref{T9.11} and the hypothesis that all gaps are open,
the RHS of \eqref{8.4} is finite. Therefore, the LHS is finite.  Next
we use this fact to prove (ii) and (iii).

As $\Delta'$ is nonvanishing at all band edges,
\begin{equation} \lb{9.38}
\sum_j F(\Delta(E_j)) <\infty \iff \sum_{j=1}^N \dist \bigl(E_j,\sigma_\ess (J)\bigr)^{3/2} <\infty
\end{equation}
which verifies (iii). By Corollary~\ref{C9.2&3},
\begin{equation} \lb{9.39}
\text{Leftmost term in \eqref{8.4}} <\infty \iff \text{(ii) holds}
\end{equation}
This completes the proof of (i)--(iii).

Conversely, if (i)--(iii) hold, then by \eqref{9.38}, \eqref{9.39},
and \eqref{8.4}, we see that the RHS of \eqref{8.4} is finite.  By Theorem~\ref{T9.11}, this implies \eqref{1.52again}.
\end{proof}

\begin{proof}[Proof of Theorem~\ref{T1.3}] Let $\beta_j$ be the $\beta$'s associated to
$\Delta(J)$, that is, $\abs{\beta_j} >1$, $\beta_j + \beta_j^{-1} =E_j$ with $E_j$ the
eigenvalues of $\Delta(J)$ in $\bbR\setminus [-2,2]$. Then $\log \abs{\beta_j}\sim
\abs{\beta_j}-1$ as $\beta\to\pm 1$ small and $\abs{\beta_j}-1 \sim (\abs{E_j}-2)^{1/2}$.
Therefore,
\begin{equation} \lb{9.40}
\text{\eqref{8.6}} <\infty \iff \sum_j (\abs{E_j}-2)^{1/2}
\iff \text{\eqref{1.49}}
\end{equation}

By Corollary~\ref{C9.2&3},
\begin{equation} \lb{9.41}
\text{\eqref{8.5}} <\infty \iff \text{\eqref{1.50}}
\end{equation}

Finally, if $\{A_n, B_n\}_{n=1}^\infty$ are the $p\times p$ blocks in $\Delta(J)$,
we have $A_n =U\abs{A_n}$ for some $U$ with $\abs{\det(U)}=1$, so
\begin{equation} \lb{9.42}
\det\abs{A_n} =\abs{\det(A_n)} = \prod_{j=(n-1)p+1}^{np}
\prod_{k=j}^{j+p=1} \biggl[ \f{a_k}{a_k^{(0)}}\biggr]
\end{equation}
by \eqref{2.2} and \eqref{2.4new}. Thus
\[
\sum_{n=1}^N \log (\det(\abs{A_n})) - p\sum_{k=1}^{Np} \log \biggl[\f{a_k}{a_k^{(0)}}\biggr]
\]
is bounded. Thus \eqref{1.51} is equivalent to \eqref{8.8}.

By Theorem~\ref{T8.2}, we see that when \eqref{1.49} holds, then \eqref{1.50} $\Leftrightarrow$
\eqref{1.51}, and if they hold, \eqref{1.58a} has a limit.

Moreover, if they hold, the hypotheses of Theorem~\ref{T1.4} hold, so \eqref{1.52} is true.
That \eqref{1.60a} holds is a theorem of Peherstorfer--Yuditskii \cite{PY}; see also the remark below.
\end{proof}

In the remainder of this section, we will describe an alternate approach to proving \eqref{1.60a}; one based
on combining the magic formula with Theorem~\ref{T5.8}.  Unfortunately, because of the strong hypothesis on
the discrete spectrum that appears in this theorem, we will not recover the full formulation from Theorem~\ref{T1.3}.

Let $\tilde \calJ$ denote the (unique) type~2 block Jacobi matrix that is equivalent
to $\calJ=\Delta_{J_0}(J)$, which is of type~3.  Further, let us use $A_j$ and $\ti A_j$ to denote the
off-diagonal block entries of $\calJ$ and $\ti\calJ$, respectively.

If we strengthen the hypothesis \eqref{1.49} to finiteness of the discrete spectrum (i.e., finiteness
of the set $\sigma(\calJ)\setminus\sigma_\ess(\calJ)$), then Theorem~\ref{T5.8} shows that \eqref{1.50}
implies the convergence of the product $\ti A_1 \cdots \ti A_n$ as $n\to\infty$.  In view of \eqref{5.7b}
and Proposition~\ref{P9.12}, this convergence is inherited by the product $A_1 \cdots A_n$.  Thus, it remains
only to connect the convergence of this matrix product to the behavior of the sequences of parameters.  This
is the job of the next lemma.  In the original version of this paper, it was only proved that the sequence $\{a_n\}$
was asymptotic to a fixed periodic sequence.  The argument for the sequence $\{b_n\}$ was provided by one of the
referees; we are most grateful for this.

\begin{lemma}
Let $J_0$ be a $p$-periodic two-sided Jacobi matrix and let $\Delta=\Delta_{J_0}$ denote its discriminant.
Let $J$ be a one-sided Jacobi matrix with parameters $\{a_n,b_n\}$.  Suppose the product $A_n\cdots A_1$
converges to a non-singular matrix as $n\to\infty$.  Here $A_n$ and $B_n$ denote the $p\times p$ block entries of $\Delta(J)$.
Then the parameters of $J$ asymptotically converge to fixed periodic parameters in the sense of \eqref{1.60a}.
\end{lemma}

\begin{proof}
By applying the same affine transformation (i.e., $x\mapsto \alpha x+\beta$) to both $J$ and $J_0$, we may assume that
the discriminant of $J_0$ takes the form $\Delta(x)=x^p+O(x^{p-2})$.  This will significantly simplify some
of the formulae that follow.  Note also that this transformation makes $a^{(0)}_1\cdots a^{(0)}_p=1$ and
$b^{(0)}_1 + \cdots + b^{(0)}_p=0$.

Let $A_n$ and $B_n$ denote the $p\times p$ block entries of $\Delta(J)$.  Then
\begin{align*}
(A_n)_{k,k}   &= a_{p(n-1)+k} \cdots a_{pn+k-1} \\
(A_n)_{k+1,k} &= a_{p(n-1)+k+1} \cdots a_{pn+k-1} \bigl[b_{p(n-1)+k+1} + \cdots + b_{pn+k}\bigr]
\end{align*}
as can be read off from Lemma~\ref{l.new1}.  Using this and the lower-triangular structure of the
matrices $A_j$, one may quickly deduce
\begin{align}
(A_1\cdots A_n)_{k,k}   &= \prod_{j=k}^{pn+k-1} a_j \label{AprodDiag}\\
(A_1\cdots A_n)_{k+1,k}
    &= \sum_{r=1}^{n} (A_1)_{k+1,k+1} \cdots (A_{r-1})_{k+1,k+1} (A_r)_{k+1,k} (A_{r+1})_{k,k} \cdots (A_n)_{k,k} \notag\\
&= \Biggl( \prod_{j=k+1}^{pn+k-1} a_j \Biggr) \sum_{j=k+1}^{pn+k} b_j \label{AprodOffDiag}
\end{align}

To see that the sequence $n\mapsto a_{pn+k}$ converges for each fixed $k\in\{1,\ldots,p\}$, one need only take ratios of \eqref{AprodDiag}
for consecutive values of $k$ and the same $n$ (and also for $(n,k=p)$ and $(n+1,k=1)$), then send $n\to\infty$.

For the parameters $b_n$, one may proceed in a similar fashion: For example, when $2\leq k \leq p-1$, the fact that
$$
a_{k} \frac{(A_1\cdots A_n)_{k+1,k}}{(A_1\cdots A_n)_{k,k}} - a_{k-1} \frac{(A_1\cdots A_n)_{k,k-1}}{(A_1\cdots A_n)_{k-1,k-1}}
    = b_{pn+k} - b_{k}
$$
shows us that $b_{pn+k}$ converges as $n\to\infty$.
\end{proof}

\section{Szeg\H{o} and Killip--Simon Theorems When Some Gaps Are Closed} \lb{s10}

Here we want to examine what might replace Theorems~\ref{T1.3} and \ref{T1.4} if $J_0$
is periodic but with some closed gaps. The Szeg\H{o}-type theorem is almost the same
as Theorem~\ref{T1.3}:

\begin{theorem}\lb{T10.1} Let $J_0$ be any two-sided periodic Jacobi matrix with Jacobi
parameters $\{a_n^{(0)}, b_n^{(0)}\}_{n=-\infty}^\infty$, and $J$ a one-sided
Jacobi matrix with Jacobi parameters $\{a_n, b_n\}_{n=1}^\infty$ and spectral measure
$d\mol$. Suppose that \eqref{1.47} holds, and that
\begin{equation} \lb{10.1}
\sum_{m=1}^N \dist (E_m,\sigma_\ess (J))^{1/2} <\infty
\end{equation}
if $\{E_m\}_{m=1}^N$ is a labelling of the eigenvalues of $J$ outside $\sigma_\ess (J)$. Then
\begin{equation} \lb{10.2}
-\int_{\sigma_\ess (J_0)} \log \biggl( \f{d\mol_\ac}{dx}\biggr) \dist (x, \bbR\setminus
\sigma_\ess (J_0))^{-1/2}\, dx <\infty
\end{equation}
implies
\begin{equation}\lb{10.3}
\lim\biggl(\, \sum_{j=1}^{pN} \log \biggl( \f{a_j}{a_j^{(0)}}\biggr)\biggr)
\end{equation}
exists and lies in $(-\infty,\infty)$. Conversely, \eqref{10.2} holds so long as
\begin{equation} \lb{10.4}
\limsup \biggl(\, \sum_{j=1}^N \log \biggl( \f{a_j}{a_j^{(0)}}\biggr)\biggr) >-\infty
\end{equation}
and then the limit in \eqref{10.3} exists and lies in $(-\infty,\infty)$.

Moreover, if \eqref{10.2} or \eqref{10.4} holds, then there is $J_1\in\calT_{J_0}$, so
\begin{equation}\lb{10.5}
d_m (J,J_1)\to 0
\end{equation}
\end{theorem}

\begin{remark} All that is missing is \eqref{1.52} which we do not claim. However,
since \eqref{10.1}/\eqref{10.2} imply (i)--(iii) of Theorem~\ref{T10.3} below,
we have \eqref{10.7}.
\end{remark}

\begin{proof} As noted, even with closed gaps, \eqref{1.50} is equivalent to \eqref{9.11}
(see Corollary~\ref{C9.2&3}). Once  one notes this, the proof of Theorem~\ref{T1.3} provides
all the results stated as Theorem~\ref{T10.1}.
\end{proof}

Theorem~\ref{T1.4} used open gaps in two ways. First, in the translation of a matrix
pseudo-Szeg\H{o} condition, \eqref{9.11} with $\alpha=\f12$ to the original spectral
measure of $J$, and second, translating a Hilbert--Schmidt bound on $\Delta (J)-S^p-S^{-p}$
to $\ell^2$ approach to the isospectral torus. The second issue can be finessed if we leave
things as a Hilbert--Schmidt condition, which reduces to a sum of translates of an
explicit positive polynomial in the $a_n$'s and $b_n$'s being finite. As for translating
\eqref{9.11} with $\alpha=\f12$, the argument that proved Corollary~\ref{C9.2&3} translates
immediately to

\begin{lemma}\lb{L10.2} Suppose $\sigma(J_0)$ has closed gaps at $\{y_j\}_{j=1}^\ell
\subset\sigma(J_0)$. Then \eqref{9.11} holds with $\alpha=\f12$ if and only if
\begin{equation} \lb{10.6}
\int_{\sigma_\ess (J_0)} \dist (x,\bbR\setminus\sigma (J_0))^{1/2}
\prod_{j=1}^\ell \abs{x-y_j}^2 \abs{\log (\omega(x))}\, dx <\infty
\end{equation}
\end{lemma}

Plugging this into our proof of Theorem~\ref{T1.4} immediately yields

\begin{theorem}\lb{T10.3} Let $J_0$ be a two-sided periodic Jacobi matrix with closed
gaps at $\{y_j\}_{j=1}^\ell \subset \sigma(J_0)$ and let $J$ be a Jacobi matrix. Then
\begin{equation} \lb{10.7}
\tr((\Delta_{J_0}(J)-S^p-S^{-p})^2) <\infty
\end{equation}
if and only if
\begin{SL}
\item[{\rm{(i)}}] \eqref{1.47} holds.
\item[{\rm{(ii)}}] \eqref{1.49} holds with $\f12$ replaced by $\f32$.
\item[{\rm{(iii)}}] \eqref{10.6} holds.
\end{SL}
\end{theorem}

\begin{example} \lb{E10.4} Take $J_0$ to be the two-sided free Jacobi matrix but think of
it as period $2$. Then
\[
\Delta_J(x)=x^2-2
\]
and a direct calculation of $J^2 -S^2-S^{-2}$ shows that \eqref{10.7} is equivalent to
the three conditions
\begin{align}
\sum_n (a_n^2 + b_n^2 + a_{n+1}^2 -2)^2 &<\infty \lb{10.8} \\
\sum_n (a_{n+1} (b_n + b_{n+1}))^2 &< \infty \lb{10.9} \\
\sum_n (a_n a_{n+1}-1)^2 &< \infty \lb{10.10}
\end{align}
If $b_n=0$ and
\begin{equation} \lb{10.11}
a_n=1 + (-1)^n (n+1)^{-\beta}
\end{equation}
then \eqref{10.8}--\eqref{10.10} hold if and only if $\beta >\f14$ while, of course,
\eqref{1.8} requires $\beta >\f12$. This is one of many known extensions of the
$(J_0$-free) Killip--Simon theorem (see, e.g., Laptev et al.\ \cite{LNS2003}, Kupin \cite{Kup2003},
and Nazarov et al.\ \cite{NPVY}).  Some of these results have MOPRL analogs which, via the
magic formula, lead to variants of Theorems~\ref{T1.4} and \ref{T10.3}.
\end{example}

\section{Eigenvalue Bounds for MOPRL} \lb{s11}

There are Birman--Schwinger kernels for MOPRL and it should be possible to extend the proofs
of most bounds on the number of eigenvalues outside $[-2,2]$ or on moments of $\abs{E_j}-2$
from the scalar to matrix case with optimal constants. But if one is willing to settle for
less than optimal constants (but still not awful constants), there is a simple method
to go from the scalar to matrix case. It depends on the following:

\begin{theorem}\lb{T11.1} Let $\calJ$ be an $\ell\times\ell$ block Jacobi matrix in the
Nevai class with Jacobi parameters $\{A_n, B_n\}_{n=1}^\infty$. Let $E_j^\pm (\calJ)$ denote
its eigenvalues counting multiplicity outside $[-2,2]$, that is, $E_1^+\geq E_2^+
\geq \cdots > 2 > - 2 > \cdots\geq E_2^- \geq E_1^-$. Let $J_\pm$ be the
scalar Jacobi matrix with $a_n\equiv 1$ and
\begin{equation} \lb{11.1}
b_n^\pm =\pm \|B_n\| \pm \|A_{n-1}-1\| \pm \|A_n-1\|
\end{equation}
and let $J_\pm^{(\ell)}$ be an $\ell$-fold direct sum of $J_\pm$. Then
\begin{equation} \lb{11.2}
\abs{E_j^\pm (\calJ)} \leq \abs{E_j^\pm (J_\pm^{(\ell)})}
\end{equation}
\end{theorem}

\begin{proof} The matrix analog of the observation of Hundertmark--Simon \cite{HunS}
extended to $2\ell\times 2\ell$ matrices (with $\ell\times\ell$ blocks) says that
\begin{equation} \lb{11.2x}
\begin{pmatrix}
-\abs{A_n -\bdone} & \bdone \\
\bdone & -\abs{A_n-\bdone}
\end{pmatrix}
\leq
\begin{pmatrix}
0 & A_n \\
A_n^\dagger & 0
\end{pmatrix}
\leq
\begin{pmatrix}
\abs{A_n -\bdone} & \bdone \\
\bdone & \abs{A_n -\bdone}
\end{pmatrix}
\end{equation}
since
\begin{equation} \lb{11.3}
\left\| \begin{pmatrix}
0 & C \\
C^\dagger & 0
\end{pmatrix} \right\|^2 =
\left\| \begin{pmatrix}
C^\dagger C & 0 \\
0 & CC^\dagger
\end{pmatrix} \right\| = \|C\|^2
\end{equation}
Thus writing $J_\pm^{(\ell)}$ as $\ell\times\ell$ blocks with each block a multiple of $\bdone$,
\begin{equation} \lb{11.4}
J_-^{(\ell)} \leq \calJ \leq J_+^{(\ell)}
\end{equation}
from which \eqref{11.2} is immediate.
\end{proof}

\begin{corollary}\lb{C11.2} For any block Jacobi matrix, $\calJ$, in Nevai class,
\begin{equation} \lb{11.5}
\sum_{j,\pm}\, (E_j^\pm (\calJ)^2-4)^{1/2} \leq 2\ell \sum_n\, \|B_n\|
+ 4\ell \sum_n\, \|A_n-1\|
\end{equation}
\end{corollary}

\begin{remark} In particular, this implies if the RHS of \eqref{11.5} is finite, so is the LHS.
\end{remark}

\begin{proof} Hundertmark--Simon \cite{HunS} proved
\begin{equation} \lb{11.6}
\sum_j\, (E_j^\pm (J_\pm)^2 -4)^{1/2} \leq \sum_n b_n^\pm
\end{equation}
from which \eqref{11.5} follows by \eqref{11.2}.
\end{proof}

\section{The Analog of Nevai's Conjecture} \lb{s12}

\begin{proof}[Proof of Theorem~\ref{T1.5}] By \eqref{1.53}, $J-J_0$ is trace class. Thus
$J^\ell -J_0^\ell =\sum_{k=0}^{\ell-1} J^k (J-J_0) J^{\ell-1-k}$ is trace class, so
$\Delta_{J_0}(J)-\Delta_{J_0}(J_0)=\Delta_{J_0}(J) -(S^p + S^{-p})$ is trace class.

It follows that if \eqref{1.53} holds and $\Delta(J)$ has matrix Jacobi parameters
$\{A_n,B_n\}_{n=1}^\infty$ that
\begin{equation} \lb{12.1}
\sum_{n=1}^\infty\, \|B_n\| + \sum_{n=1}^\infty \, \|A_n-1\|<\infty
\end{equation}

By Corollary~\ref{C11.2}, the eigenvalues $\Delta(J)$ obey
\begin{equation} \lb{12.2}
\sum_{j=1}^\infty\, (\abs{E_j^\pm}-2)^{1/2} <\infty
\end{equation}
\eqref{12.1} also implies
\begin{equation} \lb{12.3}
\sum_{n=1}^\infty\, \abs{\log (\det\abs{A_n})} <\infty
\end{equation}
We can apply Theorem~\ref{T8.2} and conclude that $Z(J)$ is finite, that is,
\begin{equation} \lb{12.4}
\int (4-E^2)^{-1/2} \log (\det(W(E)))\, dE >-\infty
\end{equation}
By Corollary~\ref{C9.2&3}, we obtain \eqref{1.50}.
\end{proof}

\section{Perturbations of Periodic OPUC} \lb{s13}

In this final section, we want to present the translations of our results to OPUC.
Since the magic formula maps periodic OPUC to MOPRL, the changes needed in the proofs
will be minor, although for the analog of Theorem~\ref{T1.4}, there is one significant
change. It is interesting to note the sequence of mappings for the OPUC periodic
Rakhmanov's theorem. We map OPUC to MOPRL using the magic formula and them map that
to MOPUC using the Szeg\H{o} map.

The OPUC version of Theorem~\ref{T1.1} is already in Last--Simon \cite{LastS}. As
for Theorem~\ref{T1.2}:

\begin{theorem}\lb{T13.1} Let $\calC_0$ be a two-sided periodic CMV matrix. Let $\calC$
be an ordinary CMV matrix with Verblunsky coefficients $\{\alpha_n\}_{n=0}^\infty$. Suppose
\begin{equation} \lb{13.1}
\Sigma_\ac (\calC) = \sigma_\ess (\calC_0)
\end{equation}
Then, as $m\to\infty$,
\begin{equation} \lb{13.2}
d_m (\alpha,\calT_{\calC_0}) \to 0
\end{equation}
\end{theorem}

\begin{proof} Let $\calC_r$ be a right limit of $\calC$. By Theorem~\ref{T6.1},
$\Delta (\calC_r)=S^p + S^{-p}$. Thus, by Theorem~\ref{T3.1}, $\calC_r\in\calT_{\calC_0}$.
\end{proof}

As for the analog of Theorem~\ref{T1.3}, if we drop discussion of $\ell^2$ convergence,
it holds, similar to Theorem~\ref{T10.1}.

\begin{theorem}\lb{T13.2} Let $\calC_0$ be a two-sided periodic CMV matrix with
Verblunsky coefficients $\{\alpha_j^{(0)}\}_{j=-\infty}^\infty$, and $\calC$ a one-sided
CMV matrix with Verblunsky coefficients $\{\alpha_j\}_{j=0}^\infty$ and spectral measure
$d\mu$. Suppose that $\sigma_\ess(\calC)=\sigma(\calC_0)$ and
\begin{equation} \lb{13.3}
\sum_{m=1}^N \dist (E_m,\sigma_\ess (\calC))^{1/2} <\infty
\end{equation}
where $\{E_m\}_{m=1}^N$ is a labelling of the eigenvalues of $\calC$ outside $\sigma_\ess (\calC_0)$. Then
\begin{equation} \lb{13.4}
-\int_{\sigma_\ess (\calC_0)} \log \biggl( \f{d\mu_\ac}{d\theta}\biggr)
\dist (\theta, \bbR\setminus\sigma_\ess (\calC_0))^{-1/2}\, \f{d\theta}{2\pi} <\infty
\end{equation}
implies
\begin{equation}\lb{13.5}
\lim_{N\to\infty}\biggl(\, \sum_{j=1}^{pN} \log \biggl( \f{\rho_j}{\rho_j^{(0)}}\biggr)\biggr)
\end{equation}
exists and lies in $(-\infty,\infty)$. Conversely, \eqref{1.50} holds so long as
\begin{equation} \lb{13.6}
\limsup_{N\to\infty} \biggl(\, \sum_{j=1}^{N} \log \biggl( \f{\rho_j}{\rho_j^{(0)}}\biggr)\biggr) >-\infty
\end{equation}
and then the limit in \eqref{13.5} exists and lies in $(-\infty,\infty)$.
\end{theorem}

\begin{remarks} 1. We have not stated that the $\alpha_n$ have a limit in $\calT_{\calC_0}$.
We suspect that the methods of \cite{PY} extend to the OPUC but have not checked this and they
do not explicitly mention it.

\smallskip
2. Of course, if this theorem is applicable and $\calC_0$ obeys the conditions of
Theorem~\ref{T13.3} below, then we have a result on $\ell^2$ convergence to $\calT_{\calC_0}$.

\smallskip
3. One can replace $\rho_j^{(0)}$ by the logarithmic capacity of $\sigma(\calC_0)$.
\end{remarks}

\begin{proof} At open gaps, $\Delta'(e^{i\theta})\neq 0$, so \eqref{13.3} is equivalent to
\begin{equation} \lb{13.7}
\sum_{\substack{ E\notin [-2,2] \\ E\in\sigma(\Delta(\calC))}}\, (\abs{E}-2)^{1/2} <\infty
\end{equation}
Moreover, by \eqref{3.3}--\eqref{3.6}, we have that
\begin{equation} \lb{13.8}
\log(\det\abs{A_n}) = \log (\det(A_n)) = \log \biggl[\,\prod_{j=(n-1)p+1}^{np}
\prod_{k=j}^{j+p^{-1}} \f{\rho_k}{\rho_k^{(0)}}\biggr]
\end{equation}
Now just follow the proof of Theorem~\ref{T10.1}.
\end{proof}

In carrying over Theorem~\ref{T1.4} to OPUC, one runs into a serious roadblock:
van Moerbeke's theorem \cite{vMoer} that the Hamiltonian flows generated by the coefficients
of the $t_j$ (given by \eqref{9.17}) are independent is not known for OPUC. Instead, we
use a weaker result of Simon \cite[Section~11.10]{OPUC2} that proves the derivatives of
coefficients of $\Delta$ span the normal bundle for a dense open set that is `most' of
the points with all gaps open. We will call the isospectral tori in this dense open set
the generic independence tori. Then by mimicking the arguments in Section~\ref{s9}, we
obtain

\begin{theorem}\lb{T13.3} Let $\calC_0$ be a two-sided CMV matrix in a generic independence
torus with Verblunsky coefficients $\{\alpha_j^{(0)}\}_{j=-\infty}^\infty$. Let $\calC$ be
a CMV matrix with Verblunsky coefficients $\{\alpha_j\}_{j=0}^\infty$. Then
\begin{equation} \lb{13.9}
\sum_{m=0}^\infty d_m (\alpha, \calT_{\calC_0})^2 <\infty
\end{equation}
if and only if
\begin{SL}
\item[{\rm{(i)}}] $\sigma_\ess (\calC)=\sigma(\calC_0)$
\item[{\rm{(ii)}}] For the eigenvalues $\{E_j\}_{j=1}^N$ not in $\sigma(\calC_0)$,
\begin{equation} \lb{13.10}
\sum_{j=1}^N \dist (E_j, \sigma (\calC_0))^{3/2} <\infty
\end{equation}
\item[{\rm{(iii)}}] If $\mu$ is the spectral measure for $\calC$, then
\begin{equation} \lb{13.11}
-\int_{\sigma(\calC_0)} \log \biggl( \f{d\mu_\ac}{d\theta}\biggr)
\dist(\theta,\partial\bbD\setminus \sigma_\ess (\calC_0))^{1/2}\, \f{d\theta}{2\pi} <\infty
\end{equation}
\end{SL}
\end{theorem}

Our last result is a periodic OPUC version of Nevai's conjecture.

\begin{theorem}\lb{T13.4} Let $\calC_0$ be a two-sided $p$-periodic CMV matrix
and let $\calC$ be a CMV matrix with Verblunsky coefficients $\{\alpha_j\}_{j=0}^\infty$. Then
\begin{equation} \lb{13.12}
\sum_{m=0}^\infty d_m (\alpha, \calT_{\calC_0}) <\infty
\end{equation}
implies \eqref{13.4}
\end{theorem}

All the above results assume the period $p$ is even. However, by using sieving
(see Example~1.6.14 of \cite{OPUC1}) to map period $p$ to period $2p$, one can extend
Theorems~\ref{T13.1} and \ref{T13.2} to $p$ odd. In particular, we obtain the $p=1$
results of \cite{BHLL,ABMV,BCL} as very special cases of ours.

\bigskip

\end{document}